\documentclass[a4paper,twoside,12pt]{report}
\usepackage{amssymb}
\usepackage{amsmath}
\usepackage{latexsym}
\input xy
\xyoption{all}

\newtheorem{teo}{Theorem}[section]
\newtheorem{prop}[teo]{Proposition}
\newtheorem{lem}[teo]{Lemma}
\newtheorem{cor}[teo]{Corollary}

\newtheorem{defi}[teo]{Definition}
\newtheorem{fac}[teo]{Fact}
\newtheorem{rem}[teo]{Remark}
\newtheorem{ej}[teo]{Example}
\newtheorem{nota}[teo]{Notation}
\newtheorem{notadefi}[teo]{Notation and Definition}

\def\dnfo{\,\raise.2em\hbox{$\,\mathrel|\kern-.9em\lower.35em\hbox{$\smile$}$}}
\def\dfo{\;\raise.2em\hbox{$\mathrel|\kern-.9em\lower.35em\hbox{$\smile$}\kern-.7em\hbox{\char'57}$}\;}

\newcommand{\td}{tr.dg}
\newcommand{\s}{SU}
\newcommand{\f}{FR}
\newcommand{\si}{\sigma}
\newcommand{\na}{\mathbb N}
\newcommand{\gr}{\mathbb G}

\newcommand{\aff}{\mathbb A}

\begin{document}
\pagenumbering{roman}

\begin{titlepage}

{\center {\Large \bf Universit\'e Paris 7 - Denis Diderot

\vspace{0.5cm}

U.F.R. de Math\'ematiques}

\vspace{1cm}

{\large  Th\`ese de Doctorat
\vspace{0.5cm}
 
Sp\'ecialit\'e:\ Logique et Fondements de l'Informatique}

\vspace{3cm}

{\Large \bf RONALD BUSTAMANTE MEDINA       }

\vspace{3cm}

  {\huge  \bf Th\'eorie des mod\`eles des corps diff\'erentiellement clos avec un automorphisme g\'en\'erique}

\vspace{1.5cm}

 { Soutenue le 10 novembre 2005

\vspace{1cm}

\begin{tabular}{ll}

   Directeur: & Zo\'e Chatzidakis   \\    

                 &          \\

   Rapporteurs:
 %&    \\

                & Anand PILLAY \\

                     &Thomas SCANLON         \\
   Jury:        %&          \\
               &  Elisabeth BOUSCAREN\\
               &  Zo\'e CHATZIDAKIS        \\
 &  Fran\c coise DELON        \\
&  Angus MACINTYRE        \\
&  Dugald MACPHERSON        \\
 & Fran\c coise POINT        \\

\end{tabular}

     }

}

\end{titlepage}
\setcounter{secnumdepth}{-1}
\chapter{Remerciements}

Depuis le d\'ebut de mes \'etudes en France Zo\'e Chatzidakis a \'et\'e pr\'esente. Sans son
soutien, diligence, patience et diposition cette th\`ese n'aurait jamais vu le jour.
Son encouragement, son enthousiasme et ses id\'ees m'ont toujours \'et\'e tr\`es utiles et je lui en suis
tr\`es reconnaissant.\\

Merci \`a Thomas Scanlon et Anand Pillay d'avoir accept\'e d'\^etre les rapporteurs de ma th\`ese, ainsi 
qu'\`a Elisabeth Bouscaren, Fran\c coise Delon, Angus Macintyre, Dugald Macpherson et Fran\c coise Ponitde m'honorer
de leur pr\'esence dans le jury.\\

Le D.E.A. a \'et\'e une exp\'erience enrichissante gr\^ace aux enseignements de Ren\'e Cori, 
Max Dickmann et Richard Lassaigne.\\

Je tiens \`a remercier Mich\`ele Wasse, Virginie Kuntzmann et Khadija Bayoud pour leur accueil, aide et disponibilt\'e.\\
  
Merci aux amis qui m'ont accueilli et aid\'e d\`es mon arriv\'ee en France: Andrea, Gianluca, Eugenio, Miriam, Xos\'e,
 Oldemar. 
Je veux remercier \'egalement tous mes coll\`egues du cinqui\`eme \'etage.\\

Merci \`a Gianmarco, Mouna, Emanuela et Zo\'e pour les repas, les th\'es, les caf\'es, la musique, le foot, le cinema 
et les journ\'ees gourmandes. Merci \'egalement \`a tous les amis du septi\`eme \'etage (grazie, gracias).\\

En Costa Rica quisiera agradecer a la gente de la Escuela de Matem\'aticas, en particular a Santiago Cambronero, Asdr\'ubal
Duarte y Bernardo Montero. Gracias a la Oficina de Asuntos Internacionales de la Universidad de Costa Rica, en particular a
F\'atima Acosta y Yamileth Damazio por toda su ayuda. Agradezco al Consejo Nacional para Investigaciones
Cient\'{\i}ficas y Tecnol\'ogicas por su apoyo.\\

Much\'{\i}simas gracias a mi familia por su apoyo y ayuda a lo lardo de estos a\~nos, muy especialmente a mi madre
Mabel Medina Arana y a mi padre Ronald Bustamante Valle. Gracias a mis hermanos Mari, Gabi, Dani, Nati y Rodolfo. 
Muchas gracias a Jen por todo.\\
Gracias a mi abuelo Don Paco, quien me inici\'o en las matem\'aticas.\\
Gracias infinitas al infinito Sib\'u-Sur\'a del mundo m\'as arriba y el mundo m\'as abajo.

\cleardoublepage
\tableofcontents
\setcounter{secnumdepth}{-1}

\chapter{Notation}

$K[{\bar X} ]_{\sigma,D} $, ring of difference-differential polynomials over $K$.\\
$(A)_{\sigma,D}  $, smallest difference-differential field containing $A$.\\
$A^{alg} $,  field-theoretic algebraic closure of $A$.\\
$acl(A)  $,  model-theoretic algebraic closure of $A$.\\
$acl(A)_{D}  $, algebraic closure of the smallest differential field containing $A$.
\\
$acl(A)_{\sigma}  $, algebraic closure of the smallest difference field containing $A$.
\\
$acl(A)_{\sigma,D}  $, algebraic closure of the smallest difference-differential field containing $A$.
\\
$E(a)_{D} $, the field generated by $E$ and $\{ D^ja;i  j \in {\mathbb N}\}  $.
\\ 
$E(a)_{\sigma} $, the field generated by $E$ and $\{ \sigma^i(a);i \in {\mathbb Z}\}  $.
\\ 
$E(a)_{\sigma,D} $, the field generated by $E$ and $\{ \sigma^i(D^ja);i \in {\mathbb Z}, j \in {\mathbb N}\}  $.
\\ 
$V^{\sigma}  ,$  image of $V$ by $\sigma$.
\\
$V(K)  $, elements of $V$ with all its coordinates in $K$.
\\
${\cal L}=\{+,-, \cdot, 0,1   \}   $ .
\\
${\cal L}_D={\cal L} \cup \{D   \}  $.
\\
${\cal L}_{\sigma}={\cal L} \cup \{\sigma   \}  $.
\\
${\cal L}_{\sigma,D}={\cal L} \cup \{\sigma,D   \}  $.
\\
$ACF$, the theory of algebraically closed fields.
\\
$qfDiag(E)$, quantifier-free diagram of $E$: the set of quantifier-free closed ${\cal L}_{\sigma,D}(E)  $-formulas satisfied by some ${\cal L}_{\sigma,D}  $-structure containing $E$.
\\
$S_n(A)  $, set of complete $n$-types with parameters in $A$.
\\
$\tau_m(V)$, $m$-th (differential) prolongation of the variety $V$.
\\
$J^m(V)_a$, $m$-th jet space of the variaty $V$ at $a$.
\\
$\Phi_m(V)$, $m$-the $(\si,D)$-prolongation of the variety $V$.
\\
${\mathcal A}_mV_a$, $m$-th arc space of the variety $V$ at $a$.
\\
$T_{\si,D}(V)_a$, $m$-th $(\si,D)$-tangent space of the variety $V$ at $a$.
\\
$A^{\sharp}$, Manin kernel of the Abelian variety $A$.

\chapter{Avant-propos}

\pagenumbering{arabic}
%\noindent%\begin{large}{\bf  Avant-propos}\end{large}
L'\'etude mod\`ele-th\'eorique des corps enrichis, c'est \`a dire, munis
d'un ou de plusieurs op\'erateurs (d\'erivation, automorphisme,
$\lambda$-fonctions ou d\'erivations de Hasse pour les corps
s\'eparablement clos), ou bien d'une valuation, a connu ces derni\`eres
ann\'ees un essor spectaculaire, d\^u en grande partie aux applications
de la th\'eorie des mod\`eles de ces corps, associ\'ee \`a la {\it
  Stabilit\'e G\'eom\'etrique},  pour r\'esoudre des questions
en G\'eom\'etrie Diophantienne. Mentionnons par exemple les travaux de
Hrushovski 
sur la Conjecture de Mordell-Lang pour les corps de fonctions (utilisant
les corps diff\'erentiellement clos, et les corps s\'eparablement clos),
ceux encore de Hrushovski sur la conjecture de Manin-Mumford et qui
donnent des bornes effectives (utilisant les corps de diff\'erence), et
enfin ceux de Scanlon sur la distance $p$-adique des points de torsion
(corps valu\'es de diff\'erence) et sur la conjecture de Denis sur les
modules de Drinfeld (corps de diff\'erence). Un des ingr\'edients
communs \`a la d\'emonstration de  ces r\'esultats, est le fait que
certains ensembles d\'efinissables contenant ceux dans lesquels nous
sommes int\'eress\'es (par exemple, le sous-groupe de torsion) sont
mono-bas\'es. En effet, un r\'esultat d\'ej\`a relativement ancien, nous
dit que ces ensembles se comportent bien :

\smallskip\noindent
{\bf Th\'eor\`eme} (Hrushovski-Pillay \cite{phgr}). Soit $G$ un groupe
stable mono-bas\'e, d\'efini sur $\emptyset$. Alors tout sous-ensemble
d\'efinissable de $G^n$ est une combinaison Bool\'eenne de cosets de
sous-groupes d\'efinissables de $G^n$. De plus, ces sous-groupes sont
d\'efinis sur $acl(\emptyset)$.

\smallskip
Ce r\'esultat dit deux choses : tout d'abord, $G$ ne contiendra pas de
famille infinie de sous-groupes d\'efinissables. Ensuite, la description
des ensembles d\'efinissables permettra souvent de d\'eduire que
l'ensemble auquel nous nous int\'eressons, est une {\it union
  finie} de cosets de sous-groupes de $G$. 

\smallskip
Jusqu'au milieu des ann\'ees 90, l'utilisation des outils de la Stabilit\'e
G\'eom\'etrique \'etait r\'eserv\'ee aux structures dont la th\'eorie
est stable : les corps alg\'ebriquement clos, s\'eparablement clos, ou
bien 
diff\'erentiellment clos. Les travaux de Hrushovski sur la conjecture de
Manin-Mumford ont montr\'e que ces techniques pouvaient aussi
s'appliquer dans des cas instables : celui des corps avec automorphisme
g\'en\'erique. On s'est alors aper\c cu que la th\'eorie de ces corps
\'etait supersimple. La simplicit\'e est une propri\'et\'e de certaines
th\'eories qui a \'et\'e remarqu\'ee par Shelah en 1980, mais peu
\'etudi\'ee jusqu'en 1995, quand Kim a montr\'e que les th\'eories
simples avaient des propri\'et\'es tr\`es int\'eressantes, notamment la
sym\'etrie de la d\'eviation (la d\'eviation est une
g\'en\'eralisation mod\`ele-th\'eorique de la notion de d\'ependance,
alg\'ebrique ou lin\'eaire par exemple). Le d\'eveloppement de
l'\'etude des th\'eories simples a montr\'e que en effet, on pouvait
souvent appliquer \`a leurs mod\`eles des techniques provenant de la
stabilit\'e.

\smallskip
Dans cette th\`ese, nous nous int\'eressons au cas des corps {\it de
  caract\'eristique $0$} munis d'une
d\'erivation $D$ et d'un automorphisme $\si$ qui commute avec $D$
(appel\'es corps diff\'erentiels de diff\'erence), et plus
particuli\`erement aux mod\`eles existentiellement clos de cette classe
(appel\'es corps diff\'erentiellement clos avec un automorphisme
g\'en\'erique). Ce sont des corps diff\'erentiels de
diff\'erence tels que tout syst\`eme fini d'\'equations diff\'erentielles
de diff\'erence \`a coefficients dans le corps et qui a une solution
dans une extension, a d\'ej\`a une solution dans le corps.
Hrushovski a montr\'e que les corps diff\'erentiels avec automorphisme
g\'en\'erique forment une classe \'el\'ementaire. On sait par ailleurs
que la th\'eorie d'un tel corps est supersimple, ce qui a plusieurs
  cons\'equences importantes pour notre \'etude. 

\smallskip
Le but de cette th\`ese est d'\'etudier la th\'eorie de ces corps
(not\'ee {\it DCFA}), et de voir dans quelle mesure les r\'esultats qui
ont \'et\'e prouv\'es pour les corps diff\'erentiellement clos  et pour
les corps avec automorphisme 
g\'en\'erique (dont les th\'eories seront not\'ees respectivement {\it
  DCF} et {\it ACFA}) peuvent \^etre g\'en\'eralis\'es \`a ces corps. De
plus, l'existence de deux op\'erateurs donne une structure plus riche,
et l'\'etude ces corps permettra peut-\^etre d'isoler des ph\'enom\`enes
nouveaux. 

\smallskip
Les travaux de cette th\`ese commencent par un chapitre de
pr\'eliminaires, qui pour l'essentiel rappelle des notions connues, mais
montre aussi quelques r\'esultats nouveaux. Nous continuons dans le
chapitre 2 avec la
description de la th\'eorie $DCFA$, et avec quelques r\'esultats assez
g\'en\'eraux. Le chapitre 3 donne une d\'emonstration de la dichotomie
corps/monobas\'e pour les ensembles de rang $1$. Enfin le chapitre 4
\'etudie les groupes d\'efinissables, et plus particuli\`erement les
sous-groupes d\'efinissables de groupes alg\'ebriques
commutatifs. Ci-dessous nous exposons les r\'esultats plus en d\'etail.

\bigskip\noindent
\begin{large}{\bf Chapitre 1}\end{large}\\
\smallskip

Ce chapitre est divis\'e en plusieurs sections. Dans la premi\`ere, nous
rappelons certaines des d\'efinitions et  propri\'et\'es des th\'eories
stables et des 
th\'eories simples que nous utiliserons par la suite, en particulier le
r\'esultat de Kim-Pillay qui caract\'erise la notion d'ind\'ependance
(ou de non-d\'eviation) dans les th\'eories simples. Nous introduisons
aussi le rang $\s$, une g\'en\'eralisation du rang U de Lascar (un rang
est une notion de dimension, mais qui peut prendre des valeurs
ordinales), ainsi que certaines de ses propri\'et\'es.
La section suivante rappelle les notions fondamentales de stabilit\'e
g\'eom\'etrique : orthogonalit\'e, 
  r\'egularit\'e, internalit\'e et analysabilit\'e, mono-bas\'eit\'e et
modularit\'e locale, poids, domination, etc. Nous y prouvons deux
petits lemmes, qui nous seront utiles pour nous r\'eduire au cas de
types r\'eguliers de rang $\omega^i$, pour lesquels nous 
n'avons pas trouv\'e de r\'ef\'erence explicite. Le premier est bien 
connu des sp\'ecialistes, et en plus grande g\'en\'eralit\'e. Le 
deuxi\`eme est connu quand $\alpha=0$, mais \`a notre connaissance
(et apr\`es avoir consult\'e F. Wagner) est nouveau. 
Nous travaillons dans une
th\'eorie supersimple qui \'elimine les imaginaires.\\ 

\smallskip\noindent
{\bf Lemme \ref{prt011}}. {\it Soit $q=tp(a,b)$ un type r\'egulier, et $p=tp(a)$. Alors
$q$ est localement modulaire si et seulement si $q$ est localement
  modulaire.}\\

\smallskip\noindent
{\bf Lemma \ref{prt012}}. {\it Soit $p=tp(a/A)$ un type r\'egulier de rang SU
  $\beta+\omega^\alpha$ et de poids $1$. Alors il existe $b\in acl(Aa)$ tel que
  $\s(a/A)=\omega^\alpha$.}\\

Dans la section 3, nous rappelons les r\'esultats d\'ej\`a connus sur
la th\'eorie des mod\`eles avec automorphisme g\'en\'erique : soit $T$
une th\'eorie dans un langage ${\mathcal L}$, qui est stable et \'elimine les quantificateurs et les
imaginaires. Nous rajoutons au langage un symbole de fonction $\si$ et
consid\'erons la th\'eorie $T_0$ des mod\`eles $M$ de $T$ dans lesquels
$\si$ est interpr\'et\'e par un automorphisme. On sait, par des
r\'esultats de Chatzidakis-Pillay (\cite{zopi}) que si la
mod\`ele-compagne (not\'ee $T_A$) de cette th\'eorie existe, alors : on
peut d\'ecrire facilement les compl\'etions de $T_A$ et leurs types ;
toutes les compl\'etions de $T_A$ sont simples, et supersimples si la
th\'eorie $T$ est superstable ; on a une bonne description de la
cl\^oture alg\'ebrique et de la relation d'ind\'ependance.

Dans les sections 4 et 5 nous nous int\'eressons aux groupes
$\omega$-stables et aux groupes d\'efinissables dans les th\'eories
simples. 
Nous rappelons quelques notions de base, comme les stabilisateurs
et les types g\'en\'eriques, ainsi que leurs propri\'et\'es.

Le reste du chapitre est plus alg\'ebrique. Nous rappelons dans la
section 6 les r\'esultats bien connus sur les corps diff\'erentiels et
les corps diff\'erentiellement clos, la topologie de Kolchin, les
vari\'et\'es et id\'eaux diff\'erentiels, 
etc. . Pour pouvoir donner une 
axiomatisation de la th\'eorie $DCFA$, nous rappelons les notions de prolongations de
vari\'et\'es alg\'ebriques, et introduisons les notions de vari\'et\'es
en forme normale et de points $(m,D)$-g\'en\'eriques. Nous montrons que
les points $(m,D)$-g\'en\'eriques d'une vari\'et\'e en forme normale ont tous le
m\^eme type dans $DCF$, et que la question de savoir si les points
$(m,D)$-g\'en\'eriques d'une vari\'et\'e en forme normale se projettent sur les points
$(m,D)$-g\'en\'eriques d'une autre vari\'et\'e en forme normale, se
r\'eduit \`a une question de dominance de morphisme de vari\'et\'es
alg\'ebriques. Ces notions nous permettent donc d'\'eviter de r\'epondre
aux deux questions suivantes :

\smallskip
-- Quand des polynomes diff\'erentiels engendrent-ils un
   id\'eal diff\'erentiel premier?

-- Si $I\subset K[X,Y]_D$ et $J\subset K[X]_D$ sont des id\'eaux diff\'erentiels
   premiers donn\'es par des syst\`emes de g\'en\'erateurs, quand avons
   nous $I\cap K[X,Y]_D=J$?

\smallskip
Ces questions sont certainement r\'esolues dans la litt\'erature, mais
nous n'avons pas trouv\'e de bonne r\'ef\'erence, et de plus notre
approche est plus dans l'esprit des axiomatisations g\'eom\'etriques. 

\smallskip
Dans la section 7, nous rappelons les d\'efinitions de base sur les
corps de diff\'erence, la $\si$-topologie, et la th\'eorie $ACFA$ des
corps alg\'ebriquement clos avec automorphisme g\'en\'erique. 

\bigskip\noindent \begin{large}{\bf Chapitre 2}\end{large} \\

\smallskip
 Le deuxi\`eme chapitre commence par une introduction g\'en\'erale aux corps diff\'erentiels de diff\'erence. 
 Un r\'esultat tr\`es important est la Noeth\'erianit\'e de la
$(\si,D)$-topologie. Ensuite nous utilisons les r\'esultat du chapitre 1
pour donner un sch\'ema d'axiomes pour la
th\'eorie  {\it DCFA}. (Les r\'esultats de Hrushovski montrant que les corps
diff\'erentiellement clos avec automorphisme g\'en\'erique forment une
classe \'el\'ementaire ne sont pas publi\'es).   \\  

%{\bf Th\'eor\`eme \ref{DCFA49}}  %\textup{(}Hrushovski\textup{)} {\it La mod\`ele-compagne de la th\'eorie 
%des corps diff\'erentiels de diff\'erence existe.  Elle est not\'ee {\it DCFA} et ses axiomes sont : 
% $(K,D,\sigma)$ est un mod\`ele de {\it DCFA} si \begin{enumerate} \item $(K,D)$ est un corps diff\'erentiellement clos.
% \item $\sigma$  est un automorphisme de $(K,D)$. \item If $U,V,W$ sont des  vari\'et\'es telles que: \begin{enumerate}
% \item $U \subset V \times V^{\sigma}  $ se projette g\'en\'eriquement sur $V$ et  $ V^{\sigma}  $. \item $W \subset 
%\tau_1 (U) $ se projette g\'en\'eriquement sur $U$. \item $\pi_1(W)^{\sigma}=\pi_2(W)$ \textup{(}on identifie 
%$\tau_1(V \times V^{\sigma})$  avec $\tau_1(V)\times \tau_1(V)^{\si}$\textup{)}.  \item Un $(1,D) $-g\'en\'erique 
%of $W $ se projette   sur un $(1,D)   $-g\'en\'erique de $\pi_1(W)$ et sur un (1,D)-g\'en\'erique de $\pi_2(W) $. 
%\end{enumerate} Alors il existe un  uplet $a \in V(K) $, tel que  $(a,Da,\si(a),\si(Da)) \in W$.  \end{enumerate} }     
Comme on l'a vu dans le chapitre 1, l'existence d'une mod\`ele-compagne de la th\'eorie 
des corps diff\'erentiels avec un 
automorphisme a quelques
cons\'equences imm\'ediates: 

\begin{enumerate}
\item Les compl\'etions de $DCFA$ sont obtenues en d\'ecrivant le type
  d'isomorphisme du corps de diff\'erence $({\mathbb Q}^{alg},\si)$.
\item Soit $K\models DCFA$, et $A\subset K$. La cl\^oture alg\'ebrique
  de $A$ est le plus petit sous-corps diff\'erentiel alg\'ebriquement
  clos de $K$ qui contient $A$ et est clos par $\si$ et $\si^{-1}$. 
\item L'ind\'ependance est d\'ecrite de la fa\c con suivante : $A$ et
  $B$ sont ind\'ependants au dessus de $C$ (not\'e : $A \dnfo_CB$) si et
  seulement si les corps $acl(CA)$ et $acl(CB)$ sont lin\'eairement
  disjoints au-dessus de $acl(C)$.
\item Toute compl\'etion de $DCFA$ est supersimple et satisfait au Th\'eor\`eme d'ind\'ependance au dessus d'un corps de diff\'erence diff\'erentiellement clos. 
\end{enumerate}

\smallskip
Comme dans $ACFA$, nous montrons ensuite que le 
Th\'eor\`eme d'ind\'ependance est vrai au dessus d'un ensemble
alg\'ebriquement clos :\\ 
\small\noindent 

{\bf Th\'eor\`eme \ref{DCFA426}} 
 {\it Si  $({\mathcal U},\si,D)$ est un mod\`ele  satur\'e  de {\it DCFA}, $E$ un sous-ensemble alg\'ebriquement 
clos de $\mathcal U$,  et  ${\bar a},  {\bar b}, {\bar c_1}, {\bar c_2} $ sont des uplets de $\mathcal U$ tels
 que :  \begin{enumerate} \item $tp(  {\bar c_1}/E  )=tp(  {\bar c_2}/E  )    $. 
 \item ${\bar a}  \dnfo_E    {\bar c_1}$, ${\bar a}  \dnfo_E    {\bar b}$ et   ${\bar b}  \dnfo_E    {\bar c_2}$.
 \end{enumerate} 
 Alors il existe $ {\bar c}  $ realisant
 $tp(  {\bar c_1}/E \cup   {\bar a}  ) \cup   tp(  {\bar c_2}/E \cup   {\bar b}  )  $  tel que
 $ {\bar c}  \dnfo_E    ({\bar a} ,  {\bar b})  $. }\\

\smallskip\noindent 
Ce th\'eor\`eme nous permet de montrer que toute compl\'etion de {\it DCFA}
 \'elimine les imaginaires, voir \ref{DCFA430}. Il nous permet aussi de montrer que le corps diff\'erentiel 
$Fix(\si)=\{x\in {\mathcal U} : \si(x)=x\}$ est stablement plong\'e, c'est-\`a-dire que si $S\subset (Fix \si)^n$ est 
${\mathcal L}_{\si,D}$-d\'efinissable (avec param\`etres dans $\mathcal U$), alors il est d\'efinissable avec
param\`etres dans $Fix \si$. De plus, $S$ est d\'efinissable dans le
langage des corps diff\'erentiels.

\smallskip Dans la section 2.3 on \'etudie le corps de constantes ${\mathcal C}=\{x:Dx=0\}$ et le 
corps fixe $Fix(\si)$ d'un mod\`ele de {\it DCFA}: \'etant donn\'e un mod\`ele $(K, \si, D)$  
de {\it DCFA}, on montre que $({\mathcal C},\si)$ est un mod\`ele de $ACFA$ (\ref{DCFA52}), 
mais qu'il n'est pas stablement plong\'e : il existe des sous-ensembles d\'efinissables de $\mathcal C$ qui 
ne sont pas d\'efinis sur $\mathcal C$. Pour $Fix \si$ on montre aussi  la chose suivante:\\ 

\smallskip\noindent 

{\bf Th\'eor\`eme \ref{DCFA55}}  {\it $((Fix \sigma)^{alg},D   )$ est
un mod\`ele de {\it DCF }. }\\

  Ce r\'esultat nous permet de d\'ecrire les structures de la forme 
$(F,D)$ o\`u $F$ est le corps fix\'e d'un mod\`ele de {\it DCFA}. Ces r\'esultats ont \'et\'e obtenus 
ind\'ependamment (et dans un cadre plus g\'en\'eral) par Pillay et Polkowska, voir \cite{pildom}.

\smallskip\noindent

 La section suivante est d\'edi\'ee \`a l'\'etude du rang $\s$ pour les compl\'etions de $DCFA$;
  entre autres r\'esultats on montre que le rang $\s$ d'un g\'en\'erique d'un mod\`ele de {\it DCFA} est $\omega^2$.
 Nous donnons aussi des bornes sur le rang $\s$ d'un \'el\'ement, et donnons quelques exemples. 

 \smallskip\noindent

Dans 2.5, nous isolons des conditions pour qu'un type dans {\it DCFA} soit stable stablement plong\'e. 

 Les sous-ensembles d\'efinissables qui sont  stables stablement plong\'es sont des ensembles qui sont stables
 pour la structure induite. On peut donc leur appliquer tous les r\'esultats valides dans les th\'eories stables. 
Les conditions obtenues nous serviront pour l'\'etude des groupes ab\'eliens dans le chapitre 4. Par exemple,
 nous montrons que dans plusieurs cas, un type sera stable stablement plong\'e si son r\'eduit au langage des corps 
de diff\'erence est stable stablement plong\'e pour la th\'eorie {\it  ACFA}, voir \ref{st3} et \ref{st4}.    
 \smallskip\noindent

 Dans la derni\`ere section de ce  chapitre nous donnons  un exemple d'un ensemble de rang $\s$ 1 mais de
 dimension infinie. Soit $(K,\si,D)$ un mod\`ele de {\it DCFA}, soit $A=\{x \in K: \si(x)=x^2+1 \land \ Dx\neq 0\}$.
 Si $E$ est un sous-corps diff\'erentiel de diff\'erence de $K$, et $a\in A$ n'est pas alg\'ebrique sur $E$, 
alors $a$ est diff\'erentiellement transcendant sur $E$. Cela entraine que l'ensemble $A$ est fortement minimal 
et nous montrons aussi qu'il est stablement plong\'e.\\

   \bigskip\noindent \begin{large}{\bf Chapitre 3}\end{large} \\

 \medskip  

Le chapitre 3 est consacr\'e \`a la d\'emonstration de la dichotomie de Zilber dans un mod\`ele de {\it DCFA}.
 On veut montrer que si $({\mathcal U},\si,D)$ est un mod\`ele de $\it DCFA$, $K$ un sous-corps diff\'erentiel
 de diff\'erence alg\'ebriquement clos   de ${\mathcal U}$ et  $a \in {\mathcal U}$ tel que $\s(a/K)=1$, 
alors ou bien $tp(a/K)$ est monobas\'e, ou bien $tp(a/K)$ est non-orthogonal \`a ${\mathcal C} \cap Fix \si$. 
 La preuve de ce r\'esultat est faite en deux parties, en consid\'erant
d'abord le cas de {\it dimension finie}, c'est \`a dire quand
$\td(acl(Ka)/K)<\infty$, puis le cas g\'en\'eral. Dans le cas de
dimension finie, nous utilisons  les techniques des espaces de jets
introduites par Pillay et Ziegler dans
(\cite{jets}) et montrons :\\ 

 \smallskip\noindent

{\bf Th\'eor\`eme \ref{J312}} {\it Soit $({\mathcal U},\si,D) $ un mod\`ele satur\'e de {\it  DCFA} et
  soit  $tp(a/K) $ un type de dimension finie. Soit $b$ tel que
 $b=Cb(qftp(a/K,b))  $. Alors $tp(b/acl(K,a))  $ est presque interne \`a $ Fix \sigma \cap {\mathcal C}.$ }\\ 

 \smallskip\noindent A partir de ce r\'esultat nous pouvons d\'eduire une dichotomie partielle :\\

  \smallskip\noindent {\bf Corollaire \ref{J313}} {\it Si $p=tp(a/K) $ a $\s$-rang 1 et est de dimension finie,
 alors ou bien $p$   est monobas\'e, ou bien il est non-orthogonal \`a  $  Fix\si  \cap{\mathcal C}    $. }\\  
  \smallskip  %si $(K,\si,D)$ est un mod\`ele de $\it DCFA$, $L$ un souscorps diff\'erentiel de diff\'erences de $K$ et 
%$a \in K$ tel que $\s(a/L)=1$ et $a$ est de dimension finie sur $L$, %alors soit $tp(a/L)$ est 1-bas\'e, soit $tp(a/L)$ est non-orthogonal a ${\mathcal C} \cap Fix \si$. 

Dans la deuxi\`eme  partie de la preuve nous utilisons les techniques des espaces d'arcs
introduites par Moosa, Pillay et Scanlon (\cite{arcs}) pour montrer : \\ 

\smallskip\noindent 
{\bf Th\'eor\`eme \ref{arc17}} {\it Soit $p$ un type r\'egulier et non localement modulaire. 
Alors il existe un  sous-groupe d\'efinissable du groupe additif  dont le type g\'en\'erique est r\'egulier 
et non-orthogonal \`a $p$. }\\  
\smallskip Pour montrer ce r\'esultat on d\'efinit d'abord la suite de prolongations 
d'une
$(\si,D)$-vari\'et\'e $V$ : grosso modo, pour chaque $l$, on prend la
cl\^oture de Zariski $V_l$ de l'ensemble $$\{(a,\si(a),\ldots,\si^l(a),
Da,\si(Da),\ldots, \si^l(D^l(a)):a\in V\},$$ avec les projections
naturelles $V_{l+1}\to V_l$. La difficult\'e essentielle r\'eside dans l'identification des
conditions n\'ecessaires et suffisantes pour qu'une telle suite de
vari\'et\'es alg\'ebriques soit la suite de prolongations associ\'ee \`a
une $(\si,D)$-vari\'et\'e. Ensuite nous d\'efinissons les espaces
d'arcs de la fa\c con suivante : si $m\in\na$, nous consid\'erons la suite d'espace
d'arcs ${\mathcal A}_mV_l$ avec les projections naturelles. Nous
rencontrons cependant une difficult\'e, car il faut enlever les points
singuliers des vari\'et\'es $V_l$, ainsi que d'autres points. Cela nous am\`ene \`a
d\'efinir la notion de points non-singuliers d'une $(\si,D)$-vari\'et\'e,
et \`a d\'efinir l'espace d'arcs ${\mathcal A}_mV_a$ seulement au-dessus
des points non-singuliers, comme \'etant la $(\si,D)$-vari\'et\'e
associ\'ee \`a une suite de prolongations locales au-dessus de $a$. Notons que les points
singuliers forment un ensemble qui est \`a priori une 
union d\'enombrable de ferm\'es de la vari\'et\'e $V$, et nous ne
pouvons donc pas les enlever de fa\c con diff\'erentiellement birationnelle. De m\^eme, nous
d\'efinissons l'espace tangent $T(V)_a$ au-dessus d'un point
non-singulier $a$ en 
utilisant comme suite de prolongations la suite des espaces tangents
au-dessus de $a$, $(a,\si(a),Da,\si(Da))$, etc. Comme dans le cas des
vari\'et\'es alg\'ebriques, on montre alors que 
les fibres des projections ${\mathcal A}_{m+1}V_a\to {\mathcal A}_mV_a$
sont des espaces homog\`enes pour $T(V)_a$. Ce r\'esultat est un des
ingr\'edients fondamentaux de la preuve.

Dans le cas d'un type de $\s$-rang $1$, la modularit\'e locale et la
monobas\'eit\'e co\"\i ncident. Nous montrons aussi que  le
g\'en\'erique d'un sous-groupe d\'efinissable  du groupe additif est de dimension finie si et seulement si il est de
rang $\s$ fini. En combinant tous ces r\'esultats, nous obtenons alors
la dichotomie :\\ 

 \smallskip\noindent 
{\bf Th\'eor\`eme \ref{arc19}} {\it Soit $p$ un type de $\s$-rang 1. Si $p$ n'est pas mono-bas\'e, 
 alors il est non-orthogonal \`a $Fix \si \cap {\mathcal C}$. }\\   

  \bigskip\noindent 

\begin{large}{\bf Chapitre 4}\end{large}\\
  
\medskip 
Dans le quatri\`eme chapitre nous \'etudions quelques classes de groupes  d\'efinissables dans un mod\`ele de {\it DCFA}. 
 %Nous commen\c cons par
%quelques remarques faciles, et montrons que si $H$ est un sous-groupe ${\mathcal L}_{\si,D}$-d\'efinissable
 %d'un groupe alg\'ebrique $G$, alors il est d'indice fini dans sa cl\^oture pour la $(\si,D)$-topologie. 
%Ce r\'esultat est obtenu gr\^ace \`a une caract\'erisation des types g\'en\'eriques du groupe. 
Dans 4.1 nous utilisons les techniques de Kowalski-Pillay  pour montrer qu'un groupe d\'efinissable est 
isog\`ene \`a un sous-groupe d\'efinissable d'un groupe alg\`ebrique. Cela ram\`ene l'\'etude des groupes 
d\'efinissables \`a celle des sous-groupes d\'efinissables d'un groupe alg\'ebrique.\\    

 \smallskip\noindent 
{\bf Th\'eor\`eme \ref{GIII6}} {\it Soient $({\mathcal U},\si,D)$ un mod\`ele de {\it DCFA},
 $K \prec {\mathcal U} $ et $G$ un groupe $K$-d\'efinissable. Alors il existe un groupe alg\'ebrique  $H$,
 un sous-groupe d\'efinissable d'indice fini $G_1$  de $G$, et un isomorphisme  d\'efinissable entre $G_1/N_1$ 
et $H_1 / N_2$, o\`u  $H_1$ est un sous-groupe d\'efinissable  de $H(\mathcal{U})$, 
 $N_1$ est un sous-groupe normal fini de $G_1$, et  $N_2$ est un  sous-groupe normal fini de $H_1$. }\\

% Quand $G$ est d\'efini sur les constantes il est d'usage de prendre $s=0$, en particulier
%si 
%$G=\gr_m$, $lD(x)=\frac{Dx}{x}$.\\ 
  
\smallskip\noindent 
%{\bf Th\'eor\`eme \ref{cold}} {\it Soit $({\mathcal U},\si,D)$ un mod\`ele de {\it DCFA}, 
%soit $K=acl(K) \subset {\mathcal U}$. Soit $G$ un groupe alg\'ebrique commutatif connexe d\'efini sur
 %$K$ et soit $s:G \to {\mathcal U}^n$  un homomorphisme alg\'ebrique tel que pour tout $g$ in $G$,
 %$(g,s(g)) \in \tau(G)$.  Soit $H_1$ un sous-groupe connexe de $G$ $ {\mathcal L}_{\si}$-d\'efinissable
 %sans quantificateurs,  $H_2$  un sous-groupe connexe de $T(G)_e$ $ {\mathcal L}_{\si,D}$-d\'efinissable 
%sans quantificateurs. 
%Soit $G_1$ un sous-groupe alg\'ebrique connexe de  $G \times \cdots \times \si^k(G)$ tel que
% $H_1=\{g \in G:(g, \si(g), \cdots, \si^k(g) ) \in G_1 \}$. 
%Soit $H'= \{h \in T(G)_e : (h,\si(h),\cdots,\si^k(h)) \in T(G_1)_{(e,\cdots,e)}\}$.  
%Supposons que $(s,s^{\si},\cdots, s^{\si^k})$ se restreint \`a une  section de $G_1$.
% Si $H_2 < H'$, alors les \'equations $x \in H_1$ et $lD_s(x) \in H_2$ engendrent un $(\si,D)$-id\'eal premier.
% C'est \`a dire, il existe  $g \in H_1$ tel que  $lD_sg \in H_2$ est un 
%$(\si,D)$-g\'en\'erique de $H_2$ sur $K$ et $g$ est un $\si$-g\'en\'erique de $H_1$ sur $K(lD_s(g))_{\si,D}$. }\\
 
\smallskip\noindent
 La deuxi\`eme  section est d\'edi\'ee \`a l'\'etude des groupes commutatifs.
 Nous nous int\'eressons surtout \`a la propri\'et\'e d'\^etre mono-bas\'e,
 et c'est pourquoi nous restreignons notre attention aux groupes commutatifs.  
Gr\^ace \`a un th\'eor\`eme de Wagner, cette \'etude se r\'eduit \`a l'\'etude des sous-groupes du groupe additif, 
sous-groupes du groupe multiplicatif,  et sous-groupes d'une vari\'et\'e Ab\'elienne simple.     
Pour les groupes additifs on montre le r\'esultat suivant:\\ 

 \smallskip\noindent 
{\bf Proposition \ref{adgr}} {\it Aucun sous-groupe d\'efinissable infini 
du groupe additif $\gr_a$ n'est   monobas\'e.  }\\ 

 Cela provient du fait que tout sous-groupe d\'efinissable de $\gr_a$ est un $Fix \si \cap {\mathcal C}$-espace vectoriel. 
 Dans le cas du groupe multiplicatif on utilise la d\'eriv\'ee logarithmique et un r\'esultat de Hrushovski 
sur les sous-groupes  d\'efinissables de $\gr_m$ dans {\it ACFA}. Nous obtenons:  

\smallskip\noindent {\it Soit $H$ un sous-groupe d\'efinissable de $\gr_m$.
 Si $H\not\subset   \gr_m({\mathcal C})$, alors $H$ n'est pas monobas\'e.  
Si $H<\gr({\mathcal C})$, alors $H$ est monobas\'e si et seulement si il est monobas\'e dans le corps 
de diff\'erence $\mathcal C$, et dans ce
cas il sera aussi stable stablement plong\'e. }\\

Les r\'esultats de  Hrushovski \cite{HMM} nous donnent alors une description compl\`ete des sous-groupes
 monobas\'es de $\gr_m(\mathcal C)$. 
 \smallskip\noindent Dans le cas des vari\'et\'es Ab\'eliennes simples on regarde le noyau de Manin, qui existe
pour n'importe quelle vari\'et\'e Ab\'elienne.\\  

\smallskip\noindent 
{\bf Proposition \ref{ab6}} {\it Soit $A$ une  vari\'et\'e Ab\'elienne. 
Alors il existe un homomorphisme $D$-d\'efinissable  $\mu: A \to \gr_a^n$, 
o\`u $n=dim(A)$ , tel que $Ker(\mu)$ a rang de Morley fini. }\\  

$Ker (\mu)$ est la $D$-cl\^oture de $Tor(A)$, on appelle $Ker(\mu)$ le noyau de Manin de $A$ et on le note $A^{\sharp}$. 
Quand $A$ est simple $A^{\sharp}$ est un sous-groupe minimal.\\  

{\bf D\'efinition}
{\it On dit qu'une vari\'et\'e  Ab\'elienne  descend aux constantes si elle est
  isomorphe \`a une vari\'et\'e Ab\'elienne d\'efinie sur le corps des constantes.
}\\

Soit $A$ une vari\'et\'e Ab\'elienne simple.
 Nous distinguerons deux cas :  quand la vari\'et\'e $A$ descend aux constantes  et quand elle ne le
fait pas. Si $A$ est d\'efinie sur ${\mathcal C}$, alors $A^\#=A({\mathcal C})$, et
sinon, alors $A^\#$ est monobas\'e (pour la th\'eorie $DCF$), et donc
sera fortement minimal. 
Dans les deux cas on utilisera la dichotomie et les  r\'esultats de la section 5 du chapitre 2.
Le th\'eor\`eme suivant nous donne  une description de la mono-bas\'eit\'e des sous-groupes 
d'une vari\'et\'e Ab\'elienne: \\

\smallskip\noindent
{\bf Th\'eor\`eme \ref{absum}}
{\it 
Soit $A$ une vari\'et\'e Ab\'elienne simple, et soit $H$ un sous-groupe  de $A({\mathcal U})$ 
 d\'efinissable sans quantificateurs sur $K=acl(K)$. Si $H\not\subset A^{\sharp}({\mathcal U})$,
 alors $H$ est non mono-bas\'e. Suppposons maintenant que  $H\subset A^{\sharp}({\mathcal U})$, 
et soit $a$ un g\'en\'erique de  $H$ sur $K$. Alors  
\begin{enumerate}  
\item Si $A$  descend aux constantes, 
alors  $H$ est mono-bas\'e si et seulement si $H$ est stable stablement plong\'e,  si et seulement si
 $tp_{ACFA}(a/K)$ est h\'er\'editairement orthogonal \`a  $(\si(x)=x)$.   
\item Si $A$ ne descend  pas aux constantes,  
alors $H$ est mono-bas\'e. De plus   
\begin{enumerate} 
\item Si quelque soit $k$, $A$ n'est pas isomorphe \`a une vari\'et\'e Ab\'elienne d\'efinie sur 
 $Fix \si^k$, alors $H$ est stable stablement plong\'e.  
 \item Supposons que $A$ est d\'efinie sur $Fix(\si)$. Alors $H$ est stable  stablement plong\'e
 si et seulement si $tp_{ACFA}(a/K)$ est stable stablement  plong\'e. 
\end{enumerate} \end{enumerate}
}
Les r\'esultats de \cite{HMM} nous donnent alors une description  compl\`ete des   
sous-groupes $H$ qui ne sont pas mono-bas\'es dans les cas (1) et
 (2)(b) ci-dessus.

\setcounter{secnumdepth}{2}
\cleardoublepage
\chapter{Preliminaries}
\label{chap:prelim}

\section{Stable and Simple Theories}\label{pre:sec1}
We mention the results and definitions on stable and simple
theories that we will need in the following chapters. We assume that the reader
is familiar with the basic notions of stability and simplicity.
For more details see  \cite{pstab} and \cite{pgeostab} for stable theories, and
\cite{kimpi} and \cite{wag} for simple theories.\\

We will focus on $\omega$-stable theories.
A theory is said to be $\omega$-stable
if, given any model $M$, and  countable set $A \subset M$, there are countable many complete types over $A$.

We will always work in a countable language, and in that case being $\omega$-stable is equivalent to being
totally transcendental, that is, any definable set of a model of the theory has Morley rank.\\

We will consider $T$ a complete  $\omega$-stable theory over a countable language ${\mathcal L}$,
and $M$ a saturated model of $T$.
We denote by $S_n(A)$ the space of $n$-types over the set $A$.
One of the main properties of $\omega$-stable theories is that every type is definable: If $A=acl^{eq}(A) \subset M^{eq}$ and
$p \in S_n(A)$, then for each formula $\varphi(x,{\bar y})$, where $x$ is an
$n$-tuple,
 there is an
${\mathcal L}(A)$-formula $d_{p}\varphi({\bar y})$ such that for all tuple ${\bar a}$ of $A$,
$\varphi(x,{\bar a}) \in p$ if and only if $d_{p}\varphi({\bar a})$ is satisfied in $M$.

The notion of canonical base in stable theories will be useful for our case.

\begin{defi}\label{PST1}
Let ${\mathbf p}$ be a global type (a type over $M$). The canonical base of ${\mathbf p}$, $Cb({\mathbf p})$ is
a tuple of $M^{eq}$ which is fixed pointwise precisely by those automorphisms of $M$ which fix ${\mathbf p}$.
\end{defi}

\begin{prop}
Any global type has a canonical base which is unique up to interdefinability.
\end{prop}

%In $\omega$-stable theories we can define a notion of independence (and thus of forking) based on
%the Morley rank. With this independence relation we can define another rank, Lascar's $U$-rank.\\

Let $A \subset M$. Let $p \in S(A)$. We say that $p$ is
stationary if for every $B \supset A$, $p$ has  a unique non-forking extension to $B$.
For a stationary type $p \in S(A)$ we define the canonical base of $p$, $Cb(p)$ as
the canonical base of the unique non-forking extension of $p$ to $M$.

\begin{defi}\label{stp}
Let $A \subset M $ and $a \in M$. The strong type of $a$ over $A$, $stp(a/A)$ is the type of
$a$ over $acl^{eq}(A)$.
\end{defi}

\begin{rem}\label{cb1} \hspace{20cm}
\begin{enumerate}
\item Let $A \subset M$, and let $p \in S(A)$ be stationary. Then $Cb(p)$ is the smallest definably closed
subset $c$ of $dcl^{eq}(A)$ such that $p$ does not fork over $c$ and the restriction of $p$ to $c$ is stationary.
\item Let $a$ be a finite tuple of $M$, and let $A \subset B \subset M$. Then $tp(a/B)$ does not fork
over $A$ if $Cb(stp(a/B)) \subset acl^{eq}(A)$.
\item  Let $A \subset M$, and let $p \in S(A)$ be stationary. Then there is $N \in \na$ such that if $a_1,\cdots,a_N $ are independent realizations  of $p$
 then $Cb(p) \subset dcl^{eq}(a_1,\cdots,a_N) $.
\end{enumerate}
\end{rem}

Now we will concentrate on simple theories, more specifically on supersimple theories.
To define simplicity we need the notions of dividing and forking of types, however in \cite{kimpi}, Kim and Pillay
proved an equivalence which allows us to define simplicity with the notion of independence between tuples.
%In order to describe supersimple theories we must define a property known as the Independence Theorem.

%\begin{defi}\label{it}
%Let $T$ be a complete theory. We say that $T$ satisfies the Independence Theorem if
%for any  saturated model $M$ of $T$,  if $N \subset M$  is a model of $T$ such that $|N|<|M|$,
% and  ${\bar a}, {\bar b}, {\bar c_1}, {\bar c_2} $ are tuples of $M$ such that:
%
%\begin{enumerate}
%\item $tp(  {\bar c_1}/N  )=tp(  {\bar c_2}/N  )    $.
%
%\item ${\bar a}  \dnfo_N    {\bar c_1}$, ${\bar a}  \dnfo_N    {\bar b}$ and
% ${\bar b}  \dnfo_N    {\bar c_2}$.
%\end{enumerate}
%
%Then there is $ {\bar c} \in M $ realizing $tp(  {\bar c_1}/N \cup   {\bar a}  ) \cup   tp(  {\bar c_2}/N \cup   {\bar b}  )  $
% such that $ {\bar c}  \dnfo_N    ({\bar a} ,  {\bar b})  $.
%\end{defi}

\begin{teo}\label{PSIM1}
Let $T$ be a complete theory. $T$ is supersimple if and only if, given a large, saturated model $M$ of $T$, there is
an independence relation $\dnfo$ between tuples of $M$ over  subsets of $M$ which satisfies the following properties:
\begin{enumerate}
\item Invariance: Let $a$ be a tuple of $M$, and $B,C \subset M$ such that $a \dnfo_C B $;
let $\tau$ be an automorphism of $M$.
Then $\tau(a) \dnfo_{\tau(C)} \tau(B)$.
\item Local: For all finite tuple $a$, and for all $B$ there is a finite subset $C$ of $acl^{eq}(B)$ such that
$a \dnfo_C B $
\item Extension: For all tuple $a$, for all set $B$ and for all $C$ containing $B$ there is a tuple $a'$ such that
$tp(a/B)=tp(a'/B)$ and $a' \dnfo_B C$.
\item Symmetry: For all tuples $a,b$ and for all $C$, $a \dnfo_C b $ if and only if $b \dnfo_C a $.
\item Transitivity: Let $a$ be a tuple and let $A \subset B \subset C$. Then $a \dnfo_B C$  and $a \dnfo_A B $
if and only if $a \dnfo_A C $.
\item \label{it}Independence Theorem:
If $N \prec M$ is such that $|N|<|M|$,
 and  ${\bar a}, {\bar b}, {\bar c_1}, {\bar c_2} $ are tuples of $M$ such that:

\begin{enumerate}
\item $tp(  {\bar c_1}/N  )=tp(  {\bar c_2}/N  )    $.

\item ${\bar a}  \dnfo_N    {\bar c_1}$, ${\bar a}  \dnfo_N    {\bar b}$ and
 ${\bar b}  \dnfo_N    {\bar c_2}$.
\end{enumerate}
Then there is $ {\bar c} \in M $ realizing $tp(  {\bar c_1}/N \cup   {\bar a}  ) \cup   tp(  {\bar c_2}/N \cup   {\bar b}  )  $
 such that $ {\bar c}  \dnfo_N    ({\bar a} ,  {\bar b})  $.

\end{enumerate}
Furthermore, the relation $\dnfo$ then coincides with non-forking.

\end{teo}

%Actually, \ref{PSIM1} is not the original definition of supersimple theory, the original definition is based on
%the concepts of dividing and forking of types. In \cite{kimpi} Kim and Pillay proved the equivalence between the two definitions
%and that
%an independence relation satisfying the conditions of \ref{PSIM1} is unique and coincides with non forking. 
%The results in \cite{kimpi} implies the following:

%\begin{prop}
%Let $T$ is a supersimple theory, $M$ a model of $T$, $a$ a tuple of $M$ and $A \subset B \subset M $.
%$tp(a/A)$ forks over $B$ if and only if $a \dfo_A B$.
%\end{prop}

Thus, in particular,
 we can define a rank in analogy to the Lascar rank for superstable theories.

\begin{defi}\label{PSU1}
Let $T$ be a supersimple theory and let $M$ be a model of $T$.
Let $A \subset M$ and  ${\bar a}$ be a tuple of  $M$.
 We define the SU-rank of $tp({\bar a}   /A) $, $\s({\bar a}/A)$, by induction
   as follows:

\begin{enumerate}
\item $\s( {\bar a} /A   )  \geq 0 $.

\item For an ordinal $\alpha$,  $\s({\bar a}   /A   )  \geq \alpha +1 $
if and only if there is  a forking extension $q$   of $tp({\bar a}   /A)  $ such that $\s(q) \geq \alpha $.

\item If $\alpha$ is a limit ordinal, then  $\s( {\bar a} /A)     \geq \alpha $ if and only if
  $\s({\bar a} /A   )  \geq \beta $ for all $\beta \in \alpha   $.

We define  $\s( {\bar a} /A   ) $ to be the smallest ordinal $\alpha  $
 such that  $\s( {\bar a} /A   )  \geq \alpha $ but
$\s({\bar a} /A)    \geq \! \! \! \!  \!  /   \, \,   \alpha +1 $.

\end{enumerate}

\end{defi}

Lascar's inequalities hold for the $\s$-rank. We need the natural sum of ordinals.
\begin{defi}\label{PSU2}
  If $\alpha  ,\beta $ are ordinal numbers, we can write in a unique way
  $\alpha= \omega^{\alpha_1}a_1+\cdots + \omega^{\alpha_k}a_k $ and
 $\beta= \omega^{\beta_1}b_1+\cdots + \omega^{\beta_l}b_l $ where
$\alpha_1 > \cdots > \alpha_k $, $\beta_1> \cdots > \beta_l $ are ordinals,
and $a_1,\cdots ,a_k,b_1,\cdots,b_l $ are positive integers.
If we allow them to be zero we can suppose $k=l$, and $\alpha_i=\beta_i$ for all $i$.
 We define the natural sum of ordinals by
 $\alpha \oplus \beta = \omega^{\alpha_1}{(a_1+b_1)}+\cdots + \omega^{\alpha_k}{(a_k+b_k)} $.
\end{defi}

\begin{fac}\label{PSU3}
Let $T$ be a supersimple theory and $M$ a model of $T $. Let   $a,b $ be tuples of $M  $, $A,B  $ subsets of $M$.
 Let $\alpha   $ be an ordinal. Then:

\begin{enumerate}
\item $\s(a/Ab) + \s(b/A) \leq \s(ab/A) \leq  \s(a/Ab) \oplus  \s(b/A).  $
\item If $\s(a/A) \geq \s(a/Ab) \oplus \alpha  $, then $\s(b/A) \geq \s(b/Aa)+\alpha$.
\item If $\s(a/A) \geq \s(a/Ab) + \omega^{\alpha}   $,
then $\s(b/A) \geq \s(b/Aa)+\omega^{\alpha}.$
\item If $a \dnfo_A b  $, then $\s(ab/A)=\s(a/A) \oplus \s(b/A) . $

\end{enumerate}

\end{fac}

For simple theories there is a notion of canonical base, which is defined with the help of amalgamation bases.
However, for the case of difference-differential fields we can avoid this definition, and use the notion of
quantifier-free  canonical base which we will define later. \\

\section{One-basedness, analyzability and local modularity}
We introduce  some helpful concepts and results concerning supersimple theories. 
They are analogous definitions for stable theories. For more details see \cite{wag}.

Let $M$ be a saturated model of a supersimple theory $T$ which eliminates imaginaries.

\begin{defi}
Let $A,B \subset M $, $p \in S(A),q \in S(B) $.
\begin{enumerate}
\item Let us suppose that $A=B$. We say that $p$ is almost orthogonal to $q$, denoted by $p \perp^a q $,
if any two realizations of $p$ and $q$ are independent over $A$.
\item We say that $p$ is orthogonal to $q$, denoted by $p \perp q  $, if for every set $C \supset A \cup B  $
 and for every two extensions $p',q' \in S(C)  $ of $p,q$ respectively, such that $p' $ does not fork
 over $A$ and $q'$ does not fork over $B$, we have that $p' \perp^a q'  $.
\item Let  $\varphi(x) \in {\mathcal L}(B)  $. We say that $p $ is orthogonal to $\varphi$,
 denoted by $p \perp \varphi  $, if for every type $q  $ over $B$ containing $\varphi$,
 we have  $p \perp q $.
\item If $\s(p)=1 $, we say that $p$ is trivial if for every $C \supset A$ and $a_1,\cdots,a_n $
realizing non-forking extensions of $p$ to $C$, the tuples $a_1,\cdots,a_n $ are independent over $C$ if and
only if they are pairwise independent.
%\item We say that $p$ has weight $1$, denoted $w(p)=1 $, if, for every realization
%$a$ of $p$ the following holds:
% for any tuples $b$ and $c$, if $tp(a/Eb) $ and $tp(a/Ec) $ fork over $E$
 %then $b \not\dfo_E c $.
\end{enumerate}
\end{defi}

%\begin{prop}
%Non-orthogonality is an equivalence relation between the types of weight 1.
%\end{prop}

\begin{prop}\label{wgeq}
Let  $a \in M  $ and $A \subset M  $.
Let us suppose that $SU(a/A)= \beta + \omega^{\alpha} \cdot n $, with $n >0 $ and
  $\omega^{\alpha+1} \leq \beta < \infty $ or $\beta =0  $.
Then $tp(a/A)  $ is non-orthogonal to a type of $\s$-rank $\omega^{\alpha}$. Moreover there is $b\in acl(Aa)$
with $\s(b/A)=\omega^{\alpha}n$.
\end{prop}

\begin{defi} \hspace{20cm}
\begin{enumerate}
\item Let $A \subset M$ and let $S$ be an $(\infty)$-definable set over $A$. We say that $S$ is 1-based if
for every $m,n \in \na $, and $a\in S^m,b \in S^n $, $a$ and $b$ are independent over
$acl(Aa) \cap acl(Ab)$.
\item A type is 1-based if the set of its realizations is 1-based.
\end{enumerate}
\end{defi}

\begin{prop}\textup{(\cite{wagbase})}\label{wgr}
\begin{enumerate}
\item The  union of $1$-based sets is $1$-based.
\item If $tp(a/A)$ and $tp(b/Aa)$ are $1$-based, so is $tp(a,b/A)$.
\end{enumerate}
\end{prop}

\begin{defi}\label{prt1}
Let $p \in S(A)$. We say that $p$ is regular if $p$ is orthogonal to all its forking extensions.
%
%for every $B \supset A$ and elements $a,b$ such that $tp(a/B)$ is a non-forking extension of $p$ and
%$tp(b/B)$ is a forking extension of $p$, we have $a \dnfo_B b$
\end{defi}

\begin{prop}\label{prt2}\textup{(5.2.12 of \cite{wag})}\\
Non almost-orthogonality over a set $A$ is an equivalence relation on regular types over $A$.
\end{prop}

\begin{defi}\label{prt3}
Let $p, \, q \in S(A)$. We say that $q$ is $p$-internal if for every realization
$a$ of $q$ there is a set $B$ such that $B \dnfo_A a$ and a tuple $c$
of realizations of $p$ such that $a \in dcl(Bc)$.
A set $X$ definable over $A$ is $p$-internal if for every tuple $a$
of $X$, $tp(a/A)$ is $p$-internal.
If we replace $dcl$ by $acl$ above we say that $q$ (or $X$) is almost $p$-internal.
\end{defi}

\begin{defi}\label{prt7}
Let $p,q \in S(A)$. We say that $q$ is $p$-analyzable if there is a realization
$a$ of $q$, an ordinal $\kappa$ and $(a_i)_{i < \kappa} \subset dcl(A,a)$
such that $tp(a_i/A \cup \{a_j:j <i \})$ is $p$-internal for all $i < \kappa$.
\end{defi}

\begin{defi}
Let $p$ and $q$ be two complete types. We say that $q$ is hereditarily orthogonal to $p$ if every extension of $q$
is orthogonal to $p$.
\end{defi}

\begin{rem}\label{prt4}
If $tp(a/A)$ is nonorthogonal to $p \in S(A)$ then $acl(Aa)$ contains a $p$-internal set:
Let $b$ realize $p$, $B$ such that $a \dnfo_A B$, $b \dnfo_A B$ and $b \dfo_B a$.
Then $Cb(Bb/Aa)$ realizes a $p$-internal type over $A$.

Let $A_p$ denote the maximal
algebraic closed subset of $acl(Aa)$ such that $tp(A_p/A)$ is almost $p$-internal.
If $b$ realizes $p$ then $a \dnfo_{A_p} b$, and therefore $a \dfo_A b$ implies $A_p \dfo_A b$.
\end{rem}

\begin{defi}\label{prt5}
Let $p \in S(A)$ be regular and let $q$ be a type over a set $X \supset A$. 
We say that $q$ is $p$-simple if there is $B \supset X$ and 
a set $Y$ of realizations of $p$ and a realization $a$ of $q$ with $a \dnfo_X B$ such that $tp(a/BY)$ is 
hereditarily orthogonal to $p$. 
%$p \not\perp B$.
%We say that $q$ is $p$-simple if there is a realization of a of $q$ such that
%every extension of $tp(a/A_p)$ is orthogonal to $p$.
%there is a set
%$C \supset A \cup B$ a realization $a$ of $q$ over $C$ and a set of realizations
%$X$ of $p$ such that every extension of $stp(a/CX)$ is orthogonal to $p$.
\end{defi}

\begin{defi}
Let $p$ be a (possibly partial) type over $A$ and $q=tp(a/B)$ a type. 
The $p$-weight of $q$, denoted by $w_p(q)$, is %the
%smallest cardinal $\lambda$ such that there is $C \supset A \cup B$, a realization $a'$ of a non forking extension
%$q'$ of $q$ to $C$, and an independent set $I$ of realizations  of $p$ such that $tp(a'/CI)$ is orthogonal
%to $p$ and $|I|=\lambda$.
the largest integer $n$ such that there are $C\supset A\cup B$, a tuple $a_1,\ldots,a_n$ of realisations of $p$ 
which are independent over $C$, and a realisation $b$ of $q$ such that $(a_1,\ldots,a_n)\dnfo_A C$, $b\dnfo_B C$ 
and $a_i \dfo_C b$ for every $i=1,\ldots,n$. 
If $p$ is the partial type $x=x$ we say weight instead of $p$-weight and it is denoted by $w(q)$.
\end{defi}

\begin{defi}\label{prt11}
Let $A$, $B$ and $C$ be sets. We say that $A$ dominates $B$ over $C$ if for every set $D$, $D \dnfo_C A$ implies
$D \dnfo_C B$.
Let $p,q$ be two types. We say that $p$ dominates $q$
if there is a set $C$ containing the domains of $p$ and $q$ and realizations $a$ and $b$ of non forking extensions of
$p$ and $q$ to $C$ respectively, such that $a$ dominates $b$ over $C$.
We say that $p$ and $q$ are equidominant if $p$ dominates $q$ and $q$ dominates $p$.
\end{defi}

\begin{rem}\label{prt12}
Equidominance is an equivalence relation between regular types.
\end{rem}

\begin{rem}
Let $A \dnfo_B C$, and let $p \in S(A)$, $q \in S(B)$ and $r \in S(C)$ be regular types.
By the independence theorem, if $p \not\perp q$ and $q \not\perp r$ then $p \not\perp r$. 
\end{rem}

\begin{rem}\label{orto}
If $p \not\perp q$ and $q \not\perp r$ then there is a conjugate $r'$ of $r$ such that $p \not \perp r'$.
\end{rem}

\begin{defi}
Let $p_1,\cdots, p_n$ be types over a set $A$. The product $p_1 \times \cdots \times p_n$ 
is the partial type of $n$ independent
realizations of $p_1, \cdots, p_n$ over $A$.
\end{defi}

\begin{prop}\label{prt13}
A type in a supersimple theory is equidominant with a finite product of regular types.
\end{prop}

\begin{defi}\label{prt6}
Let $p$ be a regular type over $A$ and let $q$ be a $p$-simple type. 
We say that $q$ is $p$-semi-regular if it is domination equivalent to a non-zero
power of $p$
\end{defi}

%\begin{rem}
%The type $tp(a/A)$ is $p$-semi-regular if and only if it is dominated by $tp(A_p/A)$.
%\end{rem}

\begin{defi}\label{prt8}
Let $p$ be a type and $A$ a set. The $p$-closure of $A$, $cl_p(A)$ is the set of all $a$
such that $tp(a/A)$ is $p$-analysable and hereditarily orthogonal to $p$.
%all (hyperimaginaries) $a$ such that $tp(a/A)$ is $p$-analyzable and
%for any realization $b$ of $p$, for every set $B \dnfo_A a$ and for every
%realization $c$ of an extension of $tp(a/A)$ over $B$, $a \dnfo_{AB} c$.
\end{defi}

\begin{defi}\label{prt9}
A type $p$ is called locally modular if for any $A$ containing the domain of $p$,
 and any tuples $a$ and $b$ of realizations of $p$,
 we have $a \dnfo_C b$ where
$C=cl_p(Aa) \cap cl_p(Ab)$.
\end{defi}

\begin{rem}
In the definition of local modularity we can replace $cl_p$ by the following $p$-closure:
define $cl'_p(A)$ as the set of all $a$ such that $tp(a/A)$ is hereditarily orthogonal to $p$.
\end{rem}Indeed, let $a$ and $b$ be as in the definition of local modularity, let 
$c_0=cl_p(a) \cap cl_p(b)$ and $c_1=cl'_p(a) \cap cl'_p(b)$. 
Then $Cb(a,b,c_0/c_1)$ 
realises a $p$-internal type over $c_0$ (since it is contained in the
algebraic closure over $c_0$ of realisations of $tp(a,b)$), and
therefore it is contained in $cl_p(c_0)=c_0$. Thus $ab\dnfo_{c_0}c_1$,
and because $a\dnfo_{c_1}b$ and $a\dnfo_{c_0}c_1$, this implies that 
$a\dnfo_{c_0}b$.
%We know that $tp(c_1/c_0a)$ is orthogonal
%to every extension of $p$, then $c_1 \dfo_{c_0a} b$, and this implies that $a \dfo_{c_1} b$ if and
%only if $a \dfo_{c_0} b$. 
%\end{rem}

\begin{prop}\label{prt10}
A type $p$ is locally modular if and only if for any two models $M$ and $N$
with $N \prec M$, and any tuple of realizations $a$ of $p$ over $M$ such that $tp(a/N)$
is $p$-semi-regular, $Cb(a/M) \subset cl_p(Na)$.
\end{prop}

\begin{lem}\label{prt011}
Let $q=tp(a,b)$ be a regular type and $p=tp(b)$. Then $p$ is locally modular if and only if $q$ is locally modular.
\end{lem}

{\it Proof}:\\

As being hereditarily orthogonal to $p$ is the same
as being hereditarily orthogonal to $q$, by definition for any set $B$,
$cl'_p(B)=cl'_q(B)$.\\

Let $(a_1,b_1)$ and $(a_2,b_2)$ be tuples (of tuples) of realisations of $q$. 
By regularity of $q$, $a_i\in cl'_p(b_i)$ for $i=1,2$, so that 
$$cl'_q(a_1,b_1)\cap cl_q(a_2,b_2)=cl'_p(b_1)\cap cl'_p(b_2)=_{\rm def}C.$$ 
It follows immediately that the local modularity of $q$ implies the local modularity of $p$. 
Conversely, assume that $p$ is locally modular. Then 
$b_1\dnfo_Cb_2$. Let $D=Cb(a_1b_1/acl(Cb_2))$. 
Then $tp(D/C)$ is almost-internal to the set of conjugates of $tp(a_1/Cb_1)$,
 and is therefore hereditarily orthogonal to $p$.
 Hence $D\subset cl'_p(C)=C$, and $a_1b_1\dnfo_Cb_2$. 
A similar reasoning gives that $Cb(a_2b_2/acl(Ca_1b_1))\subset C$.\\ 
%Then for any $c,d$, $c \dnfo_C d$ where $C=cl_p(c) \cap cl_p(d)$ if and only if
%$c \dnfo_D d$ where $D=cl_q(c) \cap cl_q(d)$.
%If $q$ is locally modular, let $c,d$ be realizations of $p$, and let $c',d'$ as independetn as possible 
%such that $(c,c'),(d,d')$ are
%realizations of $q$. Then $(c,c') \dnfo_{D'}(d,d')$ where $D'=cl_q(cc') \cap cl_q(dd')$.
%Assume that $p$ is locally modular. Let $(c,c'),(d,d')$ be realizations of $q$. Then $c \dnfo_C d$.
$\Box$

\begin{lem}\label{prt012}
Let $T$ be a supersimple theory which eliminates imaginaries, $A=acl(A)$ a subset of 
some model $M$ of $T$ and $a$ a tuple in $M$.
Assume that $tp(a/A)$ has SU-rank $\beta+\omega^\alpha=\beta\oplus \omega^\alpha$ and has weight 
$1$. Then there is $b\in acl(Aa)$ such that $SU(b/A)=\omega^\alpha$.
%Let $p=tp(a/A)$ be regular type of $\s$-rank $\beta+\omega^{\alpha}$. Then
%there is $b \in acl(Aa)$ with $\s(b/A)=\omega^{\alpha}$.
\end{lem}

{\it Proof}:\\

By \ref{wgeq}, there is some $C=acl(C)\supset A$ independent from $a$ 
over $A$ and a tuple $c$ such that $\s(c/C)=\omega^\alpha$, and $c$ and 
$a$ are not independent over $C$. Let $B$ be the algebraic closure of 
$Cb(Cc/acl(Aa))$. Then $B$ is contained in the algebraic closure of 
finitely many (independent over $Aa$) realisations of $tp(Cc/acl(Aa))$, 
say $C_1c_1,\ldots,C_nc_n$. Let $D=acl(C_1,\ldots,C_n)$. Then $D$ is 
independent from $a$ over $A$, and each $c_i$ is not independent from 
$a$ over $D$. Since $tp(a/A)$ has weight $1$, so does $tp(a/D)$, and 
therefore for each $1<i\leq n$, $c_1$ and $c_i$ are not 
independent over $D$. Thus $\s(c_i/Dc_1)<\omega^\alpha$, and therefore 
$SU(c_1,\ldots,c_n/D)<\omega^\alpha 2$. As $D$ is independent from $a$ 
over $A$, and $B \subset acl(D,c_1,\cdots,c_n)\cap acl(Aa)$, we get 
$\s(B/A)<\omega^\alpha 2$. 
Since $\s(c/C)=\omega^{\alpha}$ and $\s(c/CB)< \omega^{\alpha}$,
then $\s(B/C)\geq \omega^{\alpha}$, and as $B \dnfo_A C$ we have $\s(B/A)\geq \omega^\alpha$. 
By Lascar's inequalities \ref{PSU3} we have $\s(a/AB) + \s(B/A) \leq \beta + \omega^{\alpha}$.
As $\s(B/A)\geq \omega^\alpha$ we have that $\s(B/A)=\delta + \omega^{\alpha}$ with $\delta \geq\omega^{\alpha}$ or 
$\delta =0$, 
and $\s(B/A)<\omega^\alpha 2$ implies that  $\delta = 0$.

$\Box$

%\begin{prop}\label{prt14}
%Let $p \in S(A)$ be a regular non-locally modular type. Then there are two tuples of realizations of $p$, $a$ and $b$
%such that $a \dfo b$ over $cl_p(Aa) \cap cl_p(Ab)$, $b=Cb(a/Ab)$ and $tp(b/Aa)$ is $p$-semi-regular.
%\end{prop}
%
%{\it Proof}:\\
%
%Let $a,b$ be realizations of $p$ such that $a \dfo b$ over $cl_p(Aa) \cap cl_p(Ab)$.
%If $tp(b/Aa)$ is non $p$-semi-regular
\section{$\omega$-stable theories with an automorphism}\label{pre:sec3}

Now we take a look at $\omega$-stable theories with an automorphism. All this material is from \cite{zopi}.

Let $T$ be an $\omega$-stable ${\mathcal L}$-theory which eliminates quantifiers and imaginaries. 
Let ${\mathcal L}_{\si}={\mathcal L} \cup \{{\si} \}$,
where $\si$ is an 1-ary function symbol.

\begin{defi}
Let $M$ be a saturated model of $T$, and $N$ a model of $T$. Let  $\si$ an automorphism of $N$.
We say that $\si$ is generic if for any algebraically closed structures $A,B \subset M$,
and $\si_1,\si_2$ automorphisms of $A$ and $B$ respectively; if $(A,\si_1) \subset (N,\si)$ and
$f:(A,\si_1) \to (B,\si_2)$ is an ${\mathcal L}_{\si}$-embedding,
then there is an
${\mathcal L}_{\si}$-embedding $g:(B,\si_2) \to (N,\si)$ such that $g \circ f$ is the identity on $A$.
\end{defi}
We consider the theory $T_0$ whose models are the ${\mathcal L}_{\si}$-structures
 of the form $(M,\si)$, where
$M$ is a model of $T$ and $\si$ is an ${\mathcal L}$-automorphism of $M$.
We denote by $T_1$ the ${\mathcal L}_\si$-theory of models of $T_0$ where $\si$ is generic.

\begin{teo}\label{Paut1}
Assume that $T_0$ has a model-companion $T_A$.
\begin{enumerate}
\item $T_A=T_1$, and $(M,\si) \models T_A$ if and only if $(M,\si)$ is generic.
\item If $(M_1,\si_1)$ and  $(M_2,\si_2)$ are models of $T_A$ containing a common
 algebraically closed substructure $(A,\si)$, then  $(M_1,\si_1) \equiv_A (M_2,\si_2)$.
\end{enumerate}
\end{teo}

Let $(M,\si)$ be a saturated model of $T_A$, and let $A \subset M$.
We denote by $acl_{\si}(A)$ the algebraic
closure in the sense of $T$ of $(\cdots, \si^{-1}(A),A, \si(A), \cdots)$.

\begin{prop}\label{Paut2}Assume $T_A$ exists and $(M,\si)$ is a model of $T_A$. Then
\begin{enumerate} 
\item $acl(A)=acl_{\si}(A)$.
\item Let $(N_1,\sigma_1),(N_2,\sigma_2)$ be two submodels of $(M,\si) $ containing a common
structure $E$.
Then  $   N_1 \equiv _E N_2      $ if and only if $ (acl(E),\sigma_1 )  \simeq  (acl(E),\sigma_2 ) $.
\item Every completion of $T_A$ is supersimple, and independence is given by:
$A\dnfo_C B$ if and 
only if $acl_\sigma(AC)$ and $acl_\sigma(BC)$ are independent over
$acl_\sigma(C)$ in the sense of the theory $T$.

\end{enumerate}
\end{prop}

%Finally we would like to talk about quantifier-free $\omega$-stable theories, which will be the case of the
%model-companion of the theory of difference-differential fields.
\begin{defi}
Let $T$ be a first order theory. We say that $T$ is quantifier-free $\omega$-stable if for any saturated model
$M$ of $T$, there are only countably many quantifier free types over a countable set.
\end{defi}
\begin{rem}\label{prcb}
Let $T$ be an $\omega$-stable theory which eliminates quantifiers and imaginaries, and let $\si$ be an automorphism of a model $M$
of $T$. 

Let $A=acl(A) \subset M$ and let $a \in M$. 
Then $qftp_{{\cal L}_\sigma}(a/A)$ is entirely determined by
$tp_T((\sigma^i(a))_{i\in {\mathbb Z}} /A)$. Let $B=dcl_T(A,\sigma^{-i}(a)|i>0)$, 
and
consider 
$tp_T(a/B)$. As $T$ is $\omega$-stable, there is some integer $n$ such
that $tp_T(a/B)$ is the unique non-forking extension of $tp_T(a/A,
\si^{-1}(a),\ldots,\si^{-n}(a))$ to $B$. Applying $\sigma^i$, this
gives that $tp_T(\sigma^i(a)/\sigma^i(B))$ is the unique non-forking
extension of $tp_T(\si^i(a)/A,\si^{i-1}(a),\ldots,\si^{i-n}(a))$ to
$\sigma^i(B)$.
% Thus, $tp_T((\sigma^i(a))_{i\in {\mathbb Z}} /A)$ is axiomatized 
%by
%$$\bigcup_{i\in {\mathbb Z}} tp_T(\si^i(a,\si^{i-1}(a),\ldots,\si^{i-n}(a))/A),$$
This implies that $T_0$ is quantifier-free-$\omega$-stable (and so
are $T_1$ and $T_A$).
%Let $A \subset M$, and let $a \in M$. 
%
%Then $qftp_{T_1}(a/A)=tp_{T}(a,\si(a),\cdots/A)$; as
% $T$  is $\omega$-stable, there is $n \in \na$ such that $tp_{T}(a,\si(a),\cdots/A)$ is
%determined by $tp_{T}(a,\si,(a), \cdots, \si^n(a)/A)$. Hence $T_1$ is quantifier-free $\omega$-stable.
\end{rem}

\begin{rem}{\bf Definition of canonical bases of quantifier-free 
types}. 
Assume that $T$ is $\omega$-stable and $T_A$ exists.
Since forking in $T_A$ is witnessed by quantifier-free
formulas (see \ref{Paut2}.3 ), any completion of $T_A$ is
supersimple, and $tp_{T_A}(a/A)$ does not fork over $C$, where
$C$ is the canonical base (in the sense of $T$) of
$tp_T((\sigma^i(a))_{i\in {\mathbb Z}} /A)$. If $T$ is an $\omega$-stable theory one  easily check  that
$C=dcl_T(\sigma^i(c)|i\in {\mathbb Z})$, where
$c=Cb(tp_T(\si^i(a,\si^{i-1}(a),\ldots,\si^{i-n}(a))/A)$ for some $n$ as in \ref{prcb}. So we will
define the canonical base 
of $qftp_{{\cal L}_\sigma}(a/A)$ to be this set $C$, and denote it by 
$Cb(qftp_{{\cal L}_\sigma}(a/A))$.

Note that $Cb(qftp_{{\cal L}_\sigma}(a/A))$ does not coincide with 
the canonical base of
$tp_{T_A}(a/A)$ as defined for simple theories. However, if $T_A$
satisfies  the 
independence theorem over algebraically closed sets, then the canonical
base of $tp_{T_A}(a/A)$ will be contained in\\
$acl(Cb(qftp_{{\cal
L}_\sigma}(a/A))$.  
%Let $T$ be an $\omega$-stable theory which eliminates quantifiers.
%Assume that $T_0$ (as above) has a model-companion $T_A$ in the language ${\mathcal L}_{\si}$. Let
%$(M,\si)$ be a model of $T_A$, let $A \subset M$ and let $a \in M$.
%Then forking in $T_A$ is determined by
%quantifier-free formulas.
%So we can define the canonical base of $p=qftp_{T_A}(a/A)$, $Cb(p)$, as the canonical base of
% $tp_T(a,\si(a), \cdots, \si^n(a)/A)$, with $n$ as above.We define the canonical base of a quantifier-free type up
%to definability.

% This canonical base is not the same as
%the canonical base of $tp_{T_A}(a/A)$ as a simple theory, but it satisfies \ref{cb1} and \ref{PST1}. 
\end{rem}

\section{$\omega$-Stable Groups}

We mention some facts and important definitions on $\omega$-stable groups. For more details on this subject the
reader may consult \cite{lecpill1}, \cite{poistable} or \cite{wagstable}.

Throughout this section $G$ will be an $\omega$-stable group, that is a group with possibly extra structure, such that
the theory of the group with this  structure is $\omega$-stable.
We will often assume that $G$ is $\emptyset$-definable. We point out that some of the results listed here hold
for stable groups.

An important property of $\omega$-stable groups is the chain condition.

\begin{prop}\label{G11}
There is no infinite strictly decreasing sequence of definable subgroups of $G$.
In particular, the intersection of any set of definable subgroups is equal to the intersection of finitely many of them,
 and thus  is definable. 
\end{prop}

If we consider now the class of definable subgroups of finite index in $G$, their intersection is definable, thus
there is a unique smallest definable subgroup of $G$ of finite index.
We say that a definable group is connected if it has no proper definable subgroup of finite index.
We call the intersection of all definable subgroups of $G$ of finite index the connected component of $G$ and
we denote it by $G^0$.
As an immediate consequence we have that, if $G$ is connected and $H$ is a normal definable subgroup of $G$, then
$G/H$ is connected.

We can define an action of $G$ on the space of types over $G$ as follows: let $p \in S(G)$ and let $g \in G$.
Then $g \cdot p= \{\varphi(x): \varphi(x)$ is a formula with parameters in $G, \varphi(gx) \in p\}$.

\begin{defi}
Let $G$ be an $\omega$-stable group. Let $p \in S(G)$. 
The left stabilizer of $p$ is $Stab(p)=\{g \in G: g \cdot p=p \}$.

If $G$ is  defined over $A$, and $a \in G$ is such that $p=tp(a/A)$ is stationary, the left
stabilizer of $p$ is $Stab(p)=\{g \in G:$ for some $a'$ realizing the 
non-forking extension of $p$ to $A \cup \{g\},$ $tp(g \cdot a'/A)=p  \}$
\end{defi}

Clearly $Stab(p)$ is a subgroup of $G$, and it is definable.\\
%We could also have defined  $Stab(p)$ as the set of $g \in G$ such that for some realization $a'$ of
%$p(x)|g$, $tp(ga')=p$.\\

We define now generic types of $\omega$-stable groups.

\begin{defi}
Let $\alpha=RM(G)$, where $RM$ denotes the Morley rank. Let $A$ be some set of parameters, and let
$p \in S(A)$. We say that $p$ is a generic type of $G$ over $A$ if $(x \in G) \in p$ and 
$RM(p)=\alpha$.
%(Here $RM$ denotes the Morley rank).
\end{defi}

$\omega$-stable groups always have generic types. If $p$ is a type over $A$ and $q$ is a non-forking
extension of $p$ to $B \supset A$, then $p$ is a generic type of $G$ over $A$ if and only if
$q$ is a generic type of $G$ over $B$. 
Let $a \in G$. Then $tp(a/A)$ is a generic of $G$ if and only if
$tp(a^ {-1}/A)$ is a generic of $G$.

\begin{lem}\label{gen1}
Let $a \in G$. Then $tp(a/A)$ is a generic type of $G$ if and only if for any $b \in G$ independent from $a$
over $A$, $a \cdot b$ is independent from $b$ over $A$, if and only if $b \cdot a$ is independent from $b$ over $A$.
\end{lem}

This is why we can talk about generics instead of left generics
and right generics.

\begin{prop}
If $G$ is connected then it has a unique generic type. In particular any $\omega$-stable group $G$ has finitely
many generic types, and these types correspond to the cosets of $G^ 0$ in $G$.
\end{prop}

\begin{lem}\label{needjet1}
Let $p(x) \in S(A)$ be a  stationary type containing $x \in G$. Define $H=Stab(p)$.
Let $a \in G$ and $b$ realize the non-forking extension of $p(x)$ to $A \cup \{a\}$. 
Then  $tp(a \cdot b/Aa)$ is stationary and  $aH$, as an imaginary, is interdefinable over $A$ 
 with $Cb(stp(a \cdot b/Aa))$.
\end{lem}

\begin{lem}\label{needjet2}
Let $c \in G$ such that $p=tp(c/A)$ is stationary. Let $H=Stab(p)$ and let $a \in G$ be a generic over
$A \cup \{c\}$. Then  $Hc$, as an imaginary, is interdefinable with $Cb(stp(a/A,c \cdot a))$ over $A \cup \{a \}$
\end{lem}

We end this section with a remarkable result on 1-based stable groups, due to Hrushovski and Pillay (\cite{phgr}).
 
\begin{teo}
$G$ is 1-based if and only if, for all $n \in \na$, every definable subset of $G^n$ is a finite Boolean
combination of cosets of definable subgroups of $G^n$.
\end{teo}

\section{Definable Groups in Supersimple Theories}

Throughout this section $T$ will be a supersimple theory, $M$ a saturated model of $T$, and $G$ an
$\infty$-definable (definable by an infinite number of formulas) group over some set of parameters $A \subset M$.
Since in simple theories we do not have Morley rank, we will use the equivalence \ref{gen1} to define generic types.
We refer to \cite{simpg} for the proofs. Some of the results in this section hold for groups definable in simple theories.

\begin{defi}
Let $p \in S(A) $. We say that $p$ is a left generic type of $G$ over $A$ if it is realized in $G$ and 
for every $a \in G$ and $b$ realizing $p$ such that $a \dnfo_A b  $, we have $b \cdot a \dnfo_A a$.
\end{defi}

%From this definition we have that a type $p$ is left-generic of and only if it is right-generic. We have also that
%$p$ is generic if and only if $p^ {-1}$ is generic. (By $p^ {-1}$ we mean $tp(a^{-1}/A)$ where $a$ realizes $p$). 

Some of the properties of generic types in $\omega$-stable groups hold in simple theories.

%The following results are proved in \cite{simpg} :

\begin{fac}

Let $G$ be an $A$-definable group.

\begin{enumerate}

\item Let $a,b \in G $. If $tp(a/Ab)  $ is left generic of $G$, then so is $tp(b \cdot a/Ab ) $.

\item Let $p \in S(A)  $ be realized in $G$, $B=acl(B) \supset A$, and $q \in S(B) $
 a non-forking extension of $p$. Then $p$ is a generic of $G$ if and only if $q$ is a generic of $G$.

\item Let $tp(a/A)  $ be generic of $G$; then so is $tp(a^{-1}/A)$.

\item There exists a generic type of $G$.

\item A type is left generic if and only if it is right generic.

\end{enumerate}
\end{fac}

%\begin{defi}
%Let $H$ be an $\infty$-definable subgroup of $G$ defined over $A$. We say that $H$ is the connected component of $G$, denoted by
%$G^0$,  if $H$ is of bounded index in $G$.
%\end{defi}
%
%Connected components of simple groups always exist, but there are important differences with the stable
%case, for example it is not always true that the connected component of a simple group has a 
%unique generic type. We also find some differences with the stable case  when we talk about stabilizers.
%But for the supersimple case or the quantifier-free $\omega$-stable case things are very close to the
%stable case.

%In order to define the stabilizer of a type $p$ in the case of simple theories we need to define $S(p)$.

\begin{defi}
Let $p \in S(A)$. $S(p)=\{g \in G: gp \cup p$ does not fork over $A\}$.
\end{defi}

Equivalently,  $a \in S(p)$ if and only if there are realizations $b,c$ of $p$, each one independent from $a$ over $A$, such that
$c=a \cdot b$. This implies in particular that $a^{-1} \in S(p)$.

%A consequence that we are working in a simple theory is that $S(p)$ is $\infty$-definable.
%Now suppose that ${\mathbf p}$ is a global type. Let $a,b \in S({\mathbf p})$ independent over $M$. Then
%$b^{-1} \cdot a \in S({\mathbf p})$.

\begin{defi}
Assume that $T$ satisfies the independence theorem over $A$. Let $ p$ be  type over $A$.
 We define the stabilizer of $p$ by $Stab(p)= S(p) S(p)$.
\end{defi}

\begin{teo}
$Stab( p)$ is an $\infty$-definable subgroup of $G$, and $p$ is a generic type of $G$
if and only if $Stab(p)$ is of bounded index in $G$.
\end{teo}
\begin{rem}
In the stable case, $S(p) =Stab(p)$.
\end{rem}

%Now suppose that $T$ is a supersimple theory
The following is a consequence of \ref{PSU3}.
%proved in \cite{wag}, chapter 5.

\begin{prop}\label{GIII5}
Let $G$ be a $\emptyset$-definable group, $H$ a $\emptyset$-definable subgroup of $G$ and let $ A=acl(A) $,

\begin{enumerate}
\item Let $p \in S(A)  $, then $ p$ is a generic of $G$ over $A$ if and only if $\s(G)=\s(p) $.
\item $\s(G)=\s(H)$ if and only if $[H:G]< \infty  $.
\item $\s(H)+\s(G/H) \leq SU(G) \leq \s(H) \oplus \s(G/H) . $ 
\end{enumerate}
\end{prop}

%The following theorem, due to Wagner (\cite{wagbase}) gives us conditions for 1-basedness, sability and stable embeddability
%for groups.

%\begin{teo} \label{wgr}
%Let $1 \longrightarrow G_1   \longrightarrow G_2    \longrightarrow G_3  \longrightarrow  1$ be a
%short exact sequence of definable groups in a simple theory. Then $G_2$ is
%1-based, stable, stably embedded if and only if $G_1$ and $G_3$ are 1-based, stable, stably embedded.
%\end{teo}

\section{Differential Fields}

In this section we introduce the basic notions of differential algebra and the theory of differentially closed fields. 
Even if some of the results hold in all characteristics we shall work in fields of characteristic zero.
We will work in the language ${\mathcal L}_D=\{0,1,+,-,\cdot,D \} $, where $D$ is a 1-ary function symbol.
%By abuse of language we will say $D$-definable and $D$-formula instead of ${\mathcal L}_D$-definable and 
%${\mathcal L}_D$-formula.

For algebraic results the references are to
\cite{kolda} and \cite{lang1}, for model-theoretic results see \cite{mark}, \cite{pp} and \cite{wood}.

\begin{defi}\label{P23}
A differential ring is a commutative ring $R$, together with an operator $D$ acting over $R$, such that, for every
 $x,y \in R $, we have:

\begin{enumerate}
\item $D(x+y)=Dx+Dy$
\item $D(xy)=xDy+yDx$

\end{enumerate}
If $R$ is a field, we say that $(R,D) $ is a differential field.
\end{defi}

\begin{defi}\label{P24}
  Let $ (R,D)$ be a differential ring. The differential polynomial ring over $R$ in $n$ indeterminates
is the ring $R[X]_D= R[X,DX,D^2X \cdots]$, where  $X=(X_1, \cdots ,X_n)$.
\end{defi} 
We extend $D$ to $R[X_1, \cdots, X_n]_D$ in the obvious way, it has then a natural structure of differential ring.

\begin{defi}\label{P25}
 Let $f \in K[X]_D $. The order of $f$, denoted by $ord(f) $, is the greatest integer $n$ such that $D^nX $ appears
in $f$ with non-zero coefficient. If there is no such $n$ we set $ord(f)=-1 $.
 
\end{defi}
Differential ideals play a key role in
the study of differential fields.

\begin{defi}\label{P26}
Let $(R,D)$ be a differential ring.
\begin{enumerate}
\item Let $I$ be an ideal of $R$. We say that $I$ is a differential ideal if it is closed under $D$.
\item If $A \subset R $, we denote by $(A)_D$ the smallest differential ideal  containing $A$, and by
$\sqrt{(A)_D} $ the smallest radical differential ideal containing $A$.
\end{enumerate} 
\end{defi}

\begin{rem}\label{P27} \hspace{20cm}
\begin{enumerate}
\item The radical of a differential ideal is a differential ideal.
\item If $I$ is a differential ideal of $R$, then   $R/I$  is a differential ring. 
\end{enumerate}
\end{rem}

We have a Finite Basis Theorem for radical differential ideals of the ring of differential polynomials over a differential
field.

\begin{teo}\label{P28}
Let $ (K,D)$ be a differential field, and $I$ a radical differential ideal of $K[X_1,\cdots,X_n]_D$. Then
 there is a finite subset $A$ of $I$ such that $I =\sqrt{(A)_D} $.
\end{teo}

This result fails for differential ideals which are not radical.

\begin{cor}\label{P29}
Let $(K,D)$ be a differential ring. Then $K[X_1,\cdots,X_n]_D $ satisfies the ascending chain condition on
 radical differential ideals.
\end{cor}

\begin{defi}\label{P210}
Let $(K,D) $ a differential ring and $L$ a differential subring of $K$.
\begin{enumerate}
\item Let $a \in K$. The differential ideal of $a$ over $L$ is $I_D(a/L)=\{f \in L[X]_D: f(a)=0 \}$.
If $I_D(a/L)={0} $ we say that $a$ is differentially transcendental (or $D$-transcendental) over $L$; 
otherwise we say that it is differentially algebraic.
\item Let $a \in K^n $. The differential ideal of $a$ over $L$ is
$I_D(a/L)=\{f \in L[X_1,\cdots,X_n]_D: f(a)=0 \}$.  
If $I_D(a/L)={0} $ we say that $a$ is differentially independent over $L$.
\end{enumerate}
\end{defi}

\begin{nota}\label{P211}
Let $ (K,D)$ a differential field. Let $A \subset K^n $, and $S \subset K[X]_D $, with $X=(X_1, \cdots,X_n) $.
\begin{enumerate}
\item $I_D(A)=\{f \in K[X]_D: f(x)=0 \; \; \forall x \in A   \}$.
\item $V_D(S)=\{x \in K^n: f(x)=0 \; \; \forall f \in S   \} $.
\end{enumerate}
\end{nota}

\begin{defi}\label{P212}
Let $(K,D)$ be a differential field. We define the $D$-topology of $K^n$ 
\textup{(}also called Kolchin topology or Zariski
 differential topology\textup{)}, as the topology with the sets of the form $V_D(I) $  as basic closed sets, 
where $I \subset K[X_1,\cdots,X_n]_D$ is a differential ideal.
\end{defi}

From \ref{P29} we deduce the following:

\begin{cor}\label{P213}
Let $(K,D)$ be a differential field. Then the $D$-topology of $K^n$ is Noetherian.
\end{cor}

A remarkable result in differential algebra is Kolchin's Irreducibility Theorem. 
Its proof can be found in \cite{mark2}, Chapter II, Appendix C.

\begin{prop}
Let $(K,D)$ be a difference field, and let $V$ be an irreducible algebraic variety defined over $K$.
Then $V$ is irreducible in the $D$-topology.
\end{prop}

\begin{defi}\label{P214}
Let $(K,D)  $ be a differential field. We say that $(K,D)$ is differentially closed if,
 for every  $f,g \in K[X]_D$, such that the order of $f$ is greater than the order of
 $g$,  there is $a \in K $ such that $f(a)=0  $ and $g(a) \neq 0 $.\\
We denote the theory of differentially closed fields by {\it DCF}.
\end{defi}

\begin{rem}\label{P215}
{\it DCF} is the model-companion of the theory of differential fields. As a consequence we have that
every differential field embeds in a model of {\it DCF}.
\end{rem}

\begin{teo}\label{P216}
Let $(K,D)$ a difference field. Then $(K,D) $ is a model of  {\it DCF}  if and only if
 it is existentially closed.
\end{teo}

\begin{teo}\label{P217}
The theory of differentially closed fields  is complete and $\omega$-stable; 
it  eliminates quantifiers and imaginaries.
\end{teo}

As {\it DCF} is $\omega$-stable, given a differential field $(K,D)$, {\it DCF} has  a prime model over $K$.
This prime model is unique up to $K$-isomorphism and is
called the differential closure of $(K,D)$.

For any differential field $(K,D)$ there is a distinguished definable subfield:
the field of constants ${\mathcal C}= \{x \in K: Dx=0 \}$.

\begin{teo}\label{P217a}
Let $(K,D)$ be a model of {\it DCF}. Then its field of constants ${\mathcal C}$ is 
an algebraically closed field and has no other definable structure, that is, any definable subset
of ${\mathcal C}^ n$, for $n \in {\mathbb N}$, is definable over ${\mathcal C}$ in 
the language of fields.
\end{teo}

\begin{fac}\label{P218}
Let $(K,D)$ be a model of  {\it DCF}, $A \subset K $. The definable closure of $A$, $dcl_{DCF}(A)$, 
is the smallest 
differential field containing $A$, 
and equals the field generated by  $(A)_D=(A,D(A), \cdots) $.\\ 
The algebraic closure of $A$, $acl_{DCF}(A)$, is $dcl_{DCF}(A)^{alg}=(A)_D^{alg}$, where $(A)^{alg}$ denotes
the field-theoretic algebraic closure of the field generated by $A$.

\end{fac}

For differentially closed fields we have a version of Hilbert's Nullstellensatz.

\begin{teo}\label{P219}
Let $(K,D)$ be a model of {\it DCF}. Let $I$ be a radical differential ideal of  $K[X_1, \cdots, X_n]_D $. 
Then $I_D(V_D(I))=I $.
\end{teo}

\begin{defi}\label{P220}
Let $ (K,D) $ be a differential field, and let $V$  be a variety in the
affine space of dimension $n$, let  $F(X)$ be a finite tuple of polynomials over $K$ 
generating $I(V)$, where $X=(X_1, \cdots , X_n)  $.
\begin{enumerate}

\item  We define the first prolongation of $V$, $\tau_1(V)  $ by the equations:

$$F(X)=0,  J_F(X)Y_1^t +F^D(X)=0  $$
where $Y_1 $ is an $n$-tuple, $F^D$ denotes the tuple of polynomials obtained by applying $D$
to the coefficients of each polynomial of $F$,
 and $J_F(X)  $ is the Jacobian matrix of $F$ (i.e. if $F =(F_1, \cdots , F_k) $ then
$ J_F(X)= (\partial F_i / \partial X_j)_{1 \leq i \leq k,1 \leq j \leq n} $).

\item For $m>1 $, we define the $m$-th prolongation of $V$ by induction on $m$:

Assume that $\tau_{m-1}(V) $ is defined by $F(X)=0$, 
$J_F(X)Y_1^t+F^D(X)=0, \cdots, $ 
$ J_F(X)Y_{m-1}^t+f_{m-1}(X,Y_1,\cdots,Y_{m-2})=0. $
Then $\tau_m(V)  $ is defined by:
$$(X,Y_1, \cdots,Y_{m-1}) \in \tau_{m-1}(V)$$ and
$$ J_F(X)Y_m^t+J_F^D(X)Y_{m-1}^t+J_{f_{m-1}}(X,Y_1,\cdots,Y_{m-2})(Y_1,\cdots,Y_{m-1})^t$$ $+
  f_{m-1}^D(X,Y_1,\cdots, Y_{m-2})=0 .$

\item Let $W \subset \tau_m(V)  $ be a variety. We say that $W$ is in normal form if,
 for every $i \in \{0, \cdots,m-1   \} $,
 whenever $G(X,Y_1, \cdots,Y_i) \in I(W) \cap K [X,Y_1,\cdots,Y_i ]  $ then 

$J_G(X,Y_1,\cdots,Y_i)(Y_1, \cdots, Y_{i+1})^t +G^D(X,Y_1,\cdots ,Y_i) \in I(W). $
\item Let $W \subset \tau_m(V)  $ be a variety in normal form. \\
A point $a$ (in some extension of $K$) is an 
$(m,D) $-generic of $W$ over $K$ if  $(a,Da, \cdots ,D^ma) $ is a generic of
 $W$ over $K$ and for every $i>m  $, 
$$\td(D^ia/K(a,\cdots,D^{i-1}a))=\td(D^ma/K(a, \cdots,D^{m-1}a) )  .$$
\end{enumerate}
\end{defi}

\begin{rem}\label{P222} \hspace{10cm}
\begin{enumerate}
\item There is a natural projection from $\tau_m(V)$ onto $\tau_{m-1}(V) $.
\item The map $\rho :\tau_{m+1}(V) \to \tau_1(\tau_m(V))$ defined by

$(x,u_1, \cdots, u_m) \mapsto ((x,u_1, \cdots,u_{m-1}),(u_1, \cdots,u_m ))$, 
defines an isomorphism between $\tau_{m+1}(V)$
and a Zariski-closed subset of $\tau_1(\tau_m(V)) $.
\end{enumerate}
\end{rem}

We give now a more geometric axiomatization of {\it DCF} due to Pierce and Pillay (\cite{pp}).
\begin{teo}\label{P224}
Let $(K,D)$ be a differential field. $K$ is differentially closed if and only if $K$ is an algebraically closed field and
for every  irreducible algebraic variety $V$, if $W$ is an irreducible algebraic subvariety of $\tau_1(V)$ 
such that the projection of $W$ onto $V$ is dominant, then
there is $a \in V(K)$ such that $(a,Da) \in W$.   
\end{teo}

We remark here that the axioms above hold just for the characteristic zero case. 
If $K$ is characteristic $p$ we must replace the condition of $K$ being algebraically closed by: 
$K$ is separably closed, and every constant is a $p$-th power. 

The following lemma (\cite{lang1}, chapter X), gives us a condition for extending the derivation of a 
differential field.

\begin{lem}\label{P221}
Let $(K,D) $ be a differential field and $\bar{a}=(a_i)_{i \in I} $ a (possibly infinite) tuple in some extension of $K$.
Let $\{F_j: j \in J \}$ be a set of generators of the ideal $I({\bar a}/K) \subset K[X_i : i \in I]$.

Let $(b_i)_{i \in I}$ be a tuple of $K({\bar a}) $ such that, for all $j \in J$
 $$\sum_{i \in I}\frac{\partial F_j}{\partial X_i}({\bar a})b_i+F_j^D({\bar a})=0.   $$ 
Then $D$ extends
to a unique derivation $D^*$ on $K(a_i: i\in I)$, such that $D^*a_i=b_i $ for all $i \in I$.
\end{lem}

\begin{lem}\label{P222.1}
Let $D:K \to \mathcal{U}$ be such that for all $a,b \in K$\\
 $(*) D(a+b)=Da+Db $\\
$(**) D(ab)=aDb+bDa $.\\
Let $a \in \mathcal{U}$.
\begin{enumerate}
\item If $a$ is transcendental over $K$, and $b \in \mathcal{U}$,
then there is $D_1: K(a) \to \mathcal{U}$ extending $D$ and satisfying $(*)$ and $(**)$ such that $D_1a=b$.
\item If $a$ is algebraic over $K$, then there is a unique extension $D_1$ of
$D$ to $K(a)$ satisfying $(*)$ and $(**)$.
\end{enumerate}
\end{lem}

{\it Proof}:\\

(1) For $f(a) \in K(a)$, set $D_1(f(a))=f'(a)b+f^D(a).$
 Since $a$ is transcendental over $K$, one checks easily that $(*)$ and $(**)$ hold.

(2) Let $f(X)=\sum_{i=0}^na_iX^i$ be the monic minimal polynomial of $a$ over $K$. We define $D_1a = -f'(a)^{-1}f^D(a)$.
 Every element of $K(a)$ can be written $\sum_{i=0}^{n-1}b_ia^i$ where the $b_i$ are in $K$. We then set 
$$D_1(\sum_{i=0}^{n-1}b_ia^i)=-\sum_{i=0}^{n-1}(D(b_i)a^i+ib_iD_1(a)a^{i-1}).$$ 
Clearly, $D_1$ satisfies $(*)$, and to check that is satisfies $(**)$, it suffices to show that $D_1(a^n)=nD_1(a)a^{n-1}$. 
Since $a^n=-\sum_{i=0}^{n-1}a_ia^i$, we have 
$$D_1(a^n)=-\sum_{i=0}^{n-1}(D(a_i)a^i+ia_iD_1(a)a^{i-1})$$ 
$$=-f^D(a)-(f'(a)-na^{n-1})D_1(a) $$
$$=na^{n-1}D_1(a).$$ 
%Let $f(X) \in K[X]$ monic such that $f(a)=0$ and $f'(a)\neq 0$. Then $f'(a)$ is invertible. We define
%$D_1a=-f'(a)^{-1}\cdot f^D(a)$.\\
$\Box$

%The following proposition hepls usto clear out what being in normal form means.

\begin{prop}\label{P223} 
Let  $({\mathcal U},D) $ be a saturated model of {\it DCF}, let $K=acl(K) \subset {\mathcal U}$,
let $V$ an irreducible affine variety, and $W$  
a subvariety of $ \tau_m(V) $ both defined over $K$.
If $W$ is in normal form, then $W$ has an $(m,D) $-generic in ${\mathcal U}$.
\end{prop}

{\it Proof} :\\

We will construct a differential field containing $K$, and which contains an $(m,D)$-generic of $W$. 
We work in some large algebraically 
closed field containing $K$, and choose a generic $(a,b_1,\cdots,b_m)$ 
of $W$ over $K$. Since $W$ is in normal form, by \ref{P221} we define
$D :K(a,b_1, \cdots,b_{m-1}) \to K(a,b_1, \cdots,b_{m})$ by setting $Da=b_1$ and $Db_i=b_{i+1}$,
and so that it satisfies $(*)$ and $(**)$.
Let $v_m \subset D^ma$ be a transcendence basis of $b_m$ over $K(a,b_1, \cdots,b_{m-1})$ 
and let $(v_{n})_{n>m}$ be a set of tuples of the same length as $v_m$ such that for all 
$n>m$ the elements of $v_n$ are algebraically independent over
$K(a,b_1, \cdots,b_{m},v_{m+1}, \cdots, v_{n-1})$.
By \ref{P222.1}, the map $D$ on $K(a,b_1,\ldots,b_{m-1})$ extends (uniquely) 
to a map $D_1$ defined on $L=K(a,b_1,\cdots,b_m,v_n)_{n>m}$ which sends 
$v_n$ to $v_{n+1}$ for $n\geq m$ and satisfies $(*)$ and $(**)$. Then 
$D_1$ is a derivation of $L$, and $a$ is an $(m,D)$-generic of $W$.\\
$\Box$

\begin{cor}
Let $(K,D)$ be a differentially closed field.
Let $V$ a variety, and $W$  
a subvariety of $ \tau_m(V) $ both defined over $K$. 
Then $W$ is in normal form if and only if $\{(x,Dx, \cdots,D^mx): x\in V   \} \cap W  $ is Zariski dense in $W$.
In particular  $\{(x,Dx, \cdots,D^mx): x\in V   \}  $ is Zariski dense in $\tau_m(V)  $
and $dim(\tau_m(V))=(m+1)dim(V)$.
%\begin{enumerate}
%\item $\{(x,Dx, \cdots,D^mx): x\in V   \}  $ is Zariski dense in $\tau_m(V)  $.
%\item If $V$ is smooth then $\tau_m(V)$ is irreducible and smooth.
%\item  $W$ is in normal form if and only if $\{(x,Dx, \cdots,D^mx): x\in V   \} \cap W  $ is Zariski dense in $W$.
%\item 
\end{cor}

\begin{rem}\label{P223a}

Let $(K,D)$ be a differentially closed field and $V$ a smooth variety in the affine space of dimension $n$ defined over 
$E=acl_D(E) \subset K$. If $W \subset \tau_m(V)$ is a variety in normal form then all
$(m,D)$-generics of $W$ have the same type over $E$.
%\begin{enumerate}
%\item $V$ is smooth if and only if $\tau_m(V)$ is a variety.
%\item Let $W \subset \tau_m(V)$ be a variety in normal form. If $V$ is smooth then all the $(m,D)$-generics 
%of $W$ have the same type.
%\end{enumerate}
\end{rem}

We introduced varieties in normal form to bypass some difficulties concerning differential ideals.

Let $W\subset \tau_m(V)$ be a variety in normal form, and let $I\subset K[X,Y_1,\ldots,Y_m]$ its defining ideal, 
which is a prime ideal. Let 
$\varphi: K[X,Y_1,\ldots,Y_m]\to K[X]_D$ be the $K[X]$-algebra embedding sending $Y_i$ to $D^iX$ for $i=1,\ldots,m$, 
and let $J$ be the differential ideal generated by $\varphi(I)$. 
%Then $J$ is not necessarily a prime differential ideal.
 %However, we have shown that if $W$ is in normal form, then in some differential field extension of $K$ 
%it has an $(m,D)$-generic, and furthermore, any two $(m,D)$-generics of $W$ over $K$ satisfy the same
% differential equations over $K$. 

Let  $L$ be a sufficiently saturated differentially closed field containing $K$, 
and consider the set $\bar W$ defined by $J$. The set 
$\bar W$ may not be irreducible for the Kolchin topology. However, it 
will have an irreducible component $W_0$ with the following property: 
$I_D(W_0)$ is the unique prime differential ideal  containing 
$\varphi(I)$ and whose intersection with $K[X,DX,\ldots,D^mX]$ equals 
$\varphi(I)$. All points in the other irreducible components of $W$ 
will satisfy some additional equations of order $m$. Furthermore, if 
$a$ is a generic of $W_0$ over $K$  in the sense of the Kolchin topology 
(i.e., $ W_0$ is the smallest Kolchin closed 
set defined over $K$ which contains $a$), then $a$ will be an 
$(m,D)$-generic of $W$ and conversely.

Thus to each variety in normal form defined over $K$ is associated in a 
canonical way an irreducible Kolchin closed set defined over $K$ (and 
therefore a unique complete type over $K$). The condition of a variety 
being in normal form is clearly expressible by first-order formulas on 
the coefficients of the defining polynomials, while it is not as 
immediate that the property of differential polynomials to generate a 
prime differential ideal is elementary in their coefficients.

\begin{lem}\label{P225}
Let $(L,D)$ be a differential field, and let $K$ be a differential subfield of $L$.
Let $a$ be a tuple of $L$, let $v \subset a $.
If the elements of $D^{m+1}v $ are algebraically independent over $K(a,\cdots,D^ma)$, then for all
 $i \in \{0,\cdots, m\} $,
the elements of $D^iv $ are algebraically independent over $K(a,\cdots,D^{i-1}a)$ \textup{(}or over $K$ if $i=0$\textup{)}. 
\end{lem}
 
{\it Proof}:\\

By reverse induction on $i$ it is enough to prove that the elements of  $D^{m}v $ are 
algebraically independent over  $K(a,\cdots,D^{m-1}a)$.

If the elements of $D^{m}v $ are algebraically dependent over $K(a,\cdots,D^{m-1}a)$, then there is a non zero polynomial
$P(X) \in K(a,\cdots,D^{m-1}a)[X]$ which is irreducible and vanishes at $D^{m}v $. 
Thus $J_P(D^mv) (D^{m+1}v)^t+P^D(D^{m}v)=0 $, and, as  $P$ is irreducible and we work in characteristic zero,
$J_P(D^{m}v)\neq 0$.
Then, since $P^D(D^mv) \in K(a,\cdots,D^{m}a)$, $D^{m+1}v$ satisfies a non-trivial equation over
$K(a,\cdots,D^{m}a) $ which contradicts our assumption. Hence the elements of $D^{m}v $
are algebraically independent over $K(a,\cdots,D^{m-1}a)$.\\
$\Box$

\begin{cor}\label{P226}
Let $K $ be a differential subfield of $(L,D) $, let $a$ be a tuple of $L$, let
$d_{n+1}=\td(K(a, \cdots,D^{n+1}a)/K(a,\cdots,D^na)) $.
Then $(d_n)_{n \in {\mathbb N}}$ is a decreasing sequence.
\end{cor}

{\it Proof} :\\

Let $n \geq 0 $. Then $$d_{n+1}=\td(D^{n+1}a/K(a,\cdots,D^{n}a))
=\td(D^{n+1}a/K(a,\cdots,D^{n-1}a)(D^{n}a)),$$
and the latter, by \ref{P225} is less than or equal  to
$\td( D^{n}a/K(a,\cdots,D^{n}a))=d_n .$\\
$\Box$

\begin{rem}\label{P227}
Since $d_n $ is a decreasing sequence in ${\mathbb N} \cup \{ \infty\}$, there is $M \in {\mathbb N} $ such that $d_n=d_M $ for all
$n \geq M $. Thus $a$ is an $(M,D)$-generic of the locus of $(a,Da,\cdots,D^Ma) $ over $K$.
\end{rem}

\begin{lem}\label{P228}
Let $(K,D ) $ a differential field and $(L,D)  $ an extension. 
Let $b$ be a tuple of $L$. 

Assume that, for $i>1  $, $\td(D^ib/K(b,Db, \cdots, D^{i-1}b))=\td(Db/K(b)).$ 
Let $a \in K(b)  $ such that, for some $n > \td(b/K(a) ) $ we have 
$\td( (Da, \cdots, D^na )/K(a))=n \td(Da/K(a)  )$.

Then $\td( (Da, \cdots, D^ia )/K(a))=i \td(Da/K(a)) $ for every $i>n$.

\end{lem}

{\it Proof} :\\

We proceed by induction on $d=\td(Da/K(a))$.
It is clear for $d=0$.

Let $v \subset Da $ be a transcendence basis for $Da$ over $K(a)$.
We can rewrite the hypothesis of the theorem as:
 the elements of $\{D^jv: 0 \leq j < n \}$ are algebraically independent over $K(a)$.
And we must prove that the elements of $\{D^jv: j \in {\mathbb N} \} $
 are algebraically independent over $K(a)$.

Since $\td(v, \cdots, D^{n-1}v/K(a)) > \td(b/K(a)) $,
 $K(a,v,\cdots,D^{n-1}v) \not\subset K(b) $.
 Let $i$ be the smallest integer such that
$D^iv \not\subset K(b) $; then $D^iv \subset K(b,Db) $.
 Let $w_0 \subset D^iv $ be a transcendence basis for $D^iv$ over $K(b)$,
 and let $w \supset w_0$ be a transcendence basis for $Db $ over $K(b)$.
Let $v_0 \subset v $ be such that $D^iv_0=w_0$, let $a_0 \subset a $
be such that $Da_0=v_0 $ and let $v_1=v \setminus v_0$.

Our hypothesis implies that $\{v,\cdots,D^{n-1}v\}$ is a transcendence basis for
$K(a, \cdots,D^na)$ over $K(a)$; so $\{v_1,\cdots,D^{n-1}v_1\}$
is a transcendence basis for $K(a, \cdots,D^na)$ over
$K(a,v_0, \cdots,D^{n-1}v_0)$.
Both fields are contained in
$K(b,w,\cdots,D^{n-i-1}w)$, thus the elements of  $\{ D^jw :j \geq n-i\}$
are algebraically independent over $K(a,\cdots,D^na)$, so 
 $\td(Da,\cdots,D^na/(K(a_0)_D(a))=\td(Da,\cdots,D^na/K(a,Da_0,\cdots,D^na_0))$
$=n \td(Da/K(a,Da_0))$.

By induction hypothesis applied to $a$ and $K(a_0)_{D} $, the elements of $\{D^jv_1:j \geq 0 \} $ 
are algebraically independent over
 $K(a_0)_D(a) $; thus  the elements of $\{D^jv: j \geq 0 \} $ are algebraically independent over $K(a)$,
since $\td(a_0,\cdots,D^ia_0/K)=(i+1)\td(a_0/K)$ for all $i>0.$\\
$\Box$

\begin{cor}\label{P229}
Let  $V,W,V_1 \subset \tau_1(V),W_1 \subset \tau_1(W)$ be irreducible varieties defined over a 
differentially closed field $K$. Let $f:V \to W $ be a rational map.
Then the following property is expressible in the first order language ${\mathcal L}_D $ 
with the parameters needed to define $f,V,W,V_1,W_1 $:  

$V $ and $W$ are smooth varieties, $V_1  $ and $W_1$ are varieties in normal form,
and a $(1,D)$-generic of $V_1  $ is sent by $f$ to a $(1,D ) $-generic of $W_1$.
\end{cor}

{\it Proof} :\\

By the results in \cite{vdd},  we know that we can express in ${\mathcal L}_D $ that $V $ is a smooth variety, 
$V_1 \subset \tau_1(V)$  is a variety in normal form,  and that a rational map between two varieties sends 
generic points onto generic points. 
Using the characterization of varieties in normal form given in \ref{P220}, for every $m \geq 0$
we can construct subvarieties $V_m \subset \tau_m(V)$ and $W_m \subset \tau_m(W)$ such that
the $(m,D)$-generics of $V_m$ are exactely the $(1,D)$-generics of $V_1$ and similarly for $W_m$ and
$W_1$. By \ref{P228} it suffices to say that the projection $W_m \to V_m$ is dominant,
where $m=dim(W)-dim(V)+1$.\\
%By \ref{P228} there is $N \in \na$ such that for every $m \geq N$ we
%can construct  subvarieties $V_m$ of $\tau_m(V)$ and $W_m$ of $\tau_m(W)$ such that the $(m,D)$-generics of $V_m$ are
% the $(1,D)$-generics of $V$ and  similarily for $W_m$ and $W_1$. By \ref{P228}, it then suffices to say that the 
%projection $W_m \to V_m$ is dominant, where $m=dim W_1-dimV_1+1$\\
%there is $N \in {\mathbb N} $  such that, 
%if  $a$ is an ($m,D$)-generic point of $V_1$, and for all $n<N $ 
%$\td(D^mf(a)/K(f(a), \cdots, D^{m-1}f(a)))$ $=\td(D^Nf(a)/K(f(a), \cdots, D^{N-1}f(a))) $,
%then $f(a)$ is an $(m,D)$-generic of $W_1$.\\
$\Box$\\

As {\it DCF} is $\omega$-stable, there is a notion of independence: Let $(K,D)$ be a differentially closed field,
let $A,B,C \subset K$. We say that $A$ is independent from $B$ over $C$ if $acl_D(AC)$ and $acl_D(BC)$ are
linearly disjoint over $acl_D(C)$.\\

An important result in the theory of differential fields is  Zilber's dichotomy.

\begin{teo}\label{DCFdich}
Let $({\mathcal U},D)$ be a differentially closed field and let $K \subset {\mathcal U}$. 
Let $p \in S(K)$ be a stationary type of $U$-rank 1. Then $p$ is either 1-based or non-orthogonal to the
field of constants. 
\end{teo}

\section{Difference Fields}

In this section we mention  basic facts and definitions about fields of characteristic 0 with an automorphism. 
As in the preceding section, some of the results that we recall here hold in any characteristic.
For the model-theoretical statements we shall work in the language ${\mathcal L}_{\sigma}=\{0,1,+,-,\cdot, \sigma   \} $, where
 $\sigma $ is a 1-ary function symbol.\\
For the proofs of the results in difference algebra the reader may consult \cite{cohn2}, 
for model-theoretic results we refer to \cite{salinas}.
\begin{defi}\label{P31}
A difference field is a field $K $ together with a field endomorphism $\sigma $. If $\sigma$ is an automorphism 
we say that $(K,\sigma)$ is an inversive difference field.
\end{defi}

\begin{fac}\label{P32}
Every difference field $K$ embeds into a smallest inversive difference field, and this field is unique up to $K$ isomorphism.\
\end{fac}
{\bf From now on we  assume all difference fields to be inversive.}

\begin{defi}\label{P33}
Let $(K,\sigma)$ be a difference field. The difference polynomial ring over $K$ in $n$ indeterminates is
the ring $K[X]_{\sigma}=K[X,\sigma(X),\si^2(X),\cdots] $, where $X=(X_1, \cdots, X_n)$.
\end{defi}

\begin{rem}\label{P34}
We  extend $\sigma$ to $K[X_1,\cdots,X_n]_{\sigma} $ in the obvious way. This map is injective but not surjective.
\end{rem}

%We could have an extension of $\si$ to an automorphism of $K[X]_{\si}$ if we have included the negative powers of 
%$\si(X)$ on the definition of $K[X]_{\si}$.

\begin{defi}\label{P35}
Let $(K,\sigma)$ a difference field, $X=(X_1,\cdots,X_n) $. Let $I $ be an ideal of $K[X]_{\sigma} $.
We say that $I$ is a reflexive $\sigma$-ideal if for every $f \in K[X]_{\sigma}$, $f \in I$ if and only if
$\sigma(f) \in I$.
If, in addition, for every $f \in K[X]_{\sigma} $ and for every $m \in \na$ $f^m \sigma(f)^n \in I $ implies $f \in I $, we say that $I$ 
is a perfect $\sigma$-ideal.

A prime ideal which is a perfect $\sigma$-ideal is called a prime $\sigma$-ideal.
\end{defi}

\begin{rem}\label{P36}
If $I$ is a $\sigma$-ideal, then $\sigma$ induces an endomorphism on $K[X]_{\sigma} / I $.
\end{rem}

\begin{defi}\label{P37}
Let $(K,\sigma)$ be a difference field and $F$ a difference subfield of $K$
 \textup{(}that is, $F$ is a  subfield of $K$ and the restriction to $F$ of $\sigma$ is an automorphism of $F$\textup{)}. Let $a \in K$. 
We say that $a$ is transformally transcendental \textup{(}or $\si$-transcendental\textup{)} over $F$ if 
$I_{\sigma}(a/F)=\{f \in K[X]_{\sigma}: f(a)=0 \}=(0) $.
 Otherwise, we say that $a$ is transformally algebraic over $F$.
\end{defi}

\begin{nota}\label{P38}
Let $(K,\sigma)$ be a difference field. Let $A \subset K^n $, and $S \subset K[X]_{\sigma}$ with 
$X=(X_1,\cdots,X_n) $.
\begin{enumerate}
\item $I_{\sigma}(A)=\{f \in K[X]_{\sigma}: f(x)=0 \;\; \forall \;\; x \in A   \}$.
\item $V_{\sigma}(S)=\{x \in K^n: f(x)=0 \;\; \forall \;\; f \in S \}$. 
\end{enumerate}
\end{nota}

\begin{rem}\label{P39}
$I_{\sigma}(A) $ is a perfect $\sigma$-ideal.
\end{rem}

\begin{prop}\label{P310}
$K[X_1,\cdots,X_n]_{\sigma} $ satisfies the ascending chain condition on perfect $\sigma$-ideals.
\end{prop}

\begin{defi}\label{P311}
Let $(K,\sigma)$ a difference field. We define the $\sigma$-topology of $K^n$ as the topology
with the sets of the form  $V_{\sigma}(S)$ as basic closed sets.
\end{defi}

\begin{cor}\label{P312}
Let $(K,\sigma)$ be a difference field. Then the $\sigma$-topology of $K^n $ is Noetherian.
\end{cor}

\begin{nota}\label{P313}
Let $(K, \sigma) $ be a difference field, and let $V$ be an algebraic set defined over $K$. 
By $V^{\sigma} $ we denote the algebraic set obtained by applying $\si$ to the coefficients of the polynomials
defining $V$.
\end{nota}

\begin{teo}\label{P314}
The theory of difference fields has a model-companion, that we shall denote by {\it ACFA}.
It is described as follows. \\
$(K,\sigma) \models$  {\it ACFA} if and only if:
\begin{enumerate}

\item $K$ is an algebraically closed field.

\item $(K,\sigma)$ is a (inversive) difference field.

\item For every irreducible algebraic variety $V$, if $W$ is an irreducible algebraic subvariety of
$V \times V^{\sigma}$, such that the projections from $W$ on $V$ and $V^{\sigma}$ are dominant,
then there is $a \in V(K) $ such that $(a,\sigma(a)) \in W$.
\end{enumerate}
\end{teo}

%\begin{prop}\label{P315}
As a direct consequence of the definition of model-companion we have that 
every difference field embeds in a model of {\it ACFA}.
%\end{prop}

\begin{nota}\label{P316}
Let $(K,\sigma) $ a difference field, $A \subset K $. We denote by $(A)_{\si} $ the smallest difference field containing $(A) $, and by
$acl_{\sigma}(A) $ the field-theoretic algebraic closure of $(A)_{\si} $.
\end{nota}

\begin{prop}\label{P317}
Let $(K,\sigma) $ be a model of {\it ACFA}, and let $A \subset K $. Then $acl_{ACFA}(A)=acl_{\sigma}(A) $.
\end{prop}

\begin{prop}\label{P318}
{\it ACFA} is  model-complete and  eliminates imaginaries.
\end{prop}

As in differential fields, we define independence in {\it ACFA} using independence in field theory, that is, if
$(K,\si)$ is a model of {\it ACFA}, and $A,B,C \subset K$, we say that $A$ is independent of $B$ over $C$ if
$acl_{\si}(AC)$ and $acl_{\si}(BC)$ are linearly disjoint from $acl_{\si}(C)$. 

\begin{prop}\label{P319}
All completions of {\it ACFA} are supersimple, independence coincides with non-forking.
\end{prop}

For models of {\it ACFA} we have a version of Hilbert's Nullstellensatz.

\begin{teo}\label{P320}
Let $(K,\sigma)$ be a model of {\it ACFA} . Let $I $ be a perfect $\sigma$-ideal of $K[X_1, \cdots,X_n]_{\sigma} $. 
Then $I_{\sigma}(V_{\sigma}(I))=I $.
\end{teo}

As {\it ACFA} is supersimple every type is ranked by the $\s$-rank, and an element of a model of {\it ACFA}
is $\si$-transcendental if and only if it has $\s$-rank $\omega$. 

If $(K,\si)$ is a model of {\it ACFA}, there is a distinguished definable subfield of $K$: the fixed field
$Fix \si=\{x\in K:\si(x)=x \}$.

\begin{prop}\label{psff}
$Fix \si$ is a pseudo-finite field. That is:
\begin{enumerate}
\item $Fix \si$ is perfect.
\item $Gal((Fix \si)^{alg}/Fix \si)={\hat {\mathbb Z}} $.
\item $Fix \si$ is pseudo-algebraically closed (PAC).
\end{enumerate}
\end{prop}

One of the consequences of this is that the definable field $Fix \si^n$, for $n \in \na$, is the unique extension
of $Fix \si$ of degree $n$. We have also that $Fix \si$ is the unique definable subfield of $K$ of $\s$-rank 1, and
 that the $\s$-rank of $Fix \si^n$ is $n$.

Pseudo-finite fields are infinite models of the theory of finite fields. 
The theory of  pseudo-finite fields is studied  in \cite{axff} and \cite{hrushpff}.

\begin{prop}\label{sistabem}
$Fix \si$ is stably embedded; that is, every definable subset of $(Fix\si)^n$ is definable with parameters in $Fix\si$.
Moreover it is definable in the pure language of fields. 
\end{prop}

\begin{prop}
Let $(K,\si)$ be a model of {\it ACFA}. Then, for all $n \in \na$, $(K, \si^n)$ is a model of {\it ACFA}
\end{prop}

As {\it DCF}, {\it ACFA} satisfies a version of Zilber's dichotomy.

\begin{teo}
Let $({\mathcal U},\si)$ be a saturated model of {\it ACFA} and let $K \subset {\mathcal U}$. 
Let $p \in S(K)$  be a type of  $\s$-rank 1. Then $p$ is either 1-based or non-orthogonal to the
fixed field. 
\end{teo}

\cleardoublepage
\chapter{Difference-Differential Fields}
\label{chap:dcfa}

This chapter is devoted to the study of difference-differential fields
of characteristic zero, first we shall give the algebraic properties of
such fields. In section two we give a proof of Hrushovski's theorem about
 the existence of a model-companion for the theory of difference-differential
 fields of characteristic zero, which we call {\it DCFA}. The original approach was to emulate 
the case of difference fields, the problem is that we cannot quantify on
 differential varieties but we can get around this using prolongations of
 differential varieties. We give also some properties of {\it DCFA}.
Next we mention some properties of the fixed field and the field of constants of a model of {\it DCFA}.
Finally we talk about forking and the $\s$-rank.

\section{Difference-Differential Algebra}

First we mention some facts concerning systems of ideals, see \cite{ide} for 
 details.

\begin{defi}\label{DCFA46a}
Let $R$ be a commutative ring, and $\mathcal{C} $ a set of ideals of $R$.
\begin{enumerate}
\item  We say that $\mathcal{C}$ is a conservative system of ideals if:
  \begin{enumerate}
  \item For every $I \subset  \mathcal{C}$, $\bigcap I \in \mathcal{C}  $.
  %\item For every $I \subset \mathcal{C}$ totally ordered by  inclusion, 
 $\bigcup I \in \mathcal{C}  $. 
  \end{enumerate} 
\item Let $\mathcal{C}$ be a conservative system of ideals. We say that $\mathcal{C}$ is divisible if for 
$I \in \mathcal{C}$
 and $a \in R$ we have $(I:a) \in \mathcal{C}$.
\item Let ${\mathcal C} $ be a divisible conservative system of ideals. We say that $\mathcal C$ 
is perfect if all its members are radical ideals.
\end{enumerate}
\end{defi}

%\begin{prop}\label{id1}
%Let $C$ be perfect, and let $I \in C$. Then $A$ is the intersection of prime ideals of $C$; if $C$
%is Noetherian this intersection can be taken finite.
%\end{prop}

\begin{teo}\label{fipri}\textup{(\cite{ide}, section 2, Theorem I)}\
Let $R$ be a commutative ring and $\mathcal{C} $ a perfect system of ideals,
let $I \in \mathcal{C} $. Then $I$ is an 
intersection of prime ideals of $\mathcal{C}$. If $R$ is Noetherian this intersection can be taken 
to be a finite intersection.
\end{teo}

%\begin{teo}\label{id2}
%Let $R$ be a ring and $S$ a subring of $R$ such that $R$ is finitely generated over $S$.
%Let $\mathcal{C}$ be a perfect set of ideals of $R$. If $C/S$ is Noetherian
%then $C$ is Noetherian.
%\end{teo}

%We expose now the basic ideas of difference-differential rings and fields.

\begin{defi}\label{DCFA41}

A difference-differential ring is a ring $R$ together with a finite set of
derivations $\Delta= \{D_1, \cdots, D_m   \}$ and  a finite set of automorphisms $A=\{\si_1,\cdots, \si_n \}$
such that all pairs in $\Delta \cup A$ commute.

If $R$ is a field we say that $(R,A,\Delta) $ is a difference-differential field.
\end{defi}

Let us denote by $\Theta$ the set of formal (commuting) products of elements of $A \cup \Delta$.
Let $(R,A, \Delta)$ be a difference-differential ring. An ideal $I$ of $R$ is said to be 
a difference-differential ideal if it is closed under the operators of $\Theta$.

The set of difference-differential ideals of $R$ is conservative, but not necessarily divisible.

\begin{defi}\label{id3}
Let $(R,A,\Delta)$ be a difference-differential ring. Let $I$ be a difference-differential ideal of $R$.
We say that $I$ is a perfect ideal if:
\begin{enumerate}
\item $I$ is radical.
%\item For every $a \in R$ and $\si \in A$, if $\si(a) \in I$ then $a \in I$.
\item For every $a \in R$ and $\si \in A$, if $a\si(a) \in I$ then $a \in I$.
\end{enumerate}
\end{defi}

A theorem from \cite{ide} (pp.798-799) tells us that the set of  perfect difference-differential ideals is
perfect (in the sense of \ref{DCFA46a}), and that it contains any perfect set of ideals.

\begin{teo}\label{DCFA47}\textup{(\cite{cohn}, section 5, Corollary I)}
Let $(S, A, \Delta) $ be a difference-differential ring which contains $\mathbb Q $ 
and is such that the set of perfect
difference-differential ideals of $S$  satisfies the ascending chain condition.
Let $(R,A',\Delta' )$ be a difference-differential ring finitely generated over $S$ as a difference-differential ring. 
Then the set of perfect difference-differential ideals of $R$ satisfies the ascending chain condition.
\end{teo}

From now on we will assume that we work in difference-differential rings with one derivation and
one automorphism. We will often write $(\si,D)$ instead of difference-differential (for example
$(\si,D)$-ideal in place of difference-differential ideal).

\begin{defi}\label{DCFA42}
Let $(R,\sigma,D)$ be a difference-differential ring. The ring of difference-differential polynomials  in $n$
indeterminates over $R$ is the ring $R[X]_{\sigma,D}$ of polynomials in
the variables $\sigma^i( D^jX ) $ for
$i,j  \in {\mathbb N}  $, where $X=(X_1,\cdots,X_n)$.
\end{defi}

As in the differential and difference cases,
 we can extend  $D$ to a derivation on $R[X]_{\sigma,D} $ and $\sigma$ to an
 endomorphism of  $R[X]_{\sigma,D} $
which commutes with $D$.

\begin{rem}\label{DCFA44}
Let $(R,\sigma,D)$ be a difference-differential ring. Let $I$ be an ideal of $R$.
\begin{enumerate}
\item  $I$ is a $(\sigma,D)$-ideal if it is a differential ideal and a $ \sigma $-ideal, in the sense of
\ref{P26} and \ref{P35}.
\item  $I$ is a perfect $(\sigma,D)$-ideal if it is a $(\sigma,D)$-ideal which is perfect as a
$\sigma $-ideal.
\end{enumerate}
\end{rem}

\begin{nota}\label{DCFA45}
  Let $(K,\sigma,D) $ a difference-differential field, $ S \subset K [ {\bar X }   ] _{\sigma,D}  $,
 $A \subset K^n  $; and let $E$ be a difference-differential subfield of  $ K $, $a \in K^n$.

\begin{enumerate}
\item $V_{\sigma,D}(S)=  \{{\bar x} \in K^n: \forall \; f({\bar X}) \in S \; f({\bar x})=0  \}  $.

\item $I_{\sigma,D}(A)= \{ f({\bar X}) \in  K[ {\bar X }] _{\sigma,D}:  \forall \;  {\bar x} \in A \; f({\bar x})=0 \} $.

\item $ I_{\sigma,D}(a/E)=   \{   f({\bar X}) \in   E [ {\bar X }   ] _{\sigma,D}:    f({\bar a})=0 \}       $.

\end{enumerate}

We define the $(\sigma,D)$-topology of $K^n$ to be the topology with the sets of the form $V_{\sigma,D}(S) $
 as a basis of closed sets.

\end{nota}

\begin{rem}\label{DCFA46}
Let $(K,\sigma,D) $ be a difference-differential field, $A \subset K^n  $. Then $I_{\sigma,D}(a/E) $ is a perfect
 $(\sigma,D)$-ideal.
\end{rem}

%\begin{rem}
%Let $(R,\sigma,D) $ be a difference-differential ring. In \cite{ide} the author prove that
%the set of perfect $(\sigma,D)$-ideals of $R$ is conservative, 
%and every perfect $(\sigma ,D) $ ideal 
%is perfect in the sense of \ref{DCFA46a}.
%\end{rem}

\begin{cor}\label{DCFA48}
Let $(K,\sigma,D)$ be a difference-differential field. Then, by \ref{DCFA47} the $(\sigma,D)$-topology 
of $K^n$ is Noetherian.
\end{cor}

%\begin{cor}\label{DCFAfinprid}
%Let $(K,\sigma,D)$ be a difference-differential field and let $I$ be a perfect $(\sigma,D)$-ideal of 
%$K[X_1,\cdots,X_n]_{(\sigma,D)}$. Then $I$ is a finite intersection of prime perfect $(\sigma,D )$-ideals.
%\end{cor}

\begin{cor}\label{DCFAfingen}
Let $(K,\sigma,D)$ be a difference-differential field and let $I$ be a perfect $(\sigma,D)$-ideal of 
$K[X_1,\cdots,X_n]_{\sigma,D}$. Then $I$, as a $(\si,D)$-perfect ideal,
 is generated by a finite number of $(\si,D)$-polynomials.
\end{cor}

\begin{cor}\label{DCFAfinprid}
Let $(K,\sigma,D)$ be a difference-differential field and let $I$ be a perfect $(\sigma,D)$-ideal of 
$K[X_1,\cdots,X_n]_{\sigma,D}$. Then $I$ is a finite intersection of prime perfect $(\sigma,D )$-ideals.
\end{cor}

\section{The Model-Companion}

We begin this section with Hrushovski's theorem on the existence of a model-companion
for the theory of difference-differential fields of characteristic zero. 
In the axiom scheme that we give we try to emulate somewhat the axioms for
{\it ACFA}.% the problem is that we cannot quantify on differential varieties, 
%but we get around this thanks to varieties in normal form introduced on 
%Chapter \ref{chap:prelim}.

\begin{teo}\label{DCFA49} \textup{(}Hrushovski\textup{)}

The model companion of the theory of difference-differential fields exists. We denote it by {\it DCFA} and it is described as follows:

$(K,D,\sigma)$ is a model of  {\it DCFA} if
\begin{enumerate}
\item $(K,D)$ is a differentially closed field.
\item $\sigma$  is an automorphism of $(K,D)$.
\item If $U,V,W$ are varieties such that:
\begin{enumerate}
\item $U \subset V \times V^{\sigma}  $ projects generically onto $V$ and
 $ V^{\sigma}  $.
\item $W \subset \tau_1 (U) $ projects generically onto $U$.
\item $\pi_1(W)^{\sigma}=\pi_2(W)$ \textup{(}we identify $\tau_1(V \times V^{\sigma})$ 
with $\tau_1(V)\times \tau_1(V)^{\si}$ and let $\pi_1:\tau_1(V \times V^{\si}) \to \tau_1(V)$
and $\pi_2:\tau_1(V \times V^{\si}) \to \tau_1(V)^{\si}$ be the natural projections\textup{)}. 
\item A $(1,D) $-generic point  of $W $ projects  
onto a $(1,D)   $-generic point of $\pi_1(W)$ and onto a (1,D)-generic point of $\pi_2(W) $.
\end{enumerate}
Then there is a tuple $a \in V(K) $, such that 
$(a,Da,\si(a),\si(Da)) \in W$.

\end{enumerate}
\end{teo}

{\it Proof}:\\ 

By \ref{P229}, these are first order properties. 
First we prove that any difference-differential field embeds in a model of {\it DCFA}.
By  quantifier elimination in {\it DCF}
any difference-differential field embeds into a model of (1) and (2).
       By the usual model-theoretic argument, it suffices to show that any instance of (3) over a difference-differential
field $(K,\si,D)$ can be realized in an extension of $(K,\si,D)$.

Let $ (K,\sigma,D)$ be a difference-differential field such that $K \models ACF$.
 Let $U,V,W $ be $K$-varieties satisfying (3). Let $({\mathcal U},D ) $ be a saturated
model of $DCF$ containing $(K,D)$.
Let $(a,b)$ be a $(1,D)$-generic of $W$ over $K$;
 then $a$ is a $(1,D)$-generic of $\pi_1(W)$ over $K$ and $b$ is a $(1,D)$-generic of $\pi_1(W)^{\sigma}$ over $K$.
Hence $tp_{DCF}(b/K)=\sigma(tp_{DCF}(a/K ))  $; thus $\sigma $ extends  to an
automorphism $\sigma'$ of $({\mathcal U},D)  $ such that $\sigma'(a)=b  $.
\smallskip

Now we shall prove that the models of {\it DCFA} are existentially closed.
Let $ (K,\sigma,D)$ be a model of {\it DCFA} contained in a difference-differential field $({\mathcal U},\sigma,D)$.
Since $x \neq 0 \leftrightarrow \exists y \; xy=1    $, it suffices to prove that every finite system of
 $(\sigma,D)$-polynomial equations
over $K$ with a solution in $\mathcal U $ has a solution in $K$.
Let $\varphi(x) $ be such a system
%a positive quantifier-free ${\mathcal L}_{\sigma,D}(K) $-formula, 
and let $a$ be a tuple of $\mathcal U$ satisfying $\varphi $. Since $\sigma$ is an 
automorphism, $\varphi $ is a finite conjunction of equations of the form 
$f(x, \cdots,\sigma^n(x))=0  $, where $f$ is a differential polynomial; 
 such an equation is equivalent, modulo the theory of difference-differential fields,
to a formula of the form:
$$\exists  y_0, \cdots , y_{k-1} f(y_0, \cdots,
 y_{k-1},\sigma(y_{k-1}) )=0 \land \bigwedge_{i=1}^{k-1}(y_i=\sigma(y_{i-1}) \land y_0=x    ).   $$

Thus, if we replace $x$ by $(y_0,\cdots, y_{k-1}  )$ and $ a$ by 
$(a,\cdots,\sigma^{k-1}(a))$, we may suppose that $\varphi $ is a finite 
conjunction of equations of the form $g(x,\sigma(x))=0 $, where  $g(X,Y) $
is a differential polynomial over $K$.

Let $m $ be sufficiently large so that $X$ and $Y$ appear in each $g(X,Y) $ with differential order less than $m$
, and such that, for $M>m $

$$\td((D^{M+1}a,D^{M+1}\sigma(a))/K(a,\sigma(a), \cdots, D^Ma,D^M\sigma(a) )  )=$$
$$\td((D^ma,D^m\sigma(a) )/ K(a,\sigma(a),\cdots , D^{m-1}a,D^{m-1}\sigma(a)     )     )    $$

and 
$$\td(D^{M+1}a/K(a, \cdots, D^{M}a) )=\td(D^{m}a/K(a, \cdots, D^{m-1}a ))    $$

Let $V$ be the locus of $b=(a,Da, \cdots, D^ma)  $ over $K$, 
$U$ the locus of $(b,\sigma(b) )$ over $K$, and let $W \subset \tau_1(V \times V^{\sigma})
  $ be the locus of $(b,Db,\sigma(b),\sigma(Db)) $ over $K$. By construction  and choice of $m$, $b$ is a
 $(1,D)$-generic of $\pi_1(W)$, $\sigma(b) $ is a $(1,D)$-generic of $\pi_2(W)$ and
 $(b, \sigma(b) )$ is a $(1,D)$-generic of $W$.
By axiom (3) there is a tuple $c=(c_0, \cdots,c_m)  $ in $K$ 
such that $(c,Dc,\si(c),\si(Dc)) \in W $. 
Thus $(c_0, \sigma(c_0))$ satisfies all the equations of differential order less than or equal to
$m$ satisfied by $(a,\sigma(a))  $; 
hence $ c_0$ satisfies $\varphi(x)$.\\
$\Box$

\begin{ej}
The following shows why we need the
$(1,D)$-generics in our axioms, generics
are not strong enough to describe differential types.
Consider the set $A$ defined by the equations $\si(x)=Dx$ and $D\si(x)=x^2$.
It is then given by a subvariety $W\subset \tau(\aff^1)\times
\tau(\aff^1)$ given by the equations $x_2=y_1$ and $x_1^2=y_2$. The
variety $W$ projects on each copy of $\tau(\aff^1)$.

Let $a\in A, a\neq 0$. From $\si(a)=Da$ one deduces that
$\si^iD^ja=\si^{i+j}a=D^{i+j}a$ for all $i,j\in\na$. Thus 
$\si^3(a)=(Da)^2=2aDa$, which
implies that   $Da=2a$. Thus there are
differential relations that
cannot be seen from the defining equations. 
\end{ej}

\begin{rem}\label{DCFA410}
If $(K,D,\sigma)$ is a model of {\it DCFA} then $(K,\sigma)$ is a model of {\it ACFA }. 
\end{rem}

{\it Proof} :\\

Take $W=\tau_1(U) $, and apply \ref{DCFA49}.\\
$\Box$\\

For {\it DCFA} we have a version of Hilbert's Nullstellensatz.

\begin{prop}\label{DCFA411}
Let $(K,\sigma,D )$ be a model of {\it DCFA }, and $I$ a perfect $(\sigma,D)$-ideal. Then $I_{\sigma,D}(V(I) )=I  $
\end{prop}

{\it Proof} :\\

Clearly $I \subset I_{\sigma,D}(V(I) )    $.
Let $f \in  K[X]_{\si,D}$, such that $f \notin I$. 
By  \ref{DCFAfinprid}, there is a prime perfect $(\sigma,D)$-ideal $J$ containing $I$ such that 
$f \notin J$. Then $K[{\bar X}]_{\sigma,D}/J $ embeds into a  difference-differential field $L$.
%ring which is a domain, thus its field of fractions embeds in a
% difference-differential field $L$
%which contains $K$.

By \ref{DCFAfingen}, $J$ is generated by a finite tuple of polynomials $P(X)$.
Let $\bar a $ be the image of $\bar X  $ in $L$. Thus we have that $ L \models P({\bar a})=0 $
 and $L \models f({\bar a }) \neq 0 $. Since $(K,\sigma,D)$ is existentially closed   
there is ${\bar b} \in K$ such that
 $P({\bar b})=0 $ and $f({\bar b }) \neq 0 $. But $I \subset J  $,
 thus ${\bar b} \in V(I)  $, which implies $f \not\in I_{\si,D}(V(I))$\\
$\Box$

\begin{defi}\label{DCFA412}
Let $E \subseteq F  $ be two difference-differential fields, let $a \in F  $.
\begin{enumerate}

\item We define $deg_{\sigma, D}(a/E)  $ to be the transcendence degree of $E(a)_{\sigma, D}  $ over $E$ if it is
finite, in this case we say that $a$ is finite-dimensional; otherwise we set $deg_{\si,D}(a/E)=\infty$ and we say that
$a$ has infinite dimension.

\item If $I_{\sigma,D}(a/E)=(0)  $ we say that $a$ is $(\sigma, D) $-transcendental over $E$,
 otherwise we say that it is $(\sigma,D)  $-algebraic over $E$.

\end{enumerate}
\end{defi}

\begin{rem}\label{DCFA413}
If $a$ is $(\sigma,D)$-algebraic over $E$ it is not always true that $deg_{\sigma,D}(a/E)$ is finite.  
\end{rem}

\begin{rem}\label{DCFA414}
There is a natural notion of $(\sigma,D)$-transcendence basis.\

\end{rem}

We mention some consequences of the results of Chapter \ref{chap:prelim}, section \ref{pre:sec3}. 

\begin{fac}\label{DCFA415}
Let $K_1, K_2$ be models of {\it DCFA}, let $E$ an algebraically closed difference-differential 
subfield of $K_1$ and $K_2$. Then $K_1 \equiv _E K_2  $.
\end{fac}

\begin{cor}\label{DCFA416}
Let $E$ be an algebraically closed difference-differential field, then {\it DCFA} $ \cup qfDiag(E)$ is complete, where
$qfDiag(E)$ denotes the set of quantifier-free formulas $\varphi$ with parameters from $E$ which are true in $E$.

\end{cor}

\begin{cor}\label{DCFA417}
Let $(K_1,\sigma_1),(K_2,\sigma_2)$ be two models of {\it DCFA } containing a common 
difference-differential field $(E,\sigma)$.
Then  $   K_1 \equiv _E K_2      $ if and only if $ (E^{alg},\sigma_1 )  \simeq_E  (E^{alg},\sigma_2 ) $.
\end{cor}

\begin{nota}\label{DCFA418}
Let $(K,\sigma,D)$ be a differential-difference field, $A \subset K$. We denote by $cl_{\sigma,D}(A) $ 
the smallest difference-differential field
containing $A$, and by $acl_{\sigma ,D}(A) $, the field-theoretic algebraic closure of $cl_{\sigma,D}(A) $.
\end{nota}

\begin{cor}\label{DCFA419}
Let $E$ be a difference-differential subfield of a model $ K$ of {\it DCFA}. Let  $a,b$ be tuples of $K$.
 Then $tp(a/E) = tp(b/E)  $ if and only if there is an $E$-isomorphism between $acl_{\sigma,D}(E(a))  $ and
 $acl_{\sigma,D}(E(b))  $ which sends $a$ to $b$.
\end{cor}

\begin{cor}\label{DCFA420}
Let $\phi ({\bar x}  )$ be a formula. Then, modulo {\it DCFA}, $\phi ({\bar x}  )  $ 
is equivalent to a disjunction of formulas of the form $\exists {\bar y} \;  \psi ({\bar x}, {\bar y} ) $, 
 where $\psi  $ is quantifier free, and for every tuple $({\bar a}, {\bar b} )   $ in a
 difference-differential field $K$ satisfying $\psi$, $ {\bar b}  \in acl_{\sigma,D} (   {\bar a}  )     $.

\end{cor}

\begin{prop}\label{DCFA421}
Let $(K,\sigma,D) $ be a model of {\it DCFA}. Let $A \subset K$.
Then the \textup{(}model-theoretic\textup{)} algebraic closure  $acl(A) $ of $A$ is $ acl_{\sigma ,D}(A)$.
\end{prop}

As with differential and difference fields, we define independence in difference-differential fields.

\begin{defi}\label{DCFA422}
Let $K $ be a model of {\it DCFA}, let $A,B,C$ be subsets of $K$. 
We say that $A$ is independent from $B$ over $C$, denoted by $A \dnfo_C B$,
 if $acl(A,C)$ is linearly disjoint from $acl(B,C) $ over $acl(C) $.

\end{defi}

%If  $({\mathcal U},\si,D)$ is a saturated model of {\it DCFA}, $M \subset {\mathcal U}$ 
%such that $M $ is a model of {\it DCF} and $\sigma$ is an automorphism of $M$, 
%and  if ${\bar a}, {\bar b}, {\bar c_1}, {\bar c_2} $ tuples in $\mathcal U$ such that:

%\begin{enumerate} 
%\item $tp(  {\bar c_1}/M  )=tp(  {\bar c_2}/M  )    $.

%\item ${\bar a}  \downarrow_M    {\bar c_1}$, ${\bar a}  \downarrow_M    {\bar b}$ and 
% ${\bar b}  \downarrow_M    {\bar c_2}$.
%\end{enumerate}

%Then there is $ {\bar c}  $ realizing $tp(  {\bar c_1}/M \cup   {\bar a}  ) \cup   tp(  {\bar c_2}/M \cup   {\bar b}  )  $
% such that $ {\bar c}  \downarrow_M    ({\bar a} ,  {\bar b})  $.

By  %\ref{PSIM1}
\ref{Paut2}.3 we have:

\begin{teo}\label{DCFA425} \hspace{20cm}
\begin{enumerate} 

\item The independence relation defined above coincides with nonforking.

\item Every completion of {\it DCFA} is supersimple.

\end{enumerate}

\end{teo}

%This characterization of simple theories can be found in  \cite{kimpi}. 

As in {\it ACFA}, we have
a stronger version of the independence theorem.

\begin{teo}\label{DCFA426}

If  $\mathcal U$ is a saturated model of {\it DCFA}, $E$ an algebraically closed subset of $\mathcal U$,
 and  ${\bar a},
 {\bar b}, {\bar c_1}, {\bar c_2} $ tuples in $\mathcal U$ such that:

\begin{enumerate}
\item $tp(  {\bar c_1}/E  )=tp(  {\bar c_2}/E  )    $.

\item ${\bar a}  \dnfo_E    {\bar c_1}$, ${\bar a}  \dnfo_E    {\bar b}$ and
  ${\bar b}  \dnfo_E    {\bar c_2}$.
\end{enumerate}

Then there is $ {\bar c}  $ realizing $tp(  {\bar c_1}/E \cup   {\bar a}  ) \cup   tp(  {\bar c_2}/E \cup   {\bar b}  )  $
 such that $ {\bar c}  \dnfo_E    ({\bar a} ,  {\bar b})  $.
\end{teo}

{\it Proof} :\\

Let $\bar c  $ be a realization of $tp(  {\bar c_1}/E   )  $ such that ${\bar c} \dnfo_E  ({\bar a},{\bar b} )  $.
 Let $A=acl(E {\bar a} ),B=acl(E {\bar b} )  ,C=acl(E {\bar c} )   $.
 Let $\phi_1:acl(E {\bar c_1} )  \to C $ and  $\phi_2:acl(E {\bar c_2} )  \to C $ two
 ${\mathcal L}_{\sigma,D}(E)  $-isomorphisms such that $\phi_i({\bar c}_i)={\bar c}  $.

Let $\sigma_0=\sigma |_{ (AB)^{alg}C} $. Since $A$ is linearly disjoint from $acl(E {\bar c_1})  $
and from $C$ over $E$, we can extend $\phi_1 $ to a ${\mathcal L}_D(A)  $-isomorphism $\psi_1  $
 between $acl(A {\bar c_1})  $ and $(AC)^{alg} (=acl_D(AC ))  $. Let $\sigma_1= \psi_1 \sigma \psi_1^{-1}  $
; $\sigma_1  $ is an automorphism of $(AC)^{alg}  $ and agrees with $\sigma$ on $A$ and $C$.
 By definition of $\sigma_1  $, $\psi_1  $ is a  ${\mathcal L}_{\sigma,D}(A)  $-isomorphism between
 $(acl(A {\bar c_1}),\si) $ and  $((AC)^{alg},\si_1 )$.
In the same way we define $\psi_2:acl(Bc_2)\to (BC)^{alg}  $ and $\sigma_2 \in Aut(BC)^{alg}  $.

Let $L=(AB)^{alg} (AC)^{alg} (BC)^{alg} $ (which is a differential field that extends $A,B,C$).
 Let us suppose that there is an ${\mathcal L}_D  $-automorphism $\tau$ of $L$ which extends 
$\sigma_0,\sigma_1,\sigma_2  $. Let $(M,\tau',D)\models$ {\it DCFA}   contain $(L,\tau,D)  $.
 Since $\tau  $ extends $\sigma_0  $, by \ref{DCFA415}, we have $tp_M(AB/E) =tp_{\mathcal U}(AB/E) $;
since $\tau  $ extend $\sigma_i  $, the $\psi_i  $'s are difference-differential field isomorphisms.
 Applying \ref{DCFA419} we have $tp_M({\bar c}/A )=tp_{\mathcal U}({\bar c_1} /A  )  $ and
$tp_M({\bar c}/B )=tp_{\mathcal U}({\bar c_1} /B  )  $. Also $ {\bar c} \dnfo_E (A,B) $.
Hence to finish the proof, all we have to do is show the existence of such a $\tau$. 
To do this, we will prove that $\sigma_0,\sigma_1  $ have a unique extension $\tau_1  $ 
to $(AB)^{alg} (AC)^{alg} $, and that there is an extension $ \tau_2 $ of $\tau_1,\sigma_2  $
 to $L$ (Note that these automorphisms will commute with $D$).

For the first part it is enough to show that $(AB)^{alg}C  $ is linearly disjoint from  $(AC)^{alg}  $ 
over   $ (AB)^{alg}C    \cap  (AC)^{alg}    $, and that $\sigma_0  $ and $\sigma_1$ agree on 
$(AB)^{alg}C    \cap  (AC)^{alg}   $. Similarly for the second part.

By Remark 2 of 1.9 in \cite{salinas}, we have
$$ (AB)^{alg}C    \cap  (AC)^{alg}  =AC  \;\;(*), \qquad  (AB)^{alg} (AC)^{alg} \cap (BC)^{alg}=BC \;\;(**).$$
Since $(AC)^{alg}$ is Galois over $AC$ it implies that $(AC)^{alg}$ and 
$(AB)^{alg}$ are linearly disjoint over $AC$; as $\si_0$ and $\si_1$ 
both extend $\si$ on $AC$, they are compatible. The same argument applies for the second part. \\
%
%Since  $\sigma $ extends $\sigma_0  $, $\sigma   $ agrees with $\sigma_1 $ on
%$AC $, because $(AC)^{alg} $ is a Galois extension of $AC$; to prove the existence of $\tau_1$
%it suffices to show that $  (AB)^{alg}C    \cap  (AC)^{alg}  =AC  $(*) . Similarly, to prove the existence of
% $\tau_2 $ it suffices to show that $ (AB)^{alg} (AC)^{alg} \cap (BC)^{alg}=BC $(**).
%Proofs for $(*)$ and $(**)$ can be found in \cite{salinas}, Remark 2 of 1.9.\\
$\Box$

\begin{rem}
As {\it DCF} is $\omega$-stable, {\it DCFA} is quantifier-free $\omega$-stable.
\end{rem}

Now we want to prove that {\it DCFA} eliminates imaginaries.
We shall need some properties of the fundamental
order for types in stable theories. 
%The proofs can be found in \cite{poi}.

%\begin{defi}
%Let $T$ be a stable theory in the language $\mathcal L $, $M $ a saturated model of $T$.
%Let $A \subset M$, $p(x) \in S_n(A)$, and let $\varphi(x,y) $ be an ${\mathcal L(A)}$-formula.
%We say that $\varphi(x,y)$ is represented in $p(x)$ if there is a tuple $a$ of $A$ such that $\varphi(x,a) \in p(x)$.
%\end{defi}

%\begin{defi}\label{DCFA427}
%Let $({\mathcal U},D)$ be a saturated differentially closed field, and let $L=acl_D(L),K)acl_D(K) \subset {\mathcal U}$.
%Let $I$ be a prime differential ideal over $L$ and $J$ a a prime differential ideal over $L$.
%We say that $I$ is smaller than $J$ for the fundamental order,
%denoted $I \leq_{fo} J $, if
%for every differential polynomial $P(X,Y)$, theres is a tuple $a$ in $L$ such that $P(X,a) \in I$ then
%there is a tuple $b \in L$ such that $P(X,b) \in J$.
%
%We say that $I$ and $J$ are equivalent in the fundamental order, denoted $I \sim_{fo} J  $,
% if  $I \leq_{fo} J $ and  $J \leq_{fo} J$.
%\end{defi}
%Since in {\it DCF} there is a bijection between types and prime differential ideals, we can extend the define
% the fondamental order on types.
Recall that a type $p(x)$ over some set $A$  {\it represents} the
${\cal L}$-formula $\phi(x,y)$ if there is a tuple $a\in A$ such that
$\phi(x,a)\in p(x)$. We denote by $\beta(p)$ the set of formulas
represented by $p$.

For convenience, we will define the fundamental order on types whose
domain is algebraically closed, so that they are stationary (and
definable by elimination of imaginaries in {\it DCF}).

\begin{defi}
Let $A$ and $B$ be algebraically closed differential subfields of some model $({\cal
U},D)$ of DCF, and let $p(x)$, $q(x)$ be types over $A$ and $B$
respectively. We write $p \leq _{fo}q$ if $\beta(q)\subseteq \beta(p)$,
and $\beta(p)\sim_{fo}$ if $\beta(p)=\beta(q)$. $\leq_{fo}$ is called
the fundamental order.
\end{defi}

\begin{fac}\label{DCFA429}
If $A \subset B$ and $q$ is an extension of $p$, one
has $q\leq_{fo}p$, and $q\sim_{fo}p$ if and only if $q$ is a
non-forking extension of $p$.
\end{fac}

If $p$ and $q$ are types in an infinite number of variables $(x_i)_{i \in I}$ we say that
$p \leq_{fo} q $ if and only if for every finite $J \subset I$, if $p'$ and $q'$ denote
the restrictions of $p$ and $q$ to the variables $(x_i)_{i \in J}$, we have $p' \leq_{fo} q'$.

\begin{rem}\label{DCFA428}
$\sim_{fo} $ is an equivalence relation on the class of types in the variables $(x_i)_{i \in I}$.
\end{rem}

%\begin{fac}\label{DCFA429}
%Let $(\mathcal{U},D)$ be a saturated differentially closed fields, let $K=acl_D(K)\subset L=acl_D(L) \subset {\mathcal U}$.
%Let $a$ be a  tuple of ${\mathcal U}$. Then $tp_(a/L)  \leq_{fo} tp(a/K) $, and $tp(a/L)  \sim_{fo} tp(a/K) $
% if and only if $tp(a/L) $
%does not fork over $K$.
%\end{fac}

\begin{prop}\label{DCFA430}
Every completion of {\it DCFA}  eliminates imaginaries.
\end{prop}

{\it Proof} :\\

Let $(K,\si,D) $ be a saturated model of {\it DCFA}, let $\alpha \in K^{eq}$. Then there 
is a $\emptyset$-definable function $f$ and a tuple $a$ in $K$ such that $f(a)=\alpha$. 

%let $f$ be an $\emptyset $-definable function,
% let $a$ be a tuple of $K$, and let $\alpha $ be the equivalence class of $a$
%for the relation $E(x,y) \longleftrightarrow f(x)=f(y)$.

Let $E=acl^{eq}(\alpha) \cap K$.
If $\alpha$ is definable over $E$, let $b$ be a tuple of $E$ over which $\alpha$
 is definable; then $b \in acl^{eq}(\alpha) $. Since we are working in a field,
 there is a tuple $c$ of $K$ which codes the (finite) set of conjugates of $b$
over $\alpha$. Hence $c$ and $\alpha$ are interdefinable.

Let us suppose that $\alpha$ is not definable over $E$, in particular,
$a$ is not a tuple of $E$. We will show that  there is a realization $b$
 of $tp(a/\alpha) $ such that $b \dnfo_E a$.

We now work in the theory DCF, and replace the tuple $a$ by the
infinite tuple $(\si^i(a))_{i\in {\mathbb Z}}$, which we also denote by $a$.

Since $tp(a/\alpha)$ is non-algebraic, it has a realization $b$ such that
$acl^{eq}(a) \cap acl^{eq}(b)=acl^{eq}(\alpha) $, and thus $$acl(Ea)\cap acl(Eb)=E \eqno{(*)}.$$
Choose such a $b$ such that, if $b'$ satisfies the same properties, then $tp_{DCF}(b'/ acl(Ea) ) \not>_{fo}
 tp_{DCF}(b/acl(Ea))  $.

Let $c$ be a tuple of $K$ such that $tp(c/acl(Ea))=tp(b/acl(Ea))$, and $c \dnfo_{Ea} b $.
Then $f(c)=f(a)$ and $c$ satisfies
$$ acl(Ec) \cap acl(Eab) \subset acl(Ec) \cap acl(Ea)=E, \eqno{(**)}  $$
and there is no $c'$ satisfying $(**)$ such that $tp_{DCF}(c'/acl(Eab)) >_{fo} tp_{DCF}(c/acl(Eab)) $.
%Let $d$ be a tuple realizing $tp(c/Eb) $ such that $d \dnfo_{Eb} a$.
%Then $d$ satisfies $(*)$ and  realists $tp(a/E)$.
% By
%\ref{DCFA429} $tp_{DCF}(d/Eb) \sim_{fo} tp_{DCF}(d/Eab)$, but
% $tp_{DCF}(c/Eab) \leq_{fo} tp_{DCF}(c/Eb) =tp_{DCF}(d/Eb)$; then our hypothesis implies
%$tp_{DCF}(c/Eab)  \sim_{fo} tp_{DCF}(d/Eab)$, and since
%both types extend $tp_{DCF}(c/Eb) $ and $tp_{DCF}(d/Eab) $ is its non-forking extension, we have
% $tp_{DCF}(c/Eab)=tp_{DCF}(d/Eab) $.
%Thus we have  that $c \dnfo_{Ea} b  $ (in {\it DCF}) and   $c \dnfo_{Eb} a  $ (in {\it DCF}),
% but $c \dnfo_{Ea}b $ in {\it DCF}.
%by elimination of imaginaries in {\it DCF}, $c \dnfo_{E} ab $ (in {\it DCF}).
%Then $c \dnfo_{E} ab $ (in {\it DCFA}), hence $a \dnfo_E b $.

Then $tp(c/acl(Eb))\geq_{fo}tp(c/acl(Eab))\sim_{fo}tp(c/acl(Ea))$. Let $\tau$ be an
${\cal L}_D(E)$-automorphism sending $b$ to $a$. Then
$tp(\tau(c)/acl(Ea))\sim_{fo}tp(c/acl(Eb))$, and $\tau(c)$ satisfies $(*)$ (by
$(**)$). Hence, by maximality of $tp(b/acl(Ea))=tp(c/acl(Ea))$, we get that
$tp(c/acl(Eb))\sim_{fo} tp(c/acl(Ea))$, and therefore $c\dnfo_{Eb} a$. By
elimination of imaginaries and $(*)$, this implies that $c\dnfo_E ab$,
and therefore $a\dnfo_Eb$.

We have shown that there is a tuple $b$ realizing $tp(a/\alpha) $ independent from $a$ over $E$.
But $\alpha$ is not $E$-definable, thus there is $a'$ realizing $tp(a/E) $ such that
$f(a) \neq f(a') $, and we may choose it independent from $b$ over $E$. Since $tp(a/E)=tp(a'/E) $,
there is a realization $c'$ of $tp(a'/E)$
such that $f(a')=f(c') $  and $c' \dnfo_E a' $; we may suppose that $c' \dnfo_E b $.
If we apply the independence theorem to $tp(a/Eb) \cup tp(a'/Ec') $ we get a contradiction.\\
$\Box$

\begin{lem}\label{DCFA431}
Let $(K,\sigma,D)$ be a model of {\it DCFA}, let $E=acl(E) \subset K  $, and let $(L, \tau, D) $ be a 
difference-differential field extending $(K, \sigma^n,D)  $, where $n$ is a positive integer. 
Then there is a difference-differential field 
$ ( M,\sigma',D)$ containing $( E, \sigma, D) $ such that $(M,(\si')^n,D) \supset (L,\tau,D)$.
\end{lem}

{\it Proof} :\\

For $i=1, \cdots, n-1  $ let $L_i $ be a difference-differential field realizing $\sigma^i(tp_{DCF}(L/E )  )  $
such that $L_0=L,L_1, \cdots , L_{n-1}  $ are linearly disjoint over $E$. Let $f_0=id_L  $ and for $i=1, \cdots, n-1
$ let $f_i:L \longrightarrow L_i  $ be an ${\mathcal L}_D $-isomorphism extending $\sigma^i  $ on $E$.

For $ i=1, \cdots , n-1 $ let $\sigma_i:L_{i-1} \longrightarrow L_i  $ be defined by $\sigma_i=f_i f_{i-1}^{-1} $, 
and let $\sigma_n: L_{n-1} \longrightarrow L_0  $ be defined by $\sigma_{n}= \tau f_{n-1}^{-1} $.

Let $x \in E  $. If $i=1, \cdots, n-1  $ then 
$\sigma_i(x)=f_i(f_{i-1}^{-1}(x))=\sigma^i(\sigma^{-(i-1)}(x)    )=\sigma(x) $
; and $\sigma_n(x)= \tau(\sigma^{-(n-1)}(x)  ) = \sigma^n(\sigma^{-(n-1)}(x))=\sigma(x) $.
Hence each $\si_i$ extends $\si$ on $E$.

Also, we have $\sigma_n \sigma_{n-1} \cdots \sigma_1=\tau (f_{n-1}^{-1} f_{n-2}) \cdots (f_1 f_0^{-1}) = \tau f_0^{-1}=\tau $.

Let $M$ be the composite of $L_0, \cdots , L_{n-1 } $. Since the  $L_i  $'s are linearly disjoint over $E$,
 $M$ is isomorphic to the quotient field of $L_0 \otimes_E \cdots \otimes_E L_{n-1}$. There is a unique
derivation on $M$ extending the derivations of the $L_i$'s  and there is a unique ${\mathcal L}_D  $-automorphism
$\sigma'  $
 of $M$ which coincides with $\sigma_i $ on $L_{i-1}$ for $i=1, \cdots, n-1$.\\ By the above $(\sigma')^n  $ extends $\tau$ .\\
$\Box$

\begin{cor}\label{DCFA432}
Let $(K,\sigma,D)$ be a model of {\it  DCFA}. Then, for all $n \in {\mathbb N}  $
  $(K,\sigma^n,D)$ is a model of {\it DCFA }.

\end{cor}

{\it Proof} :\\

Let $\Sigma  $ be a finite system of $ (\sigma^n,D)$-equations over $K$, and let $(L,\tau,D) $
 be an extension of $ (K,\sigma^n,D) $ containing a solution of $\Sigma$. By \ref{DCFA431}
 there is an extension $(M,\sigma',D)  $ of $(K,\sigma,D)  $ such that $(M,(\sigma')^n,D ) $ 
is an extension of $(L,\tau,D ) $. Thus, $M$ contains a solution of $\Sigma$, and since $(K,\sigma,D)$ 
is existentially closed, 
$K$ contains a solution of $\Sigma $.\\
$\Box$

\section{The  Field of Constants and the Fixed Field}\label{DCFA5}

In this section we study two special subfields of a model $(K,\si,D)$ of {\it DCFA}: the differential field
$(Fix \si, D)$ and the difference field $({\mathcal C},\si)$ 
where $Fix\si$ is the fixed field of $K$ and ${\mathcal C}$ is the
field of constants of $K$.

Throughout this section $(K,\si,D)$ will denote a model of {\it DCFA}.

%\begin{nota}\label{DCFA51}
%Let $(K,\sigma, D)$ be a model of {\it DCFA}.
% The field of constants of $K$ is ${\mathcal C}= \{x \in K : Dx=0   \}  $; and the fixed field of $K$ is 
%$Fix \sigma= \{x \in K : \sigma (x) =x   \} $.

%\end{nota}

\begin{prop}\label{DCFA52}
$(  {\mathcal C },\sigma)$ is a model of {\it ACFA }.
\end{prop}

{\it Proof} :\\

Since $\sigma$ commutes with $D$,  $({\mathcal C}, \sigma)   $ is a difference field.

Now let $U,V  $ be varieties defined over $\mathcal C$, with $U \subset V \times V^{\sigma}   $ such that
 $U$ projects generically over  $V$ and $V^{\sigma}$.
Let $W=U\times ( {\bar 0}) \subset \tau_1 (V \times V^{\sigma})$.
Then, by \ref{DCFA49} there is $a \in K$ such that $(a,Da,\si(a),\si(Da)) \in W $.
Thus $Da =0  $ and $(a,\si(a)) \in U$.\\
$\Box$
 
\begin{rem}\label{DCFA53}
Clearly the fixed field of $\mathcal C  $ is ${\mathcal C} \cap Fix \sigma$, and is pseudofinite  by \ref{psff}.
Hence ${\mathcal C} \cap Fix \sigma   \prec Fix \si$.
\end{rem}

\begin{rem}\label{DCFA54}
$Fix \sigma  $ is a differential field, however it is not differentially closed since it is not algebraically closed as a
field. 
 Clearly, it is also a difference field, thus $$acl(Fix \sigma) =acl_D(Fix \sigma) =(Fix \sigma) ^{alg} . $$

\end{rem}

%\begin{rem}\label{I54}
%$Fix \sigma  $ is a differential field, however it is easy to see that it is not a model of {\it DCF }. 
% Clearly, it is also a difference field, thus $$acl(Fix \sigma) =acl_D(Fix \sigma) =(Fix \sigma) ^{alg} . $$
%
%\end{rem}

\begin{teo}\label{DCFA55} $((Fix \sigma)^{alg},D   )$ is a model of {\it DCF }.

\end{teo}

{\it Proof} :\\

Let $V,W$ be two irreducible affine varieties defined over $(Fix \si)^{alg}$ such that $W \subset \tau_1(V)$
and $W$ projects dominantly onto $V$.
Let $k \in \na$ be such that both $V$ and $W$ are defined over $Fix \si^k$. 
Let $U=\{(x,x):x \in V\}$.
Then $U \subset V \times V^{\si^k}=V \times V$.

Let $W'= \{(y,y):y \in W\}$. Then $W' \subset \tau_1(U)$. By \ref{DCFA432} $(K,\si^k,D)$ is a model 
of {\it DCFA}; thus, applying \ref{DCFA49} to $V,U$ and $W'$ there is $a \in V(K)$ such that
$(a,\si^k(a)) \in U$ and $(a,Da, \si^k(a),D(\si^k(a))) \in W'$. Thus $a = \si^k(a)$ and $(a,Da) \in W$.
By \ref{P224}, $((Fix \si)^{alg},D)$ is differentially closed.\\
$\Box$

Using (the proof of) \ref{DCFA55}, we can also
 axiomatize the theory of the structures $(F,D)$, where $F$ is the fixed field 
of a model of {\it DCFA}, as follows:
\begin{enumerate}
\item $F$ is a pseudo-finite field.
\item For every  irreducible algebraic variety $V$ defined over $F$, 
if $W$ is an irreducible algebraic subvariety of $\tau_1(V)$ defined over $F$, 
such that the projection of $W$ onto $V$ is dominant, then
there is $a \in V(F)$ such that $(a,Da) \in W$.
\end{enumerate}
For such a structure $(F,D)$ we can describe its completions, the types, the algebraic closure in
 the same way as we did for {\it DCFA}. For instance, if $F_1$ and $F_2$ are two models of
this theory and $E$ is a common substructure, $F_1 \equiv_E F_2$ if and only if
there is an isomorphism $\varphi:E^{alg} \cap F_1 \to E^{alg}\cap F_2$ which fixes $E$. If we add enough constants 
(for  a pseudo finite field $F$ we add a set of constants $A \subset F$ such that $FA^{alg}=F^{alg}$), 
the
generalized independence theorem will hold.
% over structures such that $E=acl(E)$ and $E^{alg}F=F^{alg}$. From this we
%can deduce the elimination of imaginaries.

Pseudo-algebraically closed structures
were studied by E. Hrushovski in a  preprint of 91, to appear in the Ravello Proceedings. In \cite{pildom} Pillay
and Polkowska
generalize Hrushovski's results and treat the differential case described above.\\

We have seen that the field of constants of a model of {\it DCF} as well as the fixed field of
a model of {\it ACFA} are stably embedded (\ref{P217a} and \ref{sistabem}). 
The same  happens in {\it DCFA} for the field $Fix \si$ but not for the field $\mathcal C$.

\begin{prop}
 $({\mathcal C},\sigma)$ is not stably embedded.
\end{prop}

{\it Proof}:\\

Let $a \in Fix \sigma \setminus {\mathcal C}$, then the set 
$\{x \in K: \exists y \ \sigma(y)=y \wedge Dx=0 \wedge y^2=x+a \}$ is contained in ${\mathcal C}$ but it is not
definable with parameters from ${\mathcal C}$.\\
$\Box$

\begin{prop}\label{I57}

Let $A $ be a definable subset of $ (Fix \sigma) ^n $. Then $A$ is definable over $Fix \sigma  $
in the language ${\mathcal L}_D $.

\end{prop}

{\it Proof} :\\

Since {\it DCFA } eliminates imaginaries, there is a canonical parameter $a$ for $A$. 
Since $A$ is fixed by $\sigma $, $a$ is fixed by $\sigma $, thus  $A$ is $(Fix \sigma)$-definable.
It is enough to show that there exist a countable subset $L$ of 
$Fix \sigma$ containing $a$ such that every ${\mathcal L}_D$-automorphism
 of $Fix \sigma$ which fixes $L$ 
extends to an elementary map of some elementary extension of $Fix \sigma$.

Let $L$ be a countable elementary ${\mathcal L}_D$-substructure of $Fix \sigma$ containing $a$. In particular $L$ is a differential field, 
and $acl(L)=L^{alg}   $.

Since  $L \prec _{\mathcal{L}_D} Fix \sigma  $, $L^{alg}$ and $Fix \sigma$ are linearly disjoint over $L$.
If $L_n$ is the unique algebraic extension of $L$ of degree $n$, then $L_n Fix \sigma  $
 is the unique algebraic extension of $Fix \sigma$ 
of degree $n$; this implies that $(Fix \sigma)^{alg}=L^{alg} Fix \sigma   $.

Let $\tau$ be a ${\mathcal L}_D$-automorphism of $Fix \sigma$ over $L$. Then we can extend $\tau$ to a 
${\mathcal L}_D$-automorphism ${\bar \tau}  $ of $L^{alg} Fix \sigma  $
  over $L^{alg}  $. We have that ${\bar \tau  }  $ commutes with $\sigma$. Thus ${\bar \tau   } $
 is a ${\mathcal L}_{\sigma,D } $-automorphism of $acl(Fix \sigma) $. 
Then, by \ref{DCFA419},    ${\bar \tau  }  $ is an elementary map.\\
$\Box$

\begin{rem}
Let $(L,L_A)$ be a pair of fields extending the pair of fields 
$(Fix \si, Fix \si\cap \mathcal{C})$ and which satisfies: $L$ is a regular 
extension of $Fix \si$, $L_A$ is a regular extension of $Fix \si \cap 
\mathcal{C}$, and $Fix \si $ and $L_A$ are linearly disjoint over $Fix \si \cap 
\mathcal{C}$.

Using the linear disjointness of $L_A$ and $Fix \si$ over $Fix \si \cap 
\mathcal{C}$ and \ref{P222.1}, the derivation $D$ of $Fix \si $ extends to a 
derivation $D_1$ on $L$ which is $0$ on $L_A$. Defining $\si$ to be the 
identity on $L$, the difference differential field $(L,D_1,id)$ embeds 
(over $Fix \si $) into an elementary extension of ${\mathcal U}$. The 
following follows easily:
\begin{enumerate}
\item The pair $(Fix \si, Fix \si \cap \mathcal{C})$ is S.P.A.C., that is, if  
$a,b$ are tuples in some extension of $Fix \si$ such that $Fix \si \subset Fix \si(a,b)$ and 
$Fix \si \cap \mathcal{C} \subset (Fix \si \cap \mathcal{C})(a)$ are regular, and
$Fix \si$ is linearly disjoint from $(Fix \si \cap \mathcal{C})(a)$ over $Fix \si \cap \mathcal{C}$; 
then there is a zero $(a',b')$ of $I(a,b/Fix \si)$ such that
$a' \in Fix \si \cap \mathcal{C}$.

This notion was introduced by H. Lejeune, see \cite{lej}.

\item The theory of the structure $Fix \si$ is model complete 
in the following languages:

\begin{enumerate}
\item The language of pairs of fields with enough constants to describe all algebraic extensions
of $Fix \si$, and with $n$-ary relation symbols for all $n$ which interpretation in 
$(Fix \si, Fix \si \cap \mathcal{C})$ is that
the elements $x_1, \cdots, x_n$ are $(Fix \si \cap \mathcal{C} )$-linearly independent. 
\item The language of differential fields with enough constants to describe all algebraic extensions
of $Fix\si$ (as in this language extensions are field extensions with an extension
of the derivation this will automatically imply linear disjointness). 
\end{enumerate}

\end{enumerate}
%If instead of the predicate symbols that we add to have linearly disjoint extension, 
%we add the derivation $D$  we have the same result.

%\item A pair of fields $(F,F_A)$, where $F_A \subset F$ is regular, is said to be S.P.A.C (for simultaneously P.A.C.), if it
%satisfies the following: Let  
%$a,b$ be tuples in some extension of $F$ such that $F \subset F(a,b)$ and $F_A \subset F_A(a)$ are regular, and
%$F$ is linearly disjoint from $F_A(a)$ over $F_A$. Let $I=I(a,b/F)$. Then there is a zero $(a',b')$ of $I$ such that
%$a' \in F_A$.
%By the discussion in 2, 
%the theory of the pair $(Fix \si, Fix\si\cap {\mathcal C})$ is S.P.A.C. (see 2.5 of \cite{lej}).

% also it is
%a beutiful  pair which, by Corollary 5.7 of \cite{pildom}, is a first-order property.
\end{rem}

\section{Forking and the SU-Rank}\label{DCFA6}

Since every completion of {\it DCFA} is supersimple, types are ranked by the $\s$-rank (\ref{pre:sec1}).
This section is devoted to the study of the $\s$-rank in {\it DCFA}. Given an element of a model of 
{\it DCFA} we will construct a sequence and we will define a rank for this sequence and we will show that this
rank bounds the SU-rank of the element. With this we prove that the SU-rank of a model of
{\it DCFA} is $\omega^2$.

\begin{rem}\label{DCFA66}
Let $E=acl(E) $, and let us suppose that  $deg_{\sigma,D}(a/E) < \infty $. Let $F=acl(F)\supset E $. 
Then $a \dnfo_E F$ if and only if $deg_{\sigma,D}(a/F)<deg_{\sigma,D}(a/E) $.
 Thus, by induction on $deg_{\sigma,D}(a/E) $ we can prove that   $\s(a/E) \leq   deg_{\sigma,D}(a/E). $
\end{rem}

In order to compute the $\s$-rank of a completion of {\it DCFA}, to any type we will associate a sequence in
${\mathbb N} \cup \{ \infty  \} $ and we will define a rank for such a sequence; in some cases this rank
will bound the $\s$-rank of the type.
Let $(I,\leq)$ be the class of decreasing sequences of ${\mathbb N} \cup \{ \infty  \} $
 indexed by ${\mathbb N } $, partially ordered as follows:
 If $(m_n), (m'_n) \in I  $, then $ (m_n) \leq  (m'_n)   $ if and only if for every $n \in {\mathbb N}  $,
 $m_n \leq m'_n $. We write $(m_n)<(m'_n)$ if $(m_n)\leq(m'_n)$ and $(m_n)\neq (m'_n)$

\begin{rem}\label{DCFA68}
If $(m_n) \in I  $, then there exist $A \in  {\mathbb N} \cup \{\infty \} $ and 
 $B,C \in {\mathbb N}  $ such that $m_n = \infty $ if and only if 
$n < A$, and  $m_n =C  $ if and only if  $n \geq A+B $.

\end{rem}

\begin{defi}\label{DCFA69}
Let $(m_n) \in I  $. We define the Foundation Rank of $(m_n)$, denoted $\f(m_n)  $ as follows.
Let $\alpha$ be an ordinal:
\begin{enumerate}

\item $\f(m_n) \geq 0 $.
\item $ \f(m_n) \geq \alpha +1  $ if there is $(m'_n) \in I  $ such that $(m_n) > (m_n') $ and $ FR(m'_n) \geq \alpha  $.
\item If $\alpha  $ is a limit ordinal, then $ \f(m_n) \geq \alpha$    if  $ \f(m_n) \geq \beta  $ for every $\beta < \alpha  $.
\item $\f(m_n) $ is the smallest ordinal $\alpha  $ such that  $\f(m_n)  \geq \alpha $ but
  $\f(m_n)    \geq \! \! \! \!  \!  /   \, \,   \alpha +1 $.

\end{enumerate}
\end{defi}

\begin{defi}\label{DCFA612}
Let $(K,\si,D)$ be a model of {\it DCFA}, $E=acl(E) \subset K  $, and $ a \in K $. 
To $a$ and $E$ we associate the sequence $(a^E_n) $ defined by:
 $$a^E_n= \td(E(a,Da, \cdots, D^{n}a)_{\sigma}/E(a,Da, \cdots , D^{n-1}a    )_{\sigma}).$$ 
\end{defi}

\begin{rem}\label{DCFA613} \hspace{20cm}    
\begin{enumerate}
\item By \ref{P226}, $(a^E_n) \in I  $.
\item Assume that either $a$ is a single element, or that $a$ is $\si$-algebraic over $E$. If $E \subset F=acl(F)  $, 
then $tp(a/E)  $ does not fork over $F$ if and only if $a \dnfo_E F$, if and only if
for all $n \in \na$, $\td(E(a,Da, \cdots, D^{n}a)_{\sigma}/E(a,Da, \cdots , D^{n-1}a)_{\sigma})=
\td(F(a,Da, \cdots, D^{n}a)_{\sigma}/F(a,Da, \cdots , D^{n-1}a    )_{\sigma}) $, 
if and only if $(a^E_n) =  (a^F_n) $. Hence $\s(a/E) \leq \f(a^E_n)  $.

\end{enumerate}
\end{rem}

\begin{prop}\label{DCFA610}
Let $(m_n) \in I $, let $A,B,C $ as in \ref{DCFA68}. If $A \neq \infty$ then
$\f(m_n)= \omega \cdot (A+C)  +\sum_{j=A}^{A+B-1}(m_j-C)$; if $A=\infty$ then $\f(m_n)=\omega^2$.
\end{prop}

{\it Proof} :\\

First we observe that if $B' > B$, then $\sum_{j=A}^{A+B'-1}(m_j-C)=\sum_{j=A}^{A+B-1}(m_j-C)$.
We proceed by induction on the ordinal $\alpha= \omega \cdot  (A+C) +\sum_{j=A}^{A+B-1}(m_j-C) $.
For $\alpha=0$ it is clear.\\
Suppose that the theorem holds for $\alpha$.

Let $(m_n)\in I $, and $A,B,C $ as in \ref{DCFA68}, such that
 $\alpha+1=   \omega \cdot (A+C) +\sum_{j=A}^{A+B-1}(m_j-C)  $; 
this implies in particular
that $B \neq 0 $ and $m_{A+B-1}>C$.\\\\
$\f(m_n)>\alpha $:

Let $(m'_n)\in I $ such that $m_n'=m_n $ for $n \neq A+B-1 $ and $m_{A+B-1}'=m_{A+B-1}-1 $, 
so that $(m_n') \in I$ and $(m_n')<(m_n) $. Let $A',B',C' $ 
be the numbers associated to $(m_n')$ by \ref{DCFA68}. 
Then $A'=A$, $C'=C$, $B' \leq B$ and $\omega \cdot  (A'+C') +\sum_{j=A'}^{A'+B'-1}(m'_j-C')=\alpha$. By induction 
hypothesis $\f(m'_n)=\alpha <\f(m_n)$.\\\\
$\f(m_n)=\alpha+1$:

Let $(m_n') \in I $ such that $(m_n')<(m_n) $. Let  $A',B',C' $ be the numbers associated to $(m_n')$ by \ref{DCFA68}. 
Then $A' \leq A $ and $C' \leq C $. We want to show that $\f(m'_n)\leq\alpha$.

If $A'<A $ or $C'<C $ we have $A'+C'<A+C $, and thus $ \omega  \cdot \ (A'+C')<\omega  \cdot(A+C)  $. 
Since $\sum_{j=A'}^{A'+B'-1}(m'_j-C') \in \mathbb{N} $,
$\alpha+1=\omega  \cdot (A+C) +\sum_{j=A}^{A+B-1}(m_j-C)>\omega \cdot (A'+C') +\sum_{j=A'}^{A'+B'-1}(m'_j-C') $ and 
by induction hypothesis the latter equals $FR(m_n') $.

If $A'=A $ and $C'=C $, then there is $k \in \{ A, \cdots,A+B-1\} $
 such that $m_k'<m_k $. In this case we have 
$\sum_{j=A}^{A+B-1}(m'_j-C)<\sum_{j=A}^{A+B-1}(m_j-C) $,
 hence $\alpha+1=\omega \cdot  (A+C) +\sum_{j=A}^{A+B-1}(m_j-C)>\omega \cdot (A+C)+\sum_{j=A}^{A+B-1}(m'_j-C) $.
This shows the result in the case $\alpha +1$
 
Assume now that $\alpha$ is a limit ordinal $< \omega^2$, and let $(m_n) \in I$ (with the associated numbers $A,B=0,C$) such that
 $\alpha=\omega \cdot (A+C)$ with $A+C \neq 0$ ($B=0$).

We shall prove that for every $k \in {\mathbb N}$ there is $(m'_n) \in I$ such that $(m'_n)< (m_n)$ and
$\f(m'_n)=\omega \cdot (A+C-1)+k$.
 
If $A \neq 0$,  let $(m_n') \in I $ be such that $m_{A-1}'=C+k $, $m'_n= \infty$ for $n<A-1$  and
$m_n'=C $ for $n>A-1$. We have $(m_n')<(m_n) $ and by induction hypothesis $\f(m_n')=\omega \cdot (A+C-1)+k  $.

If $A=0 $, then $C \neq 0$. Let $(m_n')\in I $ such that $m_n'=C-1 $ for $n \geq k$ and $m_n'=C$ if $n < k $.
 Then 
$(m_n')<(m_n) $ and by induction hypothesis $\f(m_n')=\omega \cdot (C-1)+k$ .
Thus $\f(m_n) \geq \alpha$.\\

$\f(m_n)=\alpha$:

Let $(m_n') \in I$ such that $(m_n')<(m_n)$, let $A',B',C'$ be the numbers associated to $(m_n')$ by \ref{DCFA68}.
 Then $A'<A $ or
$C'<C $, hence $A'+C'<A+C $,
 and $\omega  \cdot  (A'+C') +\sum_{j=A'}^{A'+B'-1}(m'_j-C')<\omega \cdot  (A+C)=\alpha $. By induction
hypothesis $\f(m_n') < \alpha$. This shows that $\f(m_n)\not\geq \alpha+1$, i.e. $\f(m_n)=\alpha$.\\
$\alpha=\omega^2$:
Let $(m_n)$ be the sequence defined by $m_n=\infty$ for all $n\in \na$. 
By induction hypothesis we know that if $(m'_n)<(m_n)$ is in $I$, then $\f(m'_n)< \omega^2$.
Hence $\f(m_n) \not \geq \omega^2+1$. On the other hand, for every
$k \in \omega$, let $(m_n^k)$ be the sequence with associated numbers $A=k,B=C=0$. 
Then $\f(m_n)>\f(m^k_n)=\omega k$.\\ 
$\Box$

%\begin{cor}\label{DCFA611}
%The maximal element of $(I, \leq ) $ has foundation rank $\omega^2 $.
%
%\end{cor}

\begin{prop}\label{DCFA616}
Let $(K,\si,D)$ be a model of {\it DCFA} and let $a \in K  $ be $(\si,D)$-transcendental over $F=acl(F) \subset K $.
 Let $(m_n) \in I  $. Then there is a difference-differential field $E \subset K$ such that $(a_n^E) =(m_n)$.
\end{prop}

{\it Proof} :\\

Define $b_0=a, b_1=\si(a)-a, \cdots, b_{k+1}=\si(b_k)-b_k, \cdots$ and $b_{\infty}=1.$. 
%Let $L^0_k= F(b_k^0)_{\sigma} $;
% then \[ F(a)_{\sigma}=L_0^0 \supset L_k^0 \supset L^0_{k+1} \supset  F   \] and $\td(L_k^0/L^0_{k+1})=1  $.
% Then $\td(F(a)_{\sigma} / L^0_{m_0}   )=m_0  $.
%Let $E_0=L^0_{m_0}  $, then $E_0(a)_{\sigma}=  F(a)_{\sigma}  $. Hence $\td(E_0(a)_{\sigma}/E_0  )=m_0  $.
%
%Now let us suppose that we have constructed  a difference-differential field $E_n$ such that for all $i \leq n $
%we have 
%$F(a)_{\sigma}  \supset E_n  \supset E_i     $ and also
%$\td(E_n(D^i(a))_{\sigma} / E_n(a, \cdots, D^{i-1}a)_{\sigma}   )=m_i  $.
%Let $$ $b^{n+1}_0=D^{n+1}a,  b^{n+1}_{k+1}= \sigma(b^{n+1}_k)- b^{n+1}_{k}   $.
Let $E=F(D^ib_{m_i}:i \in \na)_{\si}$. For all $i$, $D^{i+1}b_{m_i} \in F(D^{i+1}b_{m_{i+1}})_{\si})$, hence
$E$ is a difference-differential field and by construction and because $a$ si $(\si,D)$-transcendental over $F$, 
$$\td(E(a, \cdots,D^{n+1}a)_{\si}/E(a,\cdots, D^na)_{\si})=$$

$$=\td(F(a,\cdots,D^{n+1}a)_{\si}/F(a,\cdots,D^na,D^{n+1}b_{m_{n+1}})_{\si})$$
and as $a$ is $(\si,D)$-transcendental over $F$, the latter equals $m_{n+1}$.

%\begin{array}
%\td(E(a, \cdots,D^{n+1}a)_{\si}/E(a,\cdots, D^na)_{\si})$&$=$\\
%$=$&\td(F(a,\cdots,D^{n+1}a,D^nb_{m_n})_{\si}/F(a,\cdots,D^na,D^nb_{m_n})_{\si})$\\
%\end{array}

% let  $L^{n+1}_k= E_n(b_k^{n+1})_{\sigma} $. Then for all $k$,
% \[ E_n(D^{n+1}a)_{\sigma}=L_0^{n+1} \supset L_k^{n+1} \supset L^{n+1}_{k+1} \supset  E_n   \]
%and $trdg(L_k^{n+1}/L^{n+1}_{k+1}     )=1  $. Thus $trdg( L^{n+1}_0 / L^{n+1}_{m_{n+1}}   )=m_{n+1}  $.
%Let $E_{n+1}=  L^{n+1}_{m_{n+1}}    $.
%We have  $E_{n+1}(D^{n+1}a)_{\sigma}=E_{n}(D^{n+1}a)_{\sigma}    $, 
%this implies   $$\td(E_{n+1}(D^{n+1}a)_{\sigma} / E_{n+1}   )=m_{n+1}, $$ 
%and  since $a$ is $(\sigma,D ) $-transcendental, 
%$$\td(E_{n+1}(D^{n+1}a)_{\sigma} / E_{n+1}(a, \cdots ,D^na  )_{\sigma}   )=m_{n+1}. $$
%Let $0<i<n+1  $, since $E_{n+1}=E_n(b^{n+1}_{m_{n+1}})_{\sigma}   $, 
%$E_{n+1}(D^{i}a)_{\sigma}=  E_n(b^{n+1}_{m_{n+1}},D^ia )_{\sigma}  $, and
% $ E_{n+1}(a, \cdots ,D^{i-1}a  )_{\sigma}   = E_n(b^{n+1}_{m_{n+1}},a, \cdots ,D^{i-1}a    )_{\sigma}     $; 
% since $ a $ is   $(\sigma,D ) $-transcendental, 
%$\td(E_{n+1}(D^{i}a)_{\sigma} / E_{n+1}(a, \cdots ,D^{i-1}a  )_{\sigma})$ $=  \td(E_{n}(D^{i}a)_{\sigma} /
% E_{n}(a, \cdots ,D^{i-1}a  )_{\sigma}   ) =m_i.$
%Let $E=\cup _{n \in {\mathbf N}} E_n  $, then $E $ is a difference-differential field and $(a^E_n)=(m_n )  $.\\

$\Box$

%\begin{prop}\label{DCFA614}
%Let $(K,\si,D)$ be a model of {\it DCFA}  , $E=acl(E) \subset K  $, and $a \in K$ $(\sigma,D ) $-transcendental over $E$.
%Then $\s(a/E)=\omega^2.  $
%\end{prop}
%
%{\it Proof} :\\
%
%By \ref{DCFA611} and \ref{DCFA613}, $SU(a/E) \leq \omega^2  $.
%\\
%Let $b_0=a $ and $b_{n+1}= \sigma(b_n)-b_n $. Let $L_n=E(b_n)_{\sigma,D}  $.
% Thus $L_n \supseteq L_{n+1}  $ and  $L_n=L_{n+1}(b_n)  $.\\
%Let $c_{n,i}=D^ib_n  $. Let $M_{n,i}=L_{n+1}(c_{n,i})_{\sigma,D} $  .  
 %We have that $L_{n} \supset M_{n,i} \supset M_{n,i+1}  \supset L_{n+1}  $ are
 %difference-differential fields, and $ M_{n,i}=  M_{n,i+1}(c_{i,n})  $,
 %then  $\td(M_{n,i}/M_{n,i+1})=1 $;
 %hence $\s(c_{n,i}/L_{n+1} c_{n,i+1})=1  $. Since $E(b_n)_{\sigma,D} =L_n  $, by \ref{PSU3},
%$\s(b_n/M_{n,i} )=i    $ , thus $\s(b_n/E b_{n+1})= \omega  $.
% This implies that $\s(a/E b_n)= \omega \cdot n  $; hence $\s(a/E) \geq \omega^2 $.\\
%$\Box$

\begin{cor}\label{DCFA615}
Let $a$ be a tuple of $K$ such that, the elements of $\{\sigma^i(D^ja):i,j \in \na \}$
 are algebraically independent over $E$.
 Let $n$ be the transcendence degree of $a$ over $E$. 
Then $\s(a/E)= \omega^2 \cdot n . $  

\end{cor}

\section{Remarks on Stability, Stable Embeddability and 1-basedness}

We know that no completion of {\it DCFA} is stable. As in the case of 
completions of {\it ACFA}, it turns out that certain definable sets, 
endowed  with the structure induced by the ambient model, are stable 
 stably embedded. In this section we discuss how to apply results 
from \cite{salinas} to obtain similar results in models of {\it DCFA}.
We also give a criterion for 1-basedness in {\it DCFA}.
%Def of stable, stably embedded (cf appendix 1 of [CH])
\begin{defi} 
A (partial) type $p$ over a set $A$ is stable stably embedded
 if whenever $a$ realises $p$ and $B\supset A$, then $tp(a/B)$ is definable.
Equivalently, let $P$ denote the set of realizations of $p$. Then $p$ is 
stable stably embedded if and only if for all set $S\cap P^n$ where $S$ 
is  definable, there is a set $S'$ definable with parameters from $P$ 
and such that $S'\cap P^n=S\cap P^n$.
\end{defi}

[Note: if $p$ is complete, this is what Shelah calls a stable type]. \\

The following result is proved in the Appendix of \cite{salinas}: 

\begin{lem}\label{st0} 
If $tp(b/A)$ and $tp(a/Ab)$ are stable stably embedded, so is $tp(a,b/A)$. 
\end{lem}

In \cite{salinas}, a certain property (called superficial stability) is 
isolated, and guarantees that certain types over algebraically closed 
sets are stationary, and therefore definable. It follows from model 
theoretic considerations that if for any algebraically closed set $B$ 
containing $A$, $tp(a/B)$ is stationary, then  $tp(a/A)$ will be 
stable and stably embedded.

\begin{lem}\label{lemsse}
Let $(K,\si)$ be a model of {\it ACFA}, $A=acl_\si(A)\subset K$ and 
$a\in K$. Then $tp(a/A)$ is  stationary if and only if $tp(a/A)\perp 
(\si(x)=x)$.

\end{lem}

{\it Proof}:\\

Indeed, write $SU(a/A)=\omega k+n$, and let $b\in
acl_\sigma(Aa)$ be such that $SU(b/A)=n$. Then $tp(b/A)\perp
(\sigma(x)=x)$, and by Theorem 4.11 of [2], $tp(acl_\sigma(Ab)/A)$ is
stationary. If $c\in acl_\sigma(Aa)$ satisfies some 
non-trivial difference equation over $acl_\sigma(Ab)$ then 
$SU(c/Ab)<\omega$ 
and therefore $c\in acl_\sigma(Ab)$. Hence, by Theorem 5.3 of [3], 
$tp(a/acl_\sigma(Ab))$ is stationary, and therefore so is $tp(a/A)$.

For the converse, there are independent 
realizations $a_1,\cdots,a_n$ of $tp(a/A)$, and elements 
$b_1,\cdots,b_m\in Fix\si$ such that $(a_1,\cdots,a_n)$ and 
$(b_1,\cdots,b_m)$ are not independent over $A$. Looking at the field of 
definition of the algebraic 
locus of $(b_1,\cdots,b_m)$ over $acl_\si(A,a_1,\cdots,a_n)$, there is 
some $b\in Fix 
\si\cap acl_\si(A,a_1,\cdots,a_n)$, $b \not\in A$. Then $tp(b/A)$ is not 
stationary: if $c \in Fix\si$ is independent from $b$ over $A$, then 
$tp(b/A)$ has two distinct non-forking extensions to $Ac$, one in which 
$\sqrt{b+c}\in Fix\si$, the other in which $\sqrt{b+c}\not\in Fix\si$.
Hence $tp(a_1,\cdots,a_n/A)$ is not stationary, and neither is 
$tp(a/A)$.\\
$\Box$

%Let $(K,\si)$ be a model  of {\it ACFA} of characteristic $0$, $a$ a tuple in 
%$K$, and $A=acl(A)\subset K$. If $tp(a/A)$ is  orthogonal to the 
%formula $\si(x)=x$, then $tp(a/A)$ is stationary (see \cite{salinas}, and 
%\cite{salinas2}). It follows that if $tp(a/A)$ is hereditarily orthogonal to 
%$\si(x)=x$, then $tp(a/A)$ is stable stably embedded. In that 
%case, it will also be $1$-based.

It is important to note that stationarity alone does not imply 
stability: if $a$ is transformally transcendental over 
$A=acl_\si(A)$, then $tp_{ACFA}(a/A)$ is stationary, but it is not stable.

These results can be used to give sufficient conditions on types in 
{\it DCFA} to be stationary, and stable stably embedded.

\begin{prop}\label{st1}
Let $(K,\si,D)$ be a model of DCFA, let $A=acl(A)\subset K$, and 
$a$ a tuple in $K$.
\begin{enumerate}
\item Assume that $tp_{ACFA}(a,Da,D^2a,\cdots/A)\perp \si(x)=x$. Then 
$tp(a/A)$ is stationary.

\item Assume that for every $n$, $tp_{ACFA}(D^na/Aa\cdots D^{n-1}a)$ is 
hereditarily orthogonal to $(\si(x)=x)$. Then $tp(a/A)$ is stable 
stably embedded. It is also 1-based.
\item If $tp(a/A)$ is not hereditarily orthogonal to $(\si(x)=x)$, then 
$tp(a/A)$ is not stable stably embedded. 
\end{enumerate}
\end{prop}

{\it Proof}:\\

1. As $tp_{ACFA}(a,Da,D^2a,\cdots/A)\perp \si(x)=x$,  \ref{lemsse}  implies that $tp_{ACFA}(a,Da,D^2a,\cdots/A)$ 
is stationary. Let $b,c$ be two realizations of non-forking extensions of $tp(a/A)$ to a set $B=acl(B) \supset A$.
As $tp_{ACFA}(a,Da,D^2a,\cdots/A)$ is stationary we have that
 $tp_{ACFA}(b,Db,D^2b,\cdots/B) =tp_{ACFA}(c,Dc,D^2c,\cdots/B)$. If $\varphi(x)$ is an
$\mathcal{L}_{\si,D}(B)$-formula satisfied by $b$, then there is a $\mathcal{L}_{\si}(B)$-formula
$\psi(x_0,\cdots,x_k)$ such that $\phi(b)=\psi(b,Db,\cdots,D^kb)$; so we have 
$\psi(b,Db,\cdots,D^kb)\in tp_{ACFA}(b,Db,D^2b,\cdots/B) =tp_{ACFA}(c,Dc,D^2c,\cdots/B)$. This implies that
$tp(b/B)=tp(c/B)$, and thus $tp(a/A)$ is stationary.\\
2. By \ref{lemsse} for all $n \in \na$ and for all $B \supset A$, $tp_{ACFA}(D^na/Ba\cdots D^{n-1}a)$ is stationary. Thus
for all $n$,  $tp_{ACFA}(D^na/Aa\cdots D^{n-1}a)$ is stable stably embedded and 1-based.
By \ref{st0}  stable stable embeddability is preserved by extensions, 
hence $tp_{ACFA}(a,Da,\cdots/A)$ is stable stably embedded, and this implies that all extensions
to algebraically closed sets are stationary. As above, we deduce that all extensions of
$tp(a/A)$ to algebraically closed sets are stationary, hence $tp(a/A)$ is stable stably embedded. 
By \ref{wgr} we have also that $tp_{ACFA}(a,Da,\cdots/A)$ est 1-based. By the definition of independence 
in difference-differential fields and the fact that $acl(A,a)=acl_{DCF}(A,a,Da,\cdots)$ $tp(a/A)$ is 1-based:
Let $A\subset B=acl(B) \subset C=acl(C)$ , and let   $b$ be tuple
 of realisations of $tp(a/A)$.
% Let us suppose that $acl(Bb) \cap  C=B$.
 By hypothesis $tp_{ACFA}(a,Da, \cdots/A)$ is
 $1$-based,  therefore  $(b,Db, \cdots)$ is independent from $C$ over $B$ in {\it ACFA}.
% As $A=acl(A)$, we have that for all $n \in \na$, $tp_{DCF}(\si^n(b)/A)$ is
% $1$-based. By \ref{wgr} $tp_{DCF}(b,\si(b),\ldots, \si^n(b)/B)$ is also
% $1$-based.
 Hence, $(b, Db,\ldots, D^nb)$ is {\it ACFA}-independant from $C$ over $B$,
 for every $n \in  \na$.
 Then for every finite subset $S$ of 
$acl(Bb)$, $B(S)$ is linearly disjoint from $C$ over $B$ (that is 
because 
every such $S$ is such that $B(S)$ contained in $acl_{\si}(B,b,Db, 
\cdots, 
D^nb)$ for some $n$). Thus by definition of linear disjointness 
$acl(Bb)$ is linearly disjoint from $C$ over $B$. So
$b$ is {\it DCFA}-independent from $C$ over $B$.

3. If $tp(a/K)$ is not hereditarily orthogonal to $\si(x)=x$ then there is $B=acl(B) \supset A$
such that $tp(a/B) \not \perp \si(x)=x$. 
Then there are independent 
realizations $a_1,\cdots,a_n$ of $tp(a/B)$, and elements 
$b_1,\cdots,b_m\in Fix\si$ such that $(a_1,\cdots,a_n)$ and 
$(b_1,\cdots,b_m)$ are not independent over $B$.

If we look at the field of 
definition of the algebraic 
locus of $(b_1,\cdots,b_m)$ over $acl(A,a_1,\cdots,a_n)$, we can find  $b\in Fix 
\si\cap acl(A,a_1,\cdots,a_n)$, $b \not\in A$.  
Then $tp(b/A)$ is not stationary: 
Let $c \in Fix\si$ be independent from $b$ over $A$, then 
$tp(b/A)$ has two distinct non-forking extensions to $Ac$, one in which 
$\sqrt{b+c}\in Fix\si$, the other in which $\sqrt{b+c}\not\in Fix\si$.
Hence $tp(a_1,\cdots,a_n/A)$ is not stationary, and neither is 
$tp(a/A)$.

% If $tp(a/K)$ is not hereditarily orthogonal to $\si(x)=x$ then there is $B=acl(B) \supset A$
%such that $tp(a/B) \not \perp \si(x)=x$. Then there is $n \in \na$ such that 
%$tp_{ACFA}(a,Da,\cdots, D^na/B) \not\perp \si(x)=x$. By \ref{lemsse}, $tp_{ACFA}(a,Da, \cdots, D^na/A)$ is not
%stationary, then $tp(a/B)$ is not stationary, hence it is not stable stably embedded.\\
%We claim that $tp_{ACFA}(a,Da,D^a,\ldots /A) \not\perp (\si(x)=x)$. 
%(1) If $tp(a/A)$ is not stationary there is some $m$ such that $tp_{ACFA}(a,Da,\ldots,D^ma/A)$ is not stationary. 
%Let $p_1,\ldots,p_n$ be (${\mathcal L}_\si$-)types of $\s$-rank $1$ or $\omega$ (computed in {\it ACFA}) such that 
%$tp_{ACFA}(a,Da,\cdots,D^ma/A)$ is domination equivalent to $p_1\times \cdots \times p_n$.
% Then one of the $p_i$'s is not stationary, and by 5.3 of \cite{zh1}, it must have $\s$-rank $1$.
% Theorem 4.11 of \cite{salinas} finishes the proof.\\
%(2) is implied by (1) and \ref{wgr}.\\
$\Box$

\begin{rem}\label{st2}
Let $A,K$ and $a$ be as above.
\begin{enumerate}
\item If $\s(a/A)=1$, then the stationarity of $tp(a/A)$ implies its 
stability and stable embeddability.

\item There are examples of types of $\s$-rank $1$ which satisfy (1) 
above but do not satisfy (2). Thus condition (2) is not 
implied by stationarity.

%\item $tp(ab/A)$ is stable stably embedded (resp. 1-based) if and 
%only if $tp(a/A)$ and $tp(b/Aa)$ are stably embedded (resp. 
%1-based). (See the appendix of \cite{salinas} and \cite{wag}).
\end{enumerate}
\end{rem}

\begin{cor}\label{st3}
Let $A=acl(A)$, and $a$ a tuple in ${\cal C}$.
 Then $tp(a/A)$ is stable stably embedded if and only if 
$tp_{ACFA}(a/A)$ is stable stably embedded. In this case, it will also be $1$-based. 
\end{cor}

\begin{prop}\label{st4}
Let $A=acl(A)\subset K$, and $a$ a tuple in $K$, 
with $\s(a/A)=1$. If $tp_{ACFA}(a/A)\perp (\si(x)=x)$ then $tp(a/A)$ is stable stably
embedded. In particular, if  $tp_{ACFA}(a/A)$ is stable stably embedded, then so is $tp(a/A)$. 
\end{prop}

{\it Proof}:\\

 %By \ref{st0}, it suffices to show it when $\s(a/A)=1$.
%Assume that $tp_{ACFA}(a/A)$ is 
%stable stably embedded. 
Suppose that $tp(a/A)$ is not stable stably embedded; then there is 
$B=acl(B)\supset A$ such that $tp(a/B)$ is not 
stationary, and therefore $tp_{ACFA}(a,Da,D^2a,\ldots /B)$ is not stationary. 

By \ref{st1} 
$tp_{ACFA}(a,Da,D^2a,\ldots /A) \not\perp (\si(x)=x)$. 
Hence, there is some algebraically closed difference field $L$ containing $A$, 
which is linearly disjoint from $acl(Aa)$ over $A$, and an element $b\in Fix\si\cap 
(Lacl(Aa))^{alg}, b \not\in L$. Looking at the coefficients of the minimal polynomial of $b$ over $Lacl(Aa)$, 
we may assume that $b\in Lacl(Aa)$.
Let $M=acl(L)$, and chose $(M',L')$ realising $tp(M,L/A)$ and independent from 
$a$ over $A$. Then $qftp_{ACFA}(L'/Aa)=qftp_{ACFA}(L/Aa)$ and there is $b'\in L'acl(Aa)$ such that $\si(b')=b'$. 
Since $\s(a/L')=1$, we get $a\in acl(L'b')=L(b')_D^{alg}$. This implies that $tp_{ACFA}(a/L')\not\perp (\si(x)=x)$,
 and gives us a contradiction.\\
$\Box$

\begin{rem}\label{st5}
As stated, the result of \ref{st4} is false if one only 
assumes $\s(a/A)<\omega$. The correct formulation in that case is as 
follows:

Assume $\s(a/A)<\omega$ and that $acl_{\sigma}(Aa)$ contains a sequence 
$a_1,\cdots,a_n$ of tuples such that, for all $i\leq n$, working in 
{\it DCFA},  $\s(a_i/Aa_1,\cdots,a_{i-1})=1$.
% or 
%$tp(a_i/Aa_1,\cdots,a_{i-1})$ is almost internal to the set of 
%conjugates of some type of $\s$-rank $1$.
Under these hypotheses, if $tp_{ACFA}(a/A)$ is stable stably embedded 
then so is $tp(a/A)$.

%Then $tp(a/A)$ is stable 
%stably embedded if and only if $tp_{ACFA}(a/A)$ is stable stably 
%embedded, if and only if for each $i$, 
%$tp_{ACFA}(a_i/Aa_1,\cdots,a_{i-1})\perp (\sigma(x)=x)$. 

\end{rem}

\begin{lem}\label{dcfdcfa}
Let $a$ be a tuple of a model of {\it DCFA}, and $A$ a subset of that model.
If $tp_{DCF}(a/A)$ is 1-based then $tp(a/A)$ is 1-based.
\end{lem}

{\it Proof}:\\

Analogue to the last statement in the proof of \ref{st1}.2\\
%Let $A=acl(A) \subset B=acl(B) \subset C=acl(C)$ , and let   $b$ be tuple
 %of realisations of $tp(a/A)$.
% Let us suppose that $acl(Bb) \cap  C=B$.
 %By hypothesis $tp_{DCF}(a/A)$ is
 %$1$-based,  therefore  $b$ is independent from $C$ over $B$ in {\it DCF}.
 %As $A=acl(A)$, we have that for all $n \in \na$, $tp_{DCF}(\si^n(b)/A)$ is
 %$1$-based. By \ref{wgr} $tp_{DCF}(b,\si(b),\ldots, \si^n(b)/B)$ is also
 %$1$-based.
 %Hence, $(b,\si(b),\ldots,\si^n(b))$ is {\it DCF}-independant from $C$ over $B$,
 %for every $n \in  \na$.
 %Then for every finite subset $S$ of 
%$acl(Bb)$, $B(S)$ is linearly disjoint from $C$ over $B$ (that is 
%because 
%every such $S$ is such that $B(S)$ contained in $acl_D(B,b,\si(b), 
%\cdots, 
%\si^n(b))$ for some $n$). Thus by definition of linear disjointness 
%$acl(Bb)$ is linearly disjoint from $C$ over $B$. So
%$b$ is {\it DCFA}-independent from $C$ over $B$.\\
$\Box$

\section{An Example}

In this section we  exhibit a set of $\s$-rank 1 which is infinite-dimensional. 
It is known that in {\it DCF} and {\it ACFA}, being finite-dimensional and having finite rank are equivalent and this
 is an important equivalence which had led, for example, to algebraic proofs of the dichotomies for those theories
 (see \cite{jets}).
%In Chapter 3 we will generalize these algebraic proofs to difference-differential fields
%\begin{ej}\label{DCFA617}
%$\sigma(x)=x^n  $
%\end{ej}
%
%Let $ A$ be the set defined by $\sigma(x)=x^n  $; by \ref{DCFA613} we know that $\s(A) \leq \omega  $.
%
%Let us consider the logarithmic derivative: $dx=\frac{Dx}{x}  $.
%
%
%{\bf Claim}: Let $K$ be a differential field , then $K(X)_D =K(X)_d$.
%
%{\it Proof of th Claim} :\\
%
%By induction we will prove that
% 
%(*) $d^kX \in K(X,DX, \cdots, D^kX) $ 
%
%(**)$D^kX  \in  K(X,dX, \cdots, d^kX)  $.
%
%For $ k =0$ it is clear.
%Let us suppose that $d^kX \in K(X,DX, \cdots, D^kX)    $ and $D^kX  \in  K(X,dX, \cdots, d^kX)  $. 
%
%Then $d^{k+1}X=d(d^kX)=\frac{D(d^kX)}{d^kX}  $, by induction hypothesis $D(d^kX)   \in  K(X,DX, \cdots, D^{k+1}X)  
%  $ thus $d^{k+1}X \in   K(X,DX, \cdots, D^{k+1}X)   $.
%
%$D^{k+1}X=(D^kX) d(D^kX) $, then $d(D^kX) \in  K(X,dX, \cdots, d^{k+1}X)   $, then $D^{k+1}X \in  K(X,dX, \cdots, d^{k+1}X)    $.\\
%$\Box$
%
%If $a \in A $, then  $d(\sigma(a))= \frac{n a^{n-1}Da}{a^n}=n da $; thus $d^2(\sigma(a))=d^2(a)  $, 
%and for all $j \in {\mathbf Z}  $, $d^2(\sigma^j(a))=d^2(a) $, so we have that, for all $i \geq 2 $ 
%and for all $j \in {\mathbf Z}$ $(\sigma^j(d^ia))=d^i(a) $.
%
%Let $ K$ be a difference-differential field; if $a$ is a generic of $A$ over $K$ we know that $(a,da,d^2a \cdots ) $
% is a transcendence basis of $K(a)_{\sigma,D}  $ over $K$. 
%
%If $j \geq 2  $, $K(d^ja)_d $ is a difference-differential subfield 
%of $K(a)_{\sigma,D}  $. Hence $\s(A)=\s(a/K) \geq  \omega $, and by our first observation $\s(a/K)=\omega $.
%
%
\begin{ej}\label{DCFA618}
$\sigma(x)=x^2+1 $.
\end{ej}

{\it Let $A$ be the set defined by $\sigma(x)=x^2 +1 $.
Let $A_1=\{x \in A:Dx=0 \}$ and let $A_2=\{x \in A:Dx \neq 0\}$. Then $A_1$ and
$A_2$ are stably embedded and strongly minimal.}\\

{\it Proof}:\\

Let $K = acl(K)$ and  let $a \in A_2$, $a \not\in K$. Let $K_0=K(a)_{\si}$ and $K_{n+1}=K_n(D^{n+1}a)$. Since
$\sigma(D^na)= \sum_{i=0}^n    \binom{n}{i}    D^iaD^{n-i}a $ for $n>0$, each
$K_n$ is a difference field.

Let us write the equation satisfied by $\sigma(D^na)$ over $K_n $ as $\sigma(D^na)=f_n(D^na) $.
Set $f_n^1(X)=f_n(X)$ and
$f_n^{k+1}(X)=(f_n^k)^{\sigma}(f_n(X))$. Then $\sigma^k(D^na)=f_n^k(D^na)$ and
 we have $f_n(X)=  2aX+b_n $ where $b_n=\sum_{i=1}^{n-1} \binom{n}{i}     D^iaD^{n-i}a$ when $n>0$.
Note that $f_0^{k+1}(0)=f_0^k(0)^2+1$, so that $f_0^k(0) \neq 0$ for all $k \geq 0$, and the numbers $f_0^k(0)$
form a strictly increasing sequence.\\\\
Given a difference field $E$, a finite $\si$-stable extension of $E$ is a finite field extension $F$ of $E$ such that
 $\si(F) \subset F$.\\\\
We shall prove the following for $n \geq 1$:\\\\
${\mathbf I_n}$: $K_{n-1}$ contains no finite subset $S$ such that $\si(S)=f_n(S)$, unless $n=1$ in which case $S=\{0\}$. \\
$\mathbf{II}_n$: $K_{n-1}^{alg}(D^na)$ has no proper finite $\sigma$-stable extensions.\\
$\mathbf{III}_n$: Any solution of $\si(x)=x$ in $K_n$ is in $K$. This implies that the solutions of
$\si(x)=(2a)^mx$ in $K_n$ are of the form $c(Da)^m$ where $c\in Fix \si \cap K$; and
the solutions of $\si^k(x)=2^ka \si(a) \cdots \si^{k-1}(a)x$ are of the form $cDa$ where $c \in Fix \si^k \cap K$.\\\\
It will be useful to consider some variants of the first two statements:\\\\
${\mathbf I'_n}$: $K_{n-1}$ contains no finite subset $S$ such that $\si^k(S)=f^k_n(S)$, unless $n=1$ in which case
$S=\{ 0 \}$.\\
$\mathbf{I}''_n$:  $K_{n-1}^{alg}$ contains no finite subset $S$ such that $\si^k(S)=f^k_n(S)$, unless
$n=1$ in which case $S=\{ 0 \}$.\\
$\mathbf{II}'_n$: $K_{n}$ has no proper  $\sigma$-stable finite extensions.\\\\
We will first show some implications between these statements. We suppose $n \geq 1$.\\\\
$\mathbf{I}_n \implies \mathbf{I}'_n$: Replace $S$ by $S\cup \si^{-1}f_n(S) \cup \cdots \cup (\si^{-1}f_n)^{k-1}(S)$.\\\\
$\mathbf{I}'_n  \wedge \mathbf{II}'_{n-1}   \implies   \mathbf{I}''_n $: By 6.1 of \cite{salinas} we know
that $K_0$ has no proper finite $\si$-stable extension so that $\mathbf{II}_0'$ holds. Let $S \subset K_{n-1}^{alg}$ be finite and 
such that $\si^k(S)=f^k_n(S)$ for some $k \in \na$. Then
$K_{n-1}(S)_{\si}=K_{n-1}(S \cup \si(S) \cup \cdots \cup \si^{k-1}(S))$. By $\mathbf{II}'_{n-1}$ $S \subset K_{n-1}$
 and this implies $n=1$, $S=\{0\}$.\\\\
$\mathbf{I}''_n \implies \mathbf{II}_n$: 
Suppose that $L$ is a finite $\sigma$-stable  extension of  $K^{alg}_{n-1}(D^na) $
(by $\mathbf{I}''_n$, $D^na$ is transcendental over $K_{n-1}$).
Then the ramification locus of $L$ over $K_n $ gives us a finite set $S \subset K^{alg}_{n-1}$
 such that $\sigma(S)=f_n(S)$
(see the proof of 4.8 in \cite{salinas}), and this contradicts $\mathbf{I}''_{n}$.\\\\

$\mathbf{II}_n  \wedge \mathbf{II}'_{n-1}   \implies   \mathbf{II}'_n $: 
As before, we know that $\mathbf{II}'_0$ holds.
Let $L$ be a finite $\si$-stable extension
of $K_n=K_{n-1}(D^na)$. By $\mathbf{II}'_{n-1}$, $L \cap K_{n-1}^{alg}=K_{n-1}$.
Hence $[LK^{alg}_{n-1}:K^{alg}_{n-1}K_n]=[L:K_n]=1$ by $\mathbf{II}_n$.\\\\
$\mathbf{I}''_n     \implies   \mathbf{III}_n $:
Suppose there is such a solution  $b \in K_{n}$. Applying $\si$ to $b$ we get $f(X),g(X) \in K_{n-1}[X]$
relatively prime with $g(X)$  monic, such that
$$\frac{f^{\sigma}(f_{n}(D^{n}a) )} {g^{\sigma}(f_{n}(D^{n}a) ) }=
\frac{f(D^{n}a)}{g(D^{n}a)}.$$
Note that, as $f(X)$ and $g(X)$ are relatively prime, $f^{\si}(f_n(X))$, and $g^{\si}(f_n(X))$ are relatively prime:
otherwise, they would have
a common root $\alpha$ in $K$, this implies that $f(\beta)=g(\beta)=0$ for $\beta=\si^{-1}f_n(\alpha)$.

We know that the left side and the right side of the equation should have the same poles, say,
$\alpha_1, \cdots, \alpha_m \in K_{n-1}^{alg}$.
Then $g(X)= \Pi_{i=1}^m (X-\alpha_i)$ and $  g^{\sigma}(f_{n}(X))=\Pi_{i=1}^m(f_{n}(X)-\sigma(\alpha_i) )$
and they have to have the same degree since $f_n$ is linear;
thus $f_{n}(\{\alpha_1, \cdots, \alpha_m\})=\sigma(\{\alpha_1, \cdots, \alpha_m\})$ which contradicts $\mathbf{I}''_n$
unless $n=1$ and $\{ \alpha_1,\cdots, \alpha_m\}=\{0\}$.
The same argument applies to $f$, then $b \in K_1$ and we have $\frac{f(X)}{g(X)}=\alpha X^l$ with 
$\alpha \in K_0$ and $l \in\mathbb{Z}$. 

Inverting $b$, we may assume that $l\in\na$. Then $\alpha$ satisfies $\si(X)=(2a)^lX$.
 Choose $N\geq 0$ minimal such that $\si^N(\alpha)\in K(a)$. Then $\si^N(\alpha)$ satisfies 
$\si(X)=(2\si^N(a))^lX$. If $N>0$, this 
implies that  $\si^N(\alpha)\in K(\si(a))$ and contradicts the 
minimality of $N$. Hence $l=0$. Let $P,Q\in K[X]$ be realitvely prime 
with $Q$ monic and such that $\alpha={P(a)\over Q(a)}$. Then 
$${P^\si(X^2+1)\over Q^\si(X^2+1)}=(2X)^l{P(X)\over Q(X)}.$$ 
Comparing the number of poles and zeroes, we get $deg(Q)=0$ and 
$deg(P)=l$. Hence, if $l=0$, then $\alpha\in K$, and we are done. If 
$l>0$, then $P^\si(f_0(0))=0=P^\si(1)$, hence $P(1)=0$; by induction, 
on then shows that for all $k>0$, $f_0^k(0)$ is a zero of $P$. Since 
the sequence $f_0^k(0)$ is strictly increasing, this is 
impossible. Hence $l=0$. 

%Before we continue we will make an observation.
%
%Let $f(X) \in K_{n-1}[X]$ be a monic polynomial and let $S \subset K_{n-1}^{alg}$ be the set of zeroes of $f(X)$.
%Applying $\si$ we get $\si(f(D^na))=f^{\si}(f_n(D^na))\in K_{n-1}(D^na)$.
%Let $T \subset K_{n-1}^{alg}$ be the set of zeroes of $f^{\si}(f_n(X))$.
%We can write $f(X) = \prod_{\alpha \in S}(X- \alpha)$, then $f^{\si}(f_n(X))= \prod_{\alpha \in S}(f_n(X)-\si(\alpha))$.
%Hence $T=f_n^{-1}(\si(S))$.
%So, if $n>0$, $T=\{\frac{\si(\alpha)-b_n}{2a}:\alpha\in S \}$ and if $n=0$, $T=\{\pm \sqrt{\si(\alpha)-1}: \alpha \in S\}$.\\\\
{\bf Proof of} $\mathbf{I}_1$: \\

$K_0=K(a)_{\sigma}$ and $f_1(X)=2aX$. Suppose that there is a finite subset
$S \subset K_0 \setminus \{ 0\} $
such that $\sigma(S)=f_1(S)$.

%If $|S|>1 $, then 
$(\si^{-1}f_1)$ defines a permutation on $S$, so $(\si^{-1}f_1)^k=id$ for some $k>0$ (if $|S|=1$, $k=1$) and this
implies that $K_0 $ contains a  solution $b$ of
$\si^k(x)=2^ka \sigma(a) \cdots \sigma^{k-1}(a)x$.
%then, dividing this two solutions, $K_0$ has a solution $b$ of $\sigma^ k(X)=X$.
%Let us write this solution $b=\frac{f(a)}{g(a)}$,
%where $f(X),g(X) \in K[X]$ relatively prime and $g(X)$ is monic.
%Then we have $\frac{\sigma^ k(f(a))}{ \sigma^ k(g(a))}=\frac{ f^{\si^k}(f_0^k(a)) }{ g^{\si^k}(f_0^k(a)) }
%=\frac{f(a)}{g(a)}$.
%Since $f$ and $g$ are relatively primes, so are $ f^{\si^k}(f_0^k(X)) $ and $ g^{\si^k}(f_0^k(X))$. But
%$ deg( f^{\si^k}(f_0^k(X)))=2^kdeg(f(X))  $ and $ deg( g^{\si^k}(f_0^k(X)))=2^kdeg(g(X))  $ which is impossible.
%Comparing the degrees we have $2^k(deg(f)-deg(g))=deg(f)-deg(g)$, then $k=0$ which is absurd.
%If $|S|=1$ then $K_0$ contains a solution of $\sigma(X)=2aX$,
Let $N \in {\mathbb N}$ be minimal such that
$\si^N(b) \in K(a)$.
Write $\si^N(b)=\frac{f(a)}{g(a)}$ with $f(X),g(X) \in K[X]$ relatively prime and $g(X)$ monic.
%Note that
%As before, $f^{\si^k}(f_0^k(X))$ and $g^{\si^k}(f_0^k(X))$ are relatively prime in $K[X]$.
%: otherwise, they would have
%a common root $\alpha$ in $K$, this implies that $f(\beta)=g(\beta)=0$ for $\beta=\si^{-k}f_0^k(\alpha)$.\\
Then the equation is 
$$ \frac{ f^{\si^k}(f_0^k(X))}{g^{\si^k}(f_0^k(X))}=2^kf_0^N(X)  \cdots f_0^{N+k-1}(X)   \frac{f(X)}{g(X)}.$$
By minimality of $N$, $\frac{f(a)}{g(a)} \not\in K(\si(a))$, but this is impossible if $N \geq 1$. Thus $N=0$ and the equation is
 $$ \frac{ f^{\si^k}(f_0^k(X))}{g^{\si^k}(f_0^k(X))}=2^kXf_0(X)  \cdots f_0^{k-1}(X)   \frac{f(X)}{g(X)}.$$
As the righthand side and lefthand side of
this equation should have the same poles, $ f^{\si^k}(f_0^k(X))$ and $g^{\si^k}(f_0^k(X))$ are relatively prime,
 and $g$ is monic  we have $g(X)=1$.

%If $N \geq 1$, then  $f^{\si^k}(\si^k(a))=2^k \sigma^{N-1}(a) \cdots \si^{N+k-2}(a)(a^2+1)f(a)$. Thus $f(a) \in K(a^2) \cap K[a]=K[a^2]$.
%So $b \in K[a^2]=K[\si(a)]$ and this contradicts the minimality of $N$. Hence $N=0$.
Then $f^{\si^k}(f_0^k(X))=2^kX f_0(X)\cdots f_0^{k-1}(X)f(X)$. So $2^kdeg(f)=degf+2^k-1$, which implies $deg(f)=1$.
Then $f(X)=cX+d$ with $c,d \in K$. Substituting in the equation we have
$\si^k(c)f_0^k(X)+\si^k(d)=2^kX^2f_0(X) \cdots f_0^{k-1}(X)c+2^kX f_0(X)\cdots f_0^{k-1}(X)d.$ Since the
lefthand side has only even degrees  and the degree of $X f_0(X)\cdots f_0^{k-1}(X)$ is odd we have $d=0$.
Finally, as $f_0^k(0)\neq 0$, the righthand side has no constant term, we obtain $c=0$.\\\\
Now we assume that $\mathbf{I}_k$ holds for all $1 \leq k < n$, where $n \geq 2$.
By what we have shown before the following statements hold:\\

$\mathbf{I}'_k$ for $1 \leq k < n$.\\
$\mathbf{II}_k$ for $1 \leq k < n$.\\
$\mathbf{II}'_k$ for $0 \leq k < n$.\\
$\mathbf{I}''_k$ for $1 \leq k < n$.\\
$\mathbf{III}_k$ for $1 \leq k < n$.\\\\
{\bf Proof of} $\mathbf{I}_{n}$:\\\\
Assume that there is a finite set $S \subset K_{n-1}$ such that $\sigma(S)=f_{n}(S)$.

When  $|S|>1$ we will show that the difference of two distinct elements of $S$ is of the
form $cDa$ with $c \in Fix \si  \cap K$. Indeed, let $a_1,a_2$ be two distinct elements of $S$. 
Reasoning as in the proof of $\mathbf{I}_1$, there is $k>0$ such that  
$\sigma^k(a_i)=f_{n}^k(a_i)$
for $i \in \{1,2\}$.
Then $b=a_1-a_2$ satisfies the equation
$\sigma^k(X)=2^ka \sigma(a) \cdots \sigma^{k-1}(a)X$, and
by $\mathbf{III}_n$ $b=cDa$ with $c \in Fix \si^k \cap K$.

Let $\{a_1,\cdots, a_m \}$ be a cycle in $S$ (i.e. $\si(a_i)=f_n(a_{i+1})$ for $1\leq i \leq m$ and $\si(a_m)=f_n(a_1)$).
Then $\si(a_1+ \cdots + a_m)=2a(a_1+ \cdots a_m)+mb_n$, hence 
$$\si(a_1+\cdots+a_m)=\si(ma_1+(a_2-a_1)+\cdots+(a_m-a_1))=\si(ma_1+d_1Da)$$ for some $d_1 \in Fix\si^k \cap K$.
Hence $\si(a_1)=2aa_1+b_n+c$ where  $c=dDa$ for some $d \in K$.
If $|S|=1$ then $\si(a_1)=2aa_1+b_n$, and we set $c=0$.

We will show that this equation has no solutions in $K_{n-1}$.\\

Case $n \neq 2$: We can write $b_n=2nDaD^{n-1}a+c_1$ where $c_1\in K_{n-2}$.
% where $c_1 \in K_{n-2}$.
We need to show that there is no
$f(X) \in K_{n-2}(X)$ such that $f^{\si}(2aX+b_{n-1})=2af(X)+2nDaX+c'$
 where $c'=c+c_1 \in K_{n-2}$.
If we take the derivative of this equation with respect to $X$ we get $2a(f')^{\si}(2aX+b_{n-1})=2af'(X)+2nDa$, i.e.
$D^{n-1}a$ satisfies the equation $$(f')^{\si}(2aX+b_{n-1})=f'(X)+\frac{nDa}{a}. \eqno{(1)}$$
We also have  $2a(f'')^{\si}(2aX+b_{n-1})=f''(X)$, and by $\mathbf{III}_{n-1}$,  
$f''(D^{n-1}a)=e(Da)^{-1}$ for some $e \in Fix \si \cap K$.
Thus $f''(X)$ is constant, so $f'(X)$ is a polynomial of degree at most  1 in $X$ and 
its leading coefficient is $e(Da)^{-1}$.
Now we look at the degrees in $a$ of the equation (1): $deg_a{b_{n-1}}=0$, and as
$deg_a(e(Da)^{-1})=0$, $deg_a(f'(X))=deg_a(f'(0))=u $.
If $u \leq 0$ we have $deg_a((f')^{\si}(2aX+b_{n-1}))=1$ and if $u>0$ we have $
deg_a((f')^{\si}(2aX+b_{n-1}))=deg(f'(0))=2u$.
 
In both cases, if we compute the degrees in (1) we get a contradiction.\\

Case $n=2$: Then $b_2=2(Da)^2$, and the equation satisfied by $a_1$ is $\si(a_1)=2aa_1+2(Da)^2+dDa$ with $d=0$ if $|S|=1$.
 We will show 
 this equation has no solutions in  $K_1$. If it has there is 
$f(X) \in K_0(X)$ such that $f^{\si}(2aX)=2af(X)+2X^2+dX$.
Taking the second derivative we get $4a^2(f'')^{\si}(2aX)=2af''(X)+4$, i.e. 
$$(f'')^{\si}(2aX)=\frac{f''(X)}{2a}+\frac{1}{a^2} \eqno{(2)}$$
Taking the third derivative we obtain $4a^2(f''')^{\si}(2aX)=f'''(X)$; by $\mathbf{III}_1$ 
$f'''(Da)=e(Da)^{-2}$, which implies $f'''(X)=eX^{-2}$ and therefore $e=0$. Thus $f''(X)=b\in K_0$.
Let $M$ be the smallest natural number such that $\si^M(b) \in K(a)$.
Write $\si^M(b)=\frac{P(a)}{Q(a)}$ where $P$ and $Q$ are relatively prime polynomials over
$K$. Then $$\frac{\si(P(a))}{\si(Q(a))}=\frac{P(a)}{2Q(a)\si^M(a)}+\frac{1}{(\si^M(a))^2}.$$
If $M \geq 1$, by minimality of $M$, $\frac{P(a)}{Q(a)} \not\in K(\si(a))$,
but this is absurd. Hence $M=0$.
So the equation is $$\frac{\si(P(a))}{\si(Q(a))}=\frac{P(a)}{2Q(a)a}+\frac{1}{a^2}.$$
Then the zeroes of $Q^{\si}(X^2+1)$ are contained in the zeroes of $X^2Q(X)$ and  comparing the degrees we have that 
$degQ < 3$. If $Q(0)=0$, then $Q^{\si}(1)=1$, 
hence $Q(1)=0$ and  $Q^{\si}(2)$=0, thus $Q(2)=0$ which is a contradiction. 
If $Q(0) \neq 0$ then the zeroes of $Q^{\si}(X^2+1)$ are contained in the zeroes of $Q(X)$ which
implies $deg(Q)=0$. Hence $Q=1$.
The equation is reduced to $$P^{\si}(a^2+1)=\frac{P(a)}{2a}+\frac{1}{a^2}$$ and comparing the degrees we get a contradiction.

Hence (2) has  no solutions in $K_1$. This finishes the proof of $\mathbf{I}_{n}$.\\
$\Box$\\\\
Let  $a \in A$ and $a \not\in K$.

If $a\in A_1$ then by 6.1 of \cite{salinas} $K(a)_{\si,D}=K(a)_{\si}$ has no finite $\si$-stable extension. Then
all extensions of $\si$ over $acl(Ka)$ are conjugates over $K(a)_{\si,D}$ (see \cite{salinas}),
 thus $qftp(a/K) \vdash tp(a/K)$ and this holds for an arbitrary
 difference-differential field $K$. This means that $tp(a/K)$ is the only non-realized type of $A_1$, and
$A_1$ is strongly minimal. By 6.1 of \cite{salinas}, we know that $A_1$ is trivial.

If $a\in A_2$, by $\mathbf{I}''_{n+1}$, $D^{n+1}a \not\in K^{alg}_n$,  $\td(K_{n+1}/K_n)=1$, and this
implies that $a$ is differentially transcendental over $K$.
%Then $\si(x)=x^2+1 \wedge Dx \neq 0 \vdash tp(a/K)$. 
 Then $tp(a/K)$ is the only non-realized type of $A_2$. As before this implies that $A_2$ is strongly
minimal.

Thus, in particular, $\s(A_2)=1$. Moreover $tp(a/K)$ is trivial, thus 1-based: Indeed, let $a_1,a_2,a_3 \in A_2$ be
such that $a_1 \dnfo_K a_2$,$a_1 \dnfo_K a_3$ and $a_3 \dnfo_K a_2$. We will show that $a_3 \dnfo_K a_1a_2$.
 By 6.1 of \cite{salinas} $tp_{ACFA}(a_3/Ka_1a_2)$ is orthogonal to $Fix \si$ and $tp_{ACFA}(Da_1Da_2\cdots/Ka_1a_2)$
is $Fix\si$-analyzable. Thus, if $a_3 \in acl(Ka_1a_2)$, then $a_3 \in acl_{\si}(Ka_1a_2)$ and by 6.1 of \cite{salinas},
$a_3\in acl_{\si}(Ka_1)$ or $a_3 \in acl_{\si}(Ka_2)$ which is absurd. \\
$\Box$

\cleardoublepage
\chapter{The Dichotomy Theorem}

\label{chap:jets}
As we mentioned in Chapter \ref{chap:prelim}, {\it DCF} and {\it ACFA} satisfy Zilber's dichotomy.
The  original proofs of these dichotomies involve all the machinery of stability. In \cite{jets} Pillay and Ziegler
give proofs of these facts using suitable jet spaces from algebraic geometry, in fact they prove stronger
results which trivially imply the dichotomies.

In the first part of chapter we adapt this method based on jet spaces
to prove an analogue for {\it DCFA}, but with the additional hypothesis of finite-dimensionality.
In the last part we use arc spaces to remove this hypothesis.

\section{Algebraic Jet Spaces}

In this section we list the main  properties of jet spaces over algebraically closed fields of characteristic zero.
As usual, we will suppose all varieties to be absolutely irreducible.
\begin{defi}\label{J11}

Let $K$ be an algebraically closed field , and  let
$V \subset \mathbb{A}^n$ be a variety over $K^n $; let $a$ be a non singular point of $V$.
Let ${\mathcal O}_{V,a}$ be the local ring of $V$ at $a$ and let ${\mathfrak M_{V,a}}$ be its maximal ideal.
%\begin{enumerate}
%The maximal ideal of $V$ in $a$ is ${\mathcal M}_{V,a}= \{f \in K[V]: f(a) =0   \}$, where $K[V]$ is the coordinate
 %ring of $V$.
%\item
Let $m>0 $. The $m$-th jet space of $V$ at $a$, $J^m(V)_a$, is the dual space of the $K$-vector space
$ {\mathfrak M}_{V,a}/  {\mathfrak M}^{m+1}_{V,a}   $.
%\end{enumerate}
\end{defi}

\begin{nota}
If the variety $V$ is $\mathbb{A}^n $, we write ${\mathfrak M}_a $ instead of ${\mathfrak M}_{V,a} $.\
\end{nota}
The following is proved in \cite{jets} (Fact 1.2).

\begin{fac}\label{J12}
Let $U,V  $ be irreducible varieties of $K^n$, $a \in V \cap U$. If $J^m(V)_a = J^m(U)_a$ for all $m > 0 $, then $V=U  $. 
\end{fac}

\begin{prop}\label{J13}
Let $V$ be an  variety, $a$ a non-singular point of $V$. Let ${\mathcal O}_a $ be the
local ring of $V$ at $a$, and ${\mathfrak M}_{V,a} $ its maximal ideal.
Let ${\mathcal M}_{V,a}=\{f \in K[V]:f(a)=0 \}$ be the maximal ideal of the coordinate ring of $K[V]$ of $V$.
Then ${\mathcal M}_{V,a}/{\mathcal M}^m_{V,a}   $ and ${\mathfrak M}_{V,a}/{\mathfrak M}_{V,a}^m$
are isomorphic K-vector 
spaces for all $m \in {\mathbb N}$.

\end{prop}

{\it Proof} :\\

This is a consequence from the fact that  ${\mathcal M}^i_{a,V} \cap K[V]={\mathfrak M}^i_{a,V}$ for all $i$. 
(cf Proposition 2.2 in \cite{eisen}) .
%We know that $\varphi:K[V]\to {\mathcal O}_a$ defined by $f \mapsto f/1 $ is a monomorphism.
%Also ${\mathfrak M}_{a,V}={\mathcal O}_a {\mathcal M}_{V,a} $ and 
%${\mathfrak M}_{a,V}^m={\mathcal O}_a {\mathcal M}^m_{V,a}  $.

%Let $f,g \in K[V] $, if $\varphi (f-g) \in {\mathfrak M}_{a,V}^m$ then
%  $f-g \in {\mathcal  M}^m_{V,a} $, thus the map
% $\varphi^*:K[V]/{\mathcal M}^m_{V,a}  \to{\mathcal O}_a /{\mathfrak M}_{a,V}^m $ is a monomorphism.

%Let $\frac{g}{f} \in {\mathfrak M}_{a,V}$, we shall find $h \in {\mathcal M}_{V,a}  $ such that
% $h-\frac{g}{f} \in {\mathcal O}_a {\mathcal M}^m_{V,a}  $.

%Let $\alpha =f(a) \neq 0 $, let $h'=\alpha^m\frac{g}{f} +(\alpha-f)^m(\alpha^m-\frac{g}{f} ) $.

%Then $h'=\alpha^m \frac{g}{f} + \alpha ^m(\alpha-f)^m-\frac{g}{f}(\alpha-f)^m = \frac{g}{f}(\alpha^m -
%(\alpha-f)^m)+\alpha^m(\alpha-f)^m
%= \frac{g}{f}(f\sum_{i=1}^m\frac{m!}{i!(m-i)!}(-1)^if^{i-1} )+\alpha^m(\alpha-f)^m =
%g(\sum_{i=1}^m\frac{m!}{i!(m-i)!}(-1)^if^{i-1} )+\alpha^m(f-\alpha^m)$.

%Then $h' \in K[V] $ and $h'(a) =0$, thus $h=\frac{h'}{\alpha^m} \in {\mathcal M}_{V,a} $ and
%$\frac{g}{f}-h=\frac{1}{\alpha^m}(\alpha-f)^m(\alpha^m-\frac{g}{f}) \in {\mathcal O}_a {\mathcal M}^m_{V,a}  $.
%Thus $\varphi^*({\mathcal M}_{V,a}/{\mathcal M}^m_{V,a})={\mathfrak M}_{a,V}/{\mathfrak M}_{a,V}^m  $\\
$\Box$

The following fact is proved in \cite{sha}, Chapter II, section 5.

\begin{fac}\label{J14}
 Let $U,V$ be two irreducible varieties defined over $L \subset K$. Let $f: U \to V $ be a finite morphism, and let $b \in V $.
If $f$ is unramified at $b$, then, for any $a \in f^{-1}(b)$ and for any positive integer $m$, the homomorphism
${\bar f}:{\mathcal O}_{V,b}/{\mathfrak M}_{V,b}^m \to {\mathcal O}_{U,a}/{\mathfrak M}_{U,a}^m $ induced by $f$ is an isomorphism.
\end{fac}

\begin{prop}\label{J15a}
Let $U,V$ be two irreducible varieties defined over $L \subset K$. Let $f: U \to V $ be a dominant generically finite-to-one morphism.
Let $a$ be a generic of $U$ over $L$. Then $f$ induces an isomorphism of $K$-vector spaces between $J^m(U)_a $ and $J^m(V)_{f(a)}$.
\end{prop}

{\it Proof} :\\

Since $f$ is separable (as we work in characteristic zero), and since $f $ is dominant and $f^{-1}(f(a)) $ is finite, 
$U$ and $V$ are irreducible and their dimensions are equal, thus $f$ is
unramified at $f(a)$. By \ref{J14}, $f$ induces an isomorphism between
${\mathcal O}_{V,f(a)}/{\mathfrak M}_{V,f(a)}^{m+1} $ and ${\mathcal O}_{U,a}/{\mathfrak M}_{U,a}^{m+1}$;
 whose restriction to ${\mathfrak M}_{V,f(a)}/{\mathfrak M}_{V,f(a)}^{m+1} $ is an isomorphism between 
${\mathfrak M}_{V,f(a)}/{\mathfrak M}_{V,f(a)}^{m+1}  $ and ${\mathfrak M}_{U,a}/{\mathfrak M}_{U,a}^{m+1}$.
Then, by \ref{J11}, $f$ induces an isomorphism between  $J^m(U)_a $ and $J^m(V)_{f(a)}$.\\
$\Box$\\

The following lemma (2.3 of \cite{jets}) allows us to consider jet spaces as algebraic varieties.

\begin{lem}\label{J111} 
 Let $K $ be an algebraically closed field and $V$ a subvariety of $K^n $ , let  $ m \in {\mathbf N}$
and let $\mathcal D$ be the set of  operators 
\[ \frac{1}{s_1! \cdots s_n!} \frac{\partial^s }{\partial x^{s_1}_{1} \cdots  \partial x^{s_n}_{n}}\]
where $0< s <m+1$ and $s =s_1+ \cdots +s_n, \; s_i \geq 0  $. 

Let $a=(a_1,\cdots,a_n) \in V$; and let $d=| {\mathcal D}| $.

Then we can identify $J^m(V)_a   $ with 
$$\{(c_h)_{h \in {\mathcal D}} \in K^d: \sum _{h \in {\mathcal D}} D P(a) c_h =0, \; P \in I(V)   \}    . $$

\end{lem}

{\it Proof} :\\

Let $p:K[X] \longrightarrow K[V] $ such that $Ker(p)=I(V) $ ; then $p^{-1}({\mathcal M}_{a,V}   )= {\mathcal M}_a $,
and  $p^{-1}({\mathcal M}^{m+1}_{a,V}   )= {\mathcal M}_a^{m+1} +I(V)  $.
This gives us   the following short exact sequence:
$$0 \longrightarrow (I(V)+{\mathcal M}_a^{m+1})/{\mathcal M}_a^{m+1} \longrightarrow {\mathcal M}_a / {\mathcal M}_a^{m+1} \longrightarrow {\mathcal M}_{a,V}  /  {\mathcal M}^{m+1}_{a,V} \longrightarrow 0  $$

We proceed to describe the dual space of ${\mathcal M}_a / {\mathcal M}_a^{m+1}   $:% For every monomial
%$(X - a  )^{s}   $ in $(X_1-a_1), \cdots, (X_n-a_n)   $ of total degree greater than 0 and less than $m+1$. 
%where $s$ us the tuple formed by the degrees of
%the $(X_i-a_i) $. 
%Let  $u_{s} $ be the linear map which assigns 1 to $  (X -  a  )^{s}$.
The monomials $(X-a)^s=(X-a_1)^{s_1}\cdots (X-a_n)^{s_n}$ with $1 \leq s_1+\cdots + s_n=s \leq m$ form
a basis for ${\mathcal M}_a / {\mathcal M}_a^{m+1}$, and for each $s$ 
we have a $K$-linear map $u_s$ which assigns 1 to $  (X -  a  )^{s}$ and $0$ to the other monomials.
The maps $u_s$ form a basis for the dual of ${\mathcal M}_a / {\mathcal M}_a^{m+1}$.

 %Since the $ (X - a )^{s}   $ contain a $K$-basis for  ${\mathcal M}_a / {\mathcal M}_a^{m+1}   $, then the $u_{s}  $ contain a $K$-basis for its dual.

Thus, the dual $ J^m(V)_a  $   of $ {\mathcal M}_{a,V}  /  {\mathcal M}^{m+1}_{a,V}$,
 consists of those linear maps 
$u:{\mathcal M}_a / {\mathcal M}_a^{m+1}   \longrightarrow K  $ that take the value 0 on 
$(I(V)  + {\mathcal M}_{a,V})  /  {\mathcal M}^{m+1}_{a,V} $.

Let $f(X) \in  K[X] $; applying Taylor's formula we can write, modulo $ {\mathcal M}^{m+1}_{a,V}   $,
 $$ f(X) =f(a) + \sum_{1 \leq|s| \leq m } D_{s} f(a)  (X -  a  )^{s} ,    $$

where  $$D_{s}=\frac{1}{s_1! \cdots s_n!} \frac{\partial^{s} }{\partial X^{s_1}_{1} \cdots  \partial X^{s_n}_{n}    }  $$

If $u =\sum_{s} c_{s} u_{s}  $, then $u$ vanishes on $ (I(V)+{\mathcal M}_a^{m+1})/{\mathcal M}_a^{m+1}  $ if and only if for every
$P( X) \in I(V)  $, we have       $$\sum_{ 1 \leq |s| \leq m   }    D_{s} P( a) c_{s}=0  .$$
$\Box$

\section{Jet Spaces in Differential and Difference Fields}

In this section we study jet spaces of varieties over differential fields and difference fields.
We recall the concepts of $D$--modules and $\sigma$-modules
(see  \cite{jets}).

\begin{defi}\label{J22}
Let $(K,D)$ be a differential field, and let $V$ be a finite-dimensional $K$-vector space.
We say that $(V,D_V)$ is a $D$-module over $K$ if $D_V$ is an additive endomorphism of $V$ such that,
for any $v \in V $ and $c \in K $, $D_V(cv)=cD_V(v)+(Dc)v $.
\end{defi}

\begin{lem}\label{J23} \textup{(\cite{jets}, 3.1)}
Let $(V,D_V)$ be a $D$-module over the differential field  $(K,D)$.
Let $(V,D_V)^{\sharp} =\{v \in V:D_V v=0\}$ . Then $(V,D_V)^{\sharp} $ is a finite-dimensional $\mathcal C$-vector space.
Moreover, if $(K,D)$ is differentially closed, then there is a $\mathcal C$-basis of $(V,D_V)^{\sharp} $ which
is a $K$-basis of $V$. (Thus every ${\mathcal C}$-basis of  $(V,D_V)^{\sharp}$ is a $K$-basis of $V$)
\end{lem}

\begin{defi}\label{J21}
A $D$-variety is an algebraic variety  $V \subset \mathbb{A}^n $ with an algebraic section 
$s:V \to \tau_1(V)$ of the projection $\pi:\tau_1(V) \to V$. 
Then, by  \ref{P224}, $(V,s)^{\sharp}=\{x \in V: Dx=s(x)   \} $ is Zariski-dense in $V$.
We  shall write $V^{\sharp}$  when $s$ is understood.
\end{defi} 
\begin{prop}
A finite-dimensional affine differential algebraic variety is differentially birationally equivalent to
 a set of the form $(V,s)^{\sharp}=\{x \in V: Dx=s(x)   \} $ where $(V,s)$ is a $D$-variety.\\
\end{prop}

%Let us recall that an ideal of a differential ring is called a differential ideal if
%it is closed under the derivation. The quotient  of a differential ring by a differential ideal has
%a natural structure of differential ring.

\begin{rem}\label{J24} Let $V \subset {\mathbb{A}}^n$ be a variety defined over $K$.

\begin{enumerate}
%\item If we write $s=(s_1, \cdots,s_n)$, where $s_i:V \longrightarrow \mathbb{A} $, then $(V,s) $ defines a
%differential variety if and only if for every $F \in I(V) $
%$$\sum_{i=1}^n \frac{\partial F}{\partial X_i}(X)s_i(X)+F^D(X) \in I(V)$$
\item Given a $D$-variety $(V,s) $ , we can extend  the derivation $D$ to the field of rational functions
of  $V$ as follows:\\

 If $f \in {\mathcal U}(V)$, then we define $Df= \sum \frac{\partial f} {\partial X_i}s_i+f^D $.

\item If $a \in V^{\sharp} $ and  $f \in {\mathfrak M}_{V,a}  $, then
$Df(a)=\sum \frac{\partial f} {\partial X_i}s_i(a)+f^D(a)=J_f(Da)+f^D(a)=D(f(a))=0$.
 Thus  ${\mathfrak M}_{V,a}$ and  ${\mathfrak M}^{m+1}_{V,a} $ are differential ideals of ${\mathcal O}_{V,a} $, 
so it gives
${\mathfrak M}_{V,a} /{\mathfrak M}^{m+1}_{V,a} $ a structure of $D$-module over $\mathcal U $. Defining
$D^*:J^m(V)_a \to J^m(V)_a$ by $D^*(v)(F)=D(v(F))-v(D(F))$ for $v \in J^m(V)_a$ and
 $F \in {\mathfrak M}_{V,a} /{\mathfrak M}^{m+1}_{V,a} $,  gives $J^m(V)_a$ a structure of $D$-module.

\end{enumerate}
\end{rem}

\begin{defi}
Let $(K,\sigma)$ be a difference field. A $\sigma$-module over $K$ is a finite-dimensional $K$-vector space $V$
together with an additive automorphism $\Sigma: V \to V$, such that, for all $c\in K$ and $v \in V$,
$\Sigma(cv)=\sigma(c)\Sigma(v)$.
\end{defi}

\begin{lem}\label{J25}\textup{(\cite{jets}, 4.2)}
Let $(V,\Sigma)$ be a $\sigma$-module over the difference field $(K,\sigma)$.
Let $(V,\Sigma)^{\flat}=\{v \in V: \Sigma(v)=v \}$. Then $(V,\Sigma)^{\flat}$ is a finite-dimensional $Fix \sigma$-vector space.
Moreover,
if $(K,\sigma)$ is a model of \it{ACFA}, then there is a $Fix \sigma$-basis of $(V,\Sigma)^{\flat}$ which is a $K$-basis
of $V$.(Thus every $Fix \si$-basis of  $(V,\Sigma)^{\flat}$ is a $K$-basis of $V$)
\end{lem}

\begin{rem}
Let $(K,\sigma)$ be a model of \it{ACFA}. Let $V,W$ be two irreducible algebraic affine varieties over $K$ such that
$W \subset V \times V^{\sigma}$, and assume that the projections from $W$ to
$V$ and $V^{\sigma}$ are dominant and generically finite-to-one.
 Let $(a,\sigma(a))$ be a generic point of $W$ over $K$. Then, by \ref{J15a}, $J^m(W)_{(a,\sigma(a))}$ induces an isomorphism $f$
of $K$-vector spaces between $J^m(V)_a$ and $J^m(V)_{\sigma(a)}$. We have also that $(J^m(V)_a,f^{-1}\sigma)$ is a
$\sigma$-module over $K$.
\end{rem}

\section{Jet Spaces in Difference-Differential Fields}

In this section we describe the jet spaces of finite-dimensional varieties defined over difference-differential fields, and we state
the results needed to prove our main theorem \ref{J312}. Finally we give two corollaries: the first is the weak dichotomy,
and the second is an application to quantifier-free definable groups.

We start with the definition of a $(\sigma,D)$-module.

\begin{defi}
Let $(K,\sigma,D)$ be a difference-differential field. A $(\sigma,D)$-module over $K$ is a finite-dimensional $K$-vector space $V$
equipped with an additive automorphism $\Sigma: V \to V$ and an additive endomorphism $D_V:V \to V$, such that
$(V,D_V)$ is a $D$-module over $K$, $(V, \Sigma)$ is a $\sigma$-module over $K$ and for all $v \in V$ we have
 $\Sigma(D_V(v))=D_V(\Sigma(v))$.
\end{defi}

The key point of our proof of \ref{J312} is the following lemma.

\begin{lem}\label{dens1}
Let $(V,\Sigma,D_V)$ be a $(\sigma,D)$-module over the difference-differential field $(K,\sigma,D)$.
Let $(V,\Sigma,D_V)^{\natural}=\{v \in V: D_V(v)=0 \wedge \Sigma(v)=v \}$
\textup{(}we shall write $V^{\natural}$ when $D_V$ and $\Sigma$ are understood\textup{)}.
Then $V^{\natural}$ is a $(Fix \sigma \cap {\mathcal C})$-vector space. Moreover, if $(K,\sigma,D)$ is a model
of \it{DCFA}, there is a $(Fix \sigma \cap {\mathcal C})$-basis of $V^{\natural}$ which is a $K$-basis of $V$.
(Thus every $(Fix \si) \cap {\mathcal C}$-basis of  $(V)^{\natural}$ is a $K$-basis of $V$)
\end{lem}

{\it Proof} :\\

It is clear that $V^{\natural}$ is a $(Fix \sigma \cap {\mathcal C})$-vector space.
By \ref{J23} and \ref{J25} it is enough to prove that there is a $(Fix \sigma \cap {\mathcal C})$-basis of $V^{\natural}$ which is
a $\mathcal C$-basis of $V^{\sharp}$.

Let $\{v_1, \cdots, v_k   \}  $ be a $\mathcal C $-basis of $ V^{\sharp} $, then $\{\Sigma(v_1), \cdots,\Sigma( v_k)\}$ is a
$\mathcal C $-basis of $V^{\sharp}$.
Let $A$ be the invertible $k \times k $ $\mathcal C  $-matrix such that $[\Sigma(v_i)]^t=A [v_i]^t  $.

Let $\{u_1, \cdots, u_k\} $ be a $\mathcal C  $-basis of $V^{\sharp}$. Then there exists an invertible $k \times k  $
$\mathcal C $-matrix $B$ such that $[u_i]^t=B[v_i]^t$; applying $\Sigma$ we get
$[\Sigma(u_i)]^t=\sigma(B)[\Sigma(v_i)]^t=\sigma(B)A[v_i]^t$. Thus $\{u_1, \cdots, u_k\} $ is in   $V^{\natural}$ if and only if
$B=\sigma(B)A $. Since $({\mathcal C},\sigma) \models ACFA  $,
the system $ X=\sigma(X)A $, where $X$ is an invertible $k\times k$ matrix,
 has a solution in $\mathcal C $. So we can suppose that
$\{u_1, \cdots, u_k\}  $ is in $V^{\natural}$.

Let $v \in  V^{\natural}$, and let  $\lambda_1, \cdots \lambda_k \in {\mathcal C} $ such that $v=\lambda_1 u_1+\cdots + \lambda_k u_k $.
Then $ v=\sigma(\lambda_1) u_1+\cdots + \sigma( \lambda_k) u_k$, thus $\lambda_i \in Fix\sigma  $ for $i=1, \cdots , k$.
Hence $\{u_1, \cdots, u_k\}  $ is a $(Fix \sigma \cap {\mathcal C}) $-basis of $V^{\natural}.$\\
$\Box$

\begin{nota}
Let $({\mathcal U},\sigma,D )$ be a saturated model of {\it DCFA}. Let $K=acl(K)$ be a difference-differential subfield of $ \mathcal U$,
and let $a \in {\mathcal U}^n $ such that $K(a)_D=K(a)$ and $\sigma(a)\in K(a)^{alg} $.

Let $V$ be the locus of $a$ over $K$, and let $W$ be the locus of $(a,\sigma(a))$ over $K$.
Then $V^{\sigma}$ is the locus of $\sigma(a) $ over $K$ and the projections $\pi_1:W \longrightarrow V $ and
$\pi_2:W \longrightarrow V^{\sigma} $ are generically finite-to-one and dominant.

We set:

 $\pi^*_1: K[V] \longrightarrow K[W]$,  $F \longmapsto F \circ \pi_1 $.

 $\pi^*_2: K[V^{\sigma}] \longrightarrow K[W]$, $G \longmapsto G \circ \pi_2 $.

$\overline{\pi^*_1}:{\mathfrak M}_{V,a}/ {\mathfrak M}^{m+1}_{V,a}  \longrightarrow  {\mathfrak M}_{W,(a,\sigma(a))}/ {\mathfrak M}^{m+1}_{W,(a,\sigma(a)) }$ the map induced by $\pi^*_1$

$\overline{\pi^*_2}:{\mathfrak M}_{V^{\sigma},\sigma(a)}/ {\mathfrak M}^{m+1}_{V^{\sigma},a}  \longrightarrow  {\mathfrak M}_{W,(a,\sigma(a))}/ {\mathfrak M}^{m+1}_{W,(a,\sigma(a)) } $ the map induced by $\pi^*_2$

$ \pi_1': J^m(W)_{(a,\sigma(a))} \longrightarrow J^m(V)_a $, $w \longmapsto w \circ    \overline{\pi^*_1} $.

$ \pi_2': J^m(W)_{(a,\sigma(a))} \longrightarrow J^m(V^{\sigma})_{\sigma(a)} $,  $w \longmapsto w \circ    \overline{\pi^*_2} $.\\
With respect to the extension of $D$ to the coordinate rings, $\pi^*_1$ and $\pi^*_2$ are differential
homomorphisms.
By \ref{J15a} $ \pi_1' $ and $ \pi_2' $ are isomorphisms of $\mathcal U$-vector spaces.

Let $f:J^m(V)_a  \longrightarrow   J^m(V^{\sigma})_{\sigma(a)}  $ be the $\mathcal U $-isomorphism defined by $f= \pi_2' \circ  (\pi_1')^{-1} $.

Since $Da \in K(a)$ there is a rational map  $s:V \to{\mathcal U}^n $ such that $s(a)=Da $ and $(V,s)$ is a $D$-variety.
 By construction $(V^{\sigma},s^{\sigma }) $ and $(W,(s,s^{\sigma})) $ are also $D$-varieties.
\end{nota}

\begin{lem}
$(J^m(V)_a, f^{-1}\sigma,D^*)$ is a $(\sigma,D)$-module.

\end{lem}

{\it Proof} :\\

All we need to prove is that $D^*$ commutes with $f^{-1}\sigma$.
Since $f= \pi'_2 \circ (\pi'_1)^{-1}$ and $\pi'_1, \pi'_2$ are isomorphisms, 
and since $\sigma$ commutes with $D^*$, it is enough to prove
that $D^*$ commutes with $\pi'_1$ and $\pi'_2$.

Let $w \in  J^m(W)_{(a,\sigma(a))}  $ and
$F \in  {\mathfrak M}_{V,a}/ {\mathfrak M}^{m+1}_{V,a }   $.

We want to prove that $D^* (\pi'_1(w))(F)= (\pi'_1 \circ D^* (w))(F)$.
We have
$D^*(\pi_1'(w))(F)=D^*(w \circ \overline{\pi_1^*})(F)=
 D((w \circ \overline{\pi_1^*})(F))-w \circ \overline{\pi_1^*}(D(F))$.

On the other hand  $\pi'_1(D^*(w))(F)=(D^*(w) \circ \overline{\pi_1^*}) (F)=
 D(w(\overline{\pi_1^*} (F)))-w(D_W(\overline{\pi_1^*}(F))$.

But clearly $ D((w \circ \overline{\pi_1^*})(F)) =  D(w(\overline{\pi_1^*} (F)))$ and
 $w \circ \overline{\pi_1^*}(D_V(F)) = w(D_V(\overline{\pi_1^*})(F))$.\\
The proof is similar for $\pi'_2$.\\
$\Box$

\begin{lem}\label{J38}
Let $K \subset K_1=acl(K_1)  $.
Let $V_1$ be the $(\si,D)$-locus of $a$ over $K_1  $, and let $c$ be the field of definition of $V_1$.
Then $c \subset Cb(qftp(a/K_1)) \subset acl(K,c)  $.

\end{lem}

{\it Proof} :\\

Clearly $ c \subset Cb(qftp(a/K_1))  $. We know that $a \dnfo_{K,c} K_1  $ in {\it DCF} ,
 also $\sigma^i(D^ja) \subset K(a)^{alg} $; then $acl_{\si,D}(K,a) \dnfo_{K,c} K_1 $ in {\it ACF},
 thus $ Cb(qftp(a/K_1)) \subset acl(K,c).  $\\
$\Box$

\begin{rem}\label{inter}
If we replace $a$ by $(a,\si(a),\cdots,\si^m(a))$ for $m$ large enough, 
$c$ and $Cb(qftp(a/K_1))$ will be interdefinable over $K$(choose $m$
for which the Morley rank of $tp_{DCF}(\si^m(a)/K(a,\cdots,\si^{m-1}(a)))$ is minimal and
for which the Morley degree of $tp_{DCF}(\si^m(a)/K(a,\cdots,\si^{m-1}(a)))$ is minimal) .
\end{rem}

\begin{lem}\label{J39}
Let $K \subset K_1 =acl(K_1)  $. Let $V_1$ be the locus of $a$ over $K_1$.
Then $J^m(V_1)_a$ is a $(\sigma,D)$-submodule of $J^m(V)_a$.
\end{lem}

{\it Proof} :\\

Clearly $J(V_1)_a$ is a $D$-submodule of $J^m(V)_a$.
Let $W_1$ be the locus of $(a,\sigma(a))$ over $K_1$. Let $f_1$ be the isomorphism between  $ J^m(V_1)_a   $ and
 $  J^m(V^{\sigma}_1)_{\sigma(a)} $ induced by
the projections from $W_1$ onto $V_1$ and $(V_1)^{\sigma }$; since these projections are the restrictions 
of the projections
from $W$ onto $V$ and $V ^{\sigma}$,
$f_1 \subset f$. So $J^m(V_1)_a$ is a $\sigma$-submodule of $J^m(V)_a$.\\
$\Box$

\begin{teo}\label{J312}
Let $({\mathcal U},\si,D) $ be a saturated model of {\it  DCFA} and let  $K=acl(K) \subset {\mathcal U} $.
Let $tp(a/K) $ be finite-dimensional (i.e. $\td(K(a)_{\si,D}/K) < \infty$). Let $b$ be such that
 $b=Cb(qftp(a/acl(K,b)))  $. Then $tp(b/acl(K,a))  $ is almost-internal to $ Fix \sigma \cap {\mathcal C}.$

\end{teo}

{\it Proof} :\\

By assumption, $  trdg(K(a)_{\sigma,D}/K )$ is finite. 
Enlarging $a$, we may assume that $a$ contains a transcendence basis of $K(a)_{\si,D}$ over $K$.
Then $\si(a),Da \in K(a)^{alg}$ and $D^2(a) \in K(a,Da)$. Hence we may assume that
$Da \in K(a)$.
%So we can replace $a$ by a transcendence
% basis of $ K(a)_{\sigma,D} $ over $K$ which contains $a$. Then we have $\sigma(a),Da \in K(a)^{alg}  $.
 %Now we replace $a$ by $(a,Da)$, so we obtain $\sigma(a) \in K(a)^{alg}$
 %and $Da \in K(a)  $.

 Let $V$ be the locus of $a$ over $K$ ,  $W$ the locus of $ (a,\sigma(a) ) $ over $K$, 
thus $V^{\sigma}  $ is the locus of $\sigma(a)$ over $K$.

%By \ref{34} $L^m(V)_a  $ is a finite dimensional $({\mathcal C}\cap Fix\sigma) $-vector space, Zariski-dense in $J^m(V)_a  $.

Let $V_1  $ be the locus of $a$ over $acl(K,b)$; let $b_1$ be the field of definition of $V_1$. By \ref{J38}
 $b \in acl(K,b_1)$.

%By \ref{34} $L^m(V_1)_a  $ is a  subspace of $L^m(V)_a  $ which is Zariski-dense in $J^m(V_1)_a  $.

By \ref{dens1} for each $m>1$ there is a $(Fix \sigma \cap {\mathcal C})$-basis of $J^m(V)_a^{\natural}$ which is a
 $\mathcal U$-basis
of $J^m(V)_a$. Choose such a basis  $d_m$  such that
 $d=(d_1,d_2, \cdots) \dnfo _{K,a} b$. Then for each $m$ we have an isomorphism between
$J^m(V)_a^{\natural}$ and $({\mathcal C} \cap Fix\sigma)^{r_m}  $ for some $r_m  $.
Thus the image of $J^m(V_1)_a^{\natural}$ in $( Fix\si  \cap{\mathcal C}  )^{r_m}  $
is a  $ ( Fix\si  \cap{\mathcal C})$-subspace of $ ( Fix\si  \cap{\mathcal C} ) ^{r_m}  $ and
therefore it is defined over some tuple
 $e_m \subset    Fix\si  \cap{\mathcal C}    $; let $e=(e_1,e_2, \cdots)  $.
If $\tau$ is an automorphism of $({\mathcal U},\sigma, D)$ fixing $K,a,d,e$, then $J^m(V_1)_a=\tau(J^m(V_1)_a)$;
 on the other hand,
%by \ref{111},
 $\tau(J^m(V_1)_a)=J^m(\tau(V_1))_a$, thus for all $m>1$, $J^m(V_1)_a =J^m(\tau(V_1))_a$ and by \ref{J12} $\tau(V_1)=V_1$,
thus $\tau(b_1)=b_1$ which implies that $b_1 \in dcl(K,a,d,e)$.
Hence $b \in acl(K,a,d,e)$. Since $e \subset Fix \si \cap \mathcal{C}$ and
$d \dnfo_{Ka} b$, this proves our assertion.\\
$\Box$

As in  \cite{camp}, we deduce the dichotomy theorem.
\begin{cor}\label{J313}
If $tp(a/K) $ is of $\s$-rank 1 and finite-dimensional, then it is either 1-based or non-orthogonal to
 $  Fix\si  \cap{\mathcal C}    $.
\end{cor}

{\it Proof}:\\

We supress the set of parameters. Let $p=tp(a)$. If $p$ is not 1-based there is a tuple of
realizations $d$ of $p$ and a tuple $c$ such that 
$c=Cb(qftp(d/c)) \not\subset acl(d)$. Then $tp(c/d)$ is non-algebraic and
by \ref{J312} it is almost-internal to $Fix \si \cap {\mathcal C}$. As $tp(c/d)$ is $p$-internal
we have $p \not \perp Fix \si \cap {\mathcal C}$.\\
$\Box$

We conclude with an application to definable groups of {\it DCFA}. 
We need quantifier-free versions of \ref{needjet1} and \ref{needjet2}.

\begin{lem}\label{J314}
Let $M$ be a simple quantifier-free stable structure which eliminates imaginaries. Let $G$ be a connected group, 
quantifier-free  definable in $M$ defined over
$A=acl(A) \subset M$. 
Let $c\in G$ and  let $H$ be the left stabilizer of $p(x)=qftp(c/A)$. Let $a \in G$ and $b$
 realize a non-forking extension of $p(x)$ to  $acl(Aa)$.
Then $aH$ is interdefinable over $A$  with $Cb(qftp(a \cdot b/A,a))$.
Likewise with right stabilizers and cosets in place of left ones, and $b \cdot a$ instead of $a \cdot b$.
\end{lem}

{\it Proof} :\\

Let $q$ be the quantifier-free type over $M$ which is the non-forking extension of $p$.
Then $aq$ is the non-forking extension to $M$ of $qftp(a \cdot b/Aa)$.
So we must prove that for every automorphism $\tau \in Aut(M/A)$, $\tau(aH)=aH$ if and only if
$\tau(aq)=aq$.

Since $q$ is $A$-definable, $\tau(q)=q$, and $\tau(aq)=\tau(a)\tau(q)=\tau(a)q$.
Thus $\tau(aq)=aq$ if and only if $a^{-1}\tau(a)q=q$. But $H= \{x \in G:xq=q \}$, then
$a^{-1}\tau(a) \in H$ if and only if $\tau(a)H=aH$, and as $H$ is $Aa$-definable, $\tau(aH)=\tau(a)H$.\\
$\Box$

\begin{lem}\label{J314i}
Let $M$ be a simple quantifier-free stable structure which eliminates imaginaries.
Let $G$ be a connected group, quantifier-free  definable in $M$ defined over
$A=acl(A) \subset M$. Let $c \in G$, let $H$ be the left stabilizer of $qftp(c/A)$ and let $a \in G$ be a generic over
$A \cup \{c\}$. Then  $Hc$ is interdefinable with $Cb(qftp(a/A,c \cdot a))$ over $A \cup \{a \}$
\end{lem}

{\it Proof} :\\

We may assume $A=\emptyset$. Let $p=qftp(c/A)$.
We know that $H$ is the right stabilizer of $p^{-1}$, on the other hand, since $a$ is a generic of $G$ we have
$c \dnfo c \cdot a$. 
By \ref{J314}, $Hc \cdot a$ is interdefinable  with $Cb(qftp(c^{-1}(c \cdot a)/c \cdot a))$. Since $H$ is
$\emptyset$-definable, $Hc \cdot a$ is interdefinable with $Hc$ over $a$.\\
$\Box$

\begin{cor}\label{J315}
Let $({\mathcal U},\si,D)$ be a model of {\it DCFA}, and let $K=acl(K) \subset {\mathcal U}$.
Let $G$ be a finite-dimensional quantifier-free definable group, defined over $K$. Let  $a \in G$ and
let  $p(x)=qftp(a/K)$. Assume that $p$ has trivial stabilizer. Then $p$ is internal to $Fix\si \cap{\mathcal C} $.
\end{cor}

{\it Proof} :\\

Let $b \in G$ be a generic over $K \cup \{a\}$. By \ref{J314i} $a$ is interdefinable with $Cb(qftp(b/K,a \cdot b))$
over $K \cup \{b\}$ and by \ref{J312}, $tp(Cb(qftp(b/K,a \cdot b))/K,b)$ is internal to ${\mathcal C} \cap Fix\sigma$.
Thus $tp(a/K,b )$  is internal to $ Fix\si  \cap{\mathcal C}   $; and since $a \dnfo_K b$,
 $tp(a/K)$  is internal to $ Fix\si  \cap{\mathcal C}    $.\\
$\Box$

\section{Arc Spaces in Difference-Differential Fields}

In \cite{arcs} Moosa, Pillay and Scanlon prove a dichotomy theorem for fields with finitely many commuting derivations.
We adapt their proof to our case.\\

Let $K$ be a field, and $K^{(m)}$ the $K$-algebra
$K[\epsilon]/(\epsilon^{m+1})$. Then, identifying $K^{(m)}$ with
$K\cdot 1\oplus K\cdot \epsilon \ldots \oplus K\cdot \epsilon^m$, we see
that the $K$-algebra $K^{(m)}$ is quantifier-free interpretable in $K$,
if one encodes elements of $K^{(m)}$ by $(m+1)$-tuples of $K$.\\

Let $V\subset \mathbb{A}^\ell$ be a variety defined over $K$.
For $m\in \na$, we
consider
the set $V(K^{(m)})$ of $K^{(m)}$-rational points of $V$.

Using the quantifier-free interpretation of $K^{(m)}$ in $K$, we may
(and will) identify $V(K^{(m)})$ with a subvariety $\mathcal{A}_mV(K)$ of
$\mathbb{A}^{(m+1)\ell}(K)$.  The variety $\mathcal{A}_mV$ is called the $m$-th arc
 bundle of
$V$. More precisely, if $f_1,\ldots,f_k\in
K[X_1,\ldots,X_\ell]$ generate the ideal $I(V)$, then the ideal of
$\mathcal{A}_mV$ is generated by the polynomials $f_{j,t}\in K[X_{i,t}|1\leq
i\leq\ell, 0\leq t\leq m]$, $1\leq j\leq k, 0\leq
t\leq m$, which are defined by the identity
$$f_j((\sum_{t=0}^m x_{i,t}\epsilon^t)_{1\leq i\leq l})=\sum_{t=0}^{m}
f_{j,t}(x_{{i,t}_{1\leq i\leq\ell, 0\leq t\leq m}})\epsilon^t.$$

If $r>m$,  the natural map
$K^{(r)}\to K^{(m)}$ then induces a map $V(K^{(r)})\to V(K^{(m)})$, which
in turn induces a morphism $\rho_{r,m}: \mathcal{A}_rV\to \mathcal{A}_mV$. 

Moreover, given a morphism of varieties $f:U\to V$  defined over $K$, the
natural morphism $U(K^{(m)})\to V(K^{(m)})$ induced by $f$ gives
rise to a morphism $\mathcal{A}_mf:\mathcal{A}_mU\to \mathcal{A}_mV$.

Let us  write $\rho_m$ for $\rho_{m,0}$. For $a\in V(K)$ the $m$-th
arc space of $V$ at $a$, $\mathcal{A}_mV_a$ is the fiber of $\rho_m$ over $a$.
The following three results appear in \cite{arcs}.

%First we introduce some notation.\\
%Let $K$ be an algebraically closed field, let
%$K^{(m)}=K[\epsilon]/(\epsilon)^{m+1}$. We view $K^{(m)}$ as a $K$-algebra
%under the map $a \mapsto A+0\epsilon+ \cdots +0\epsilon $.

%\begin{defi}\label{arc00}
%Let $V$ be a variety in $\mathbb{A}^l$ defined ovre $K$. Let $(f_1, \cdots, f_k)$ be a tuple of polynomials over $K[X_1, \cdots, X_l]$ generating $V$.
%The $m$-th arc space of $V$, $\mathcal{A}_mV$ is the variety of
%$\mathbb{A}^{l(m+1)}$ defined by the polynomials $(f_{(j,t)})$ of
%$K[(X_{i,s})_{1 \leq i \leq l,0 \leq s \leq m}]$; where $(f_{i,j})$ are determined
%by
%$$f_j((\sum_{t=0}^mx_{i,t}\epsilon^t)_{1 \leq i \leq l })=\sum_{t=0}^mf_{j,t}\epsilon^t $$
%for $1 \leq j \leq k$, in the ring
%$K[(X_{i,s})_{1 \leq i \leq l,0 \leq s \leq m},\epsilon]/(\epsilon^m)$.
%\end{defi}

%As in \cite{arcs} we identify $\mathcal{A}_mV(K)$ with $V(K^{(m)})$.

%For $r>m$, the quotient map $K^{(r)} \to K^{(m)}$ induces a map $\rho_{r,m}:\mathcal{A}_rV \to \mathcal{A}_mV$.
%We set $\rho_r=\rho_{r,0}: \mathcal{A}_rV \to V$.

%Given a morphism of varieties $f: U \to V$, we define
%$\mathcal{A}_mf:\mathcal{A}_mU \to \mathcal{A}_mV$ as follows: assume
%that $U \subset \mathbb{A}^l$ and $V \subset \mathbb{A}^r$.
%Let $f=(f_1,\cdots,f_r)$. If we see $b \in \mathcal{A}_mU(K)$ as an
%element of $\mathbb{A}^l(K^{(m)})$, we have
%$\mathcal{A}_mf(b)=(f_1(b), \cdots, f_r(b))$.

\begin{lem}\label{arc01}
Let $U,V$ be two algebraic varieties, and let $f:U \to V$ be a morphism, all defined over $K$.
Let $m \in \na$ and $a \in \mathcal{A}_mU(K)$ be such that for all $m$, $\bar{a}=\rho_m(a)$
and $\bar{f(a)}=\rho_m(f(a))$ are non-singular. Let $U'$ be the fiber of
$\rho_{m+1,m}:\mathcal{A}_{m+1}U \to \mathcal{A}_mU$ over $a$ and $V'$ the fiber of
$\rho_{m+1,m}:\mathcal{A}_{m+1}V \to \mathcal{A}_mV$ over $\mathcal{A}_mf(a)$.
Let $\bar{a}=\rho_m(a)$. Then there are biregular maps $\varphi_U:U' \to T(U)_{\bar{a}}$  and
$\varphi_V:V' \to T(V)_{f(\bar{a})}$  such that the following diagram is commutative:

$$\xymatrix{
% U \ar@/_/[ddr]_y \ar@/^/[drr]^x
 %  \ar@{.>}[dr]|-{(x,y)}            \\
  & U' \ar[d]^{\varphi_U} \ar[r]^{\mathcal{A}_m(f)}
                 & V' \ar[d]^{\varphi_V}       \\
  & T(U)_{\bar a} \ar[r]^{df_{\bar a}}   & T(V)_{f({\bar a})}                }$$

\end{lem}

\begin{lem}\label{arc02}
Let $U,V$ be algebraic varieties defined over $K$, and let $f:U \to V$ be a dominant map defined over $K$.
Let $a \in U(K)$ be non-singular such that $f(a)$ is non-singular and the rank of $df_a$  equals $dimV$. Then for every $m \in\na$
the map $\mathcal{A}_m(f):\mathcal{A}_mU_a(K) \to \mathcal{A}_mV_{f(a)}(K)$ is surjective.
\end{lem}

\begin{lem}\label{arc03}
Let $U,V,W$ be algebraic varieties defined over $K$ such that $U,V \subset W$. Let $a\in U(K) \cap V(K)$ be non-singular.
 Then
$U=V$ if and only if $\mathcal{A}_mU_a(K)=\mathcal{A}_mV_a(K)$ for all $m \in \na$.
\end{lem}
Let $\nabla_m:V \to \tau_m(V)$ be defined by $x \mapsto (x,Dx,\cdots,D^mx)$ and
let $\pi_{l,m}:\tau_l(V) \to \tau_m(V)$ be the natural projection for $l \geq m$.
$S_m(V)$ will
denote the Zariski closure of $\{(x, \cdots, \si^m(x)): x \in V\}$.
Let $q_m:V \to S_m(V)$ be defined by $x \mapsto (x, \cdots, \si^m(x))$ and let
$p_{l,m}:S_l(V) \to S_m(V)$ be the natural projections for $l \geq m$.\\

We now define a notion of difference-differential prolongation.

Let $\Phi_m(V)=\tau_m(S_m(V))$, let $\psi_m:V \to \Phi_m(V)$ be such that
$x \mapsto \nabla_m(q_m(x))$ and for $l \geq m$ let
$t_{l,m}:\Phi_l(V) \to \Phi_m(V)$ be defined by $t_{l,m}=\pi_{l,m} \circ p_{l,m}$.
Let us denote $\pi_l=\pi_{l,0}$, $p_l=p_{l,0}$, $t_l=t_{l,0}$,
$\Phi(V)=\Phi^1(V)=\Phi_1(V)$ and $\Phi^{m+1}(V)=\Phi(\Phi^m(V))$.
We define $\psi=\psi^1=\psi_1:V \to \Phi(V)$ and
$\psi^{m+1}(V)=\psi(\psi^m):V \to \Phi^{m+1}(V)$.

Let $(\mathcal{U},\si,D)$ be a saturated model of {\it DCFA}, let $K$ be a
difference-differential subfield of ${\mathcal U}$.
We can identify
$\tau_m(\mathcal{A}_rV)(K)$ with $\mathcal{A}_r\tau_m(V)(K)$.

We extend $\si$ and $D$ to $K^{(m)}$ by defining $\si(\epsilon)=\epsilon$
and $D \epsilon=0$. Then
we can identify $\mathcal{A}_r(S_m(V))(K)$ with $S_m(\mathcal{A}_r(V))(K)$.
We can, then, identify $\mathcal{A}_r(\Phi_m(V))(K)$ with
$\Phi_m(\mathcal{A}_m(V))(K)$.\\

%A model $(L,\si,D)$ of {\it DCFA} can be characterized as an algebraically closed difference-differential field such that, given a projective system of
% dominant maps of varieties $\{\mu_{l,m}:V_l \to V_m \}$ where $V_l$ is
 %a closed subvariety of $\Phi^{l-m}(V_m)$ and $\mu_{l,m}$ is the restriction
 %of $t_{l-m}$ to $V_l$ ; there is $a \in V_0(K)$ such that
 %$\psi^l(a) \in V_l(K)$ for all $l$ (This is a direct consequence of the axioms of
 %{\it DCFA}).

Let $V$ be a $(\si,D)$-variety given as a $(\si,D)$-closed subset of an
algebraic variety $\bar{V}$. We define $\Phi_m(V)$ as the Zariski closure
of $\psi_m(V)$ in $\Phi_m(\bar{V})$. Thus $V$ is determined by the prolongation sequence
$\{ t_{l,m}:\Phi_l(V) \to \Phi_m(V)):l \geq m \}$, since
$V(\mathcal{U})=\{a \in {\bar V}({\mathcal U}): \psi_l(a) \in \Phi_l(V) \forall l \}$.
We call this sequence the prolongation sequence of $V$. \\

\begin{prop}\label{arc04}
Let $\{ V_l \subset \Phi_l(\bar{V}):l \geq 0\}$ be a sequence of algebraic varieties
and $\{t_{m,l}:V_m \to V_l, m\geq l\}$ a sequence of morphisms such that:
\begin{enumerate}
\item $t_{l+1,l} \upharpoonright V_{l+1} \to V_l$ is dominant.
\item After embedding $\Phi_l(\bar{V})$ in $\Phi^l(\bar{V})$ and
$\Phi_{l+1}(\bar{V})$ in $\Phi^{l+1}(\bar{V})$,
\begin{enumerate}
\item $V_{l+1}$ is a %closed 
subvariety of $\Phi(V_l)$.
\item Let $\pi'_1:\Phi(V_l)\to \tau(V_l)$ and $\pi'_2:\Phi(V_l)\to
\tau(V_l^\sigma)$ be the projections induced by $\Phi(V_l)\subset
\tau(V_l)\times \tau(V_l^\sigma)$; then $\pi'_1(V_{l+1})^\sigma$ and
$\pi'_2(V_{l+1})$ have the same Zariski closure.
%$V_l$ is a closed subvariety of $\Phi(V_l)=
%\tau(S(V_l)) \subset \tau(V_l \times V_l^{\si}) \simeq \tau(V_l) \times
%\tau(V_l^{\si})$, the projections $\pi'_1,\pi'_2$ from $V_{l+1}$
% to $\tau(V_l)$ and $\tau(V_l^{\si})$ are such that $\si$ applied to the Zariski closure of
% $(\pi'_1(V_{l+1}))$ equals  the Zariski closure of$\pi'_2(V_{l+1})$; and the projections from
% the Zariski closure of $\pi(V_{l+1})$ to the Zariski closure of $V_l$ and onto he Zariski closure of $V_l^{\si}$ are dominant.
\end{enumerate}
\end{enumerate}
Then there is a (unique) $(\si,D)$-variety $V$ with prolongation sequence $\{t_{m,l}:V_m \to V_l, m\geq l\}$ .
\end{prop}
{\it Proof}:\\

We work now in a saturated model $\mathcal U$ of {\it DCFA}.
For each $l$, as the maps $\pi_{m,j}$ are dominant, the system $\{p_{m,l}(V_m), \pi_{m,j}:m >j \geq l\}$
defines a differential subvariety  $W_l$ of $\bar{V} \times \cdots \times \bar{V}^{\si^l}$.

Condition (1) implies that for $m$ sufficiently large, an $(m,D)$-generic
of $p_{m,l+1}(V)$ is sent by $p_{l+1,l}$ to an $(m,D)$-generic of 
$p_{m,l}(V)$. Hence,
a $D$-generic of $W_{l+1}$ is sent by $p_{l+1,l}$ to a $D$-generic of 
$W_l$.

By conditions (2) (b) and (1), the  map $t'_{l+1,l}:V_{l+1}\to V_l^\sigma$ 
induced
by $\Phi(V_l)\to V_l^\sigma$ is dominant. Hence, considering the 
natural
projection
 $p'_{l+1,l}:S_{l+1}(\bar
V)\to S_l(\bar V)^\sigma$, and reasoning as above, we obtain that
$p'_{l+1,l}$ sends a $D$-generic of $W_{l+1}$ to a $D$-generic of
$W_l^\sigma$.

%Since $V_{l+1} \subset \tau(V_l) \times \tau(V_l)^{\si}$, $\si$ applied to
%the Zariski closure of projection of $V_{l+1}$ to $\tau(V_l)$ equals the
%Zariski closure of the projection of $V_{l+1}$ to $\tau(V_l)^{\si}$ and $t_{l+1,l}:V_{l+1,l} \to V_l$ is
%dominant, this gives us a dominant map $t'_{l+1,l}:V_l \to V_l^{\si}$. Hence the restrictions of
%$t_{l+1,l}$ and $t'_{l+1,l}$ to $W_{l+1}$ are dominant onto $W_l$ and $W_l^{\si}$ respectively.
Hence by the axioms of {\it DCFA}, for every $l$ there is $a$ such that $\psi_l(a)$ is a generic of $V_l$ over $K$. 
By saturation, there is $a$ such that for all $l$ $\psi_l(a)$ is a generic of $V_l$. Then
$\{t_{m,l}:V_m \to V_l, m\geq l\}$  is the prolongation sequence of the $(\si,D)$-locus of $a$ over $K$.\\
$\Box$
%Then there is a unique $(\si,D)$-subvariety $V$ of $\bar{V}$ such
%that $\Phi_l(V)=V_l$. In particular, there is a tuple $a$ in $\mathcal{U}$ such that, for all
%$l \in \na$, $\psi_l(a)$ is a generic of $V_l$ over $K$.

%(since 2. holds for every $l$, it guarantees that
%an $(1,D)$-generic of $V_{l+1}$ projects onto an $(1,D)$-generic
%of $\pi'_1(V_{l+1})$ and a $(1,D)$-generic of $\pi'_2(V_{l+1})$).

\begin{defi}
Let $V$ be a $(\si,D)$-subvariety of the algebraic variety $\bar{V}$.
We say that a point $a \in V$ is non-singular if, for all $l$, $\psi_l(a)$ is
a non-singular point of $\Phi_l(V)$, the maps $dt_{l+1,l}$ and $dt'_{l+1,l}$ at $\psi_{l+1}(a)$ have rank
equal to $dim V_l$ and the maps $d \pi_1'$ and $d\pi_2'$ (as defined above) at $\psi_{l+1}(a)$ have rank equal to 
the dimension of the Zariski closure of $\pi_1'(\Phi_{l+1}(V))$.
\end{defi}

\begin{prop}\label{arc1}
Let $(K,\si,D)$ be a model of {\it DCFA}. Let $V$ be a $(\si,D)$-variety
given as a closed subvariety of an algebraic variety $\bar{V}$.
Let $m \in \na$ and $a\in V(K)$ a non-singular point.
Then $\{\mathcal{A}_m(t_{r,s}): \mathcal{A}_m\Phi_r(V)_{\psi_r(a)} \to \mathcal{A}_m\Phi_s(V)_{\psi_s(a)}, r \geq s \}$
 form the
 $(\si,D)$-prolongation sequence of
a $(\si,D)$-subvariety of $\mathcal{A}_m\bar{V}_a$.
We define the $m$-th arc space of $V$ at $a$, $\mathcal{A}_mV_a$, to be this subvariety.
%Moreover, if $a \in V(K)$  is a non-singular point of
%$V$ 
%and $d(t_{r,s})_{\psi_r(a)}$ has full rank for
%$r \geq s$,
We have also that $\Phi_r(\mathcal{A}_mV_a)=\mathcal{A}_m\Phi_r(V)_{\psi_r(a)}$ for all $r$.
\end{prop}

{\it Proof}:\\

Since we can identify $\mathcal{A}_m\Phi_r(\bar{V})$ with
$\Phi_r(\mathcal{A}_m{\bar V})$, we look at $\mathcal{A}_m\Phi_r(V)_{\psi_r(a)}$ as an
algebraic subvariety of $\Phi_r(\mathcal{A}_m \bar{V})_{\psi_r(a)}$.
We have that $\Phi_{r+1}(V) \subset \Phi(\Phi_r(V))$ for all $r$. Since
$\mathcal{A}$ preserves inclusion we have
$\mathcal{A}_m\Phi_{r+1}(V)_{\psi_{r+1}(a)} \subset \mathcal{A}_m\Phi(\Phi_r(V)_{\psi_{r}(a)})=\Phi(\mathcal
{A}_m\Phi_r(V)_{\psi_{r}(a)})$. This shows conditions $1$ and $2(a)$.

Moreover, the maps $t_{r,s}:\Phi_r(V) \to \Phi_s(V)$
are dominant, and  by \ref{arc02}, the maps
$\mathcal{A}(t_{r,s}):\mathcal{A}_m\Phi_r(V)_{\psi_r(a)} \to \mathcal{A}_m\Phi_s(V)_{\psi_s(a)}$,
are dominant.
Applying $\mathcal{A}_m$ to the dominant maps $\pi_1':\Phi_{r+1}(V) \to \tau(\Phi_r(V))$ 
and $\pi_2':\Phi_{r+1}(V) \to \tau(\Phi_r(V))^{\si}$,
using the hypothesis on $a$ and \ref{arc02}, we get
$$\mathcal{A}_m\pi_1'(\mathcal{A}_m(\Phi_{r+1}(V)_{\psi_{r+1}(a)}))=
\mathcal{A}_m(\pi_1'(\Phi_{r+1}(V))_{\pi_1'(\psi_{r+1}(a))})$$
and
$$\mathcal{A}_m\pi_2'(\mathcal{A}_m(\Phi_{r+1}(V)_{\psi_{r+1}(a)}))=
\mathcal{A}_m(\pi_2'(\Phi_{r+1}(V))_{\pi_2'(\psi_{r+1}(a))})$$
and since $\pi_1'(\Phi_{r+1}(V))^{\si}$ and $\pi_2'(\Phi_{r+1}(V))$ have the same Zariski closure,
and $\si(\pi_1'\psi_{r+1}(a))=\pi_2'\psi_{r+1}(a)$ we get condition $2(b)$.

%and by the identification  $\mathcal{A}_m\Phi_r(\bar{V})_{\psi_r(a)}$ with
%$\Phi_r(\mathcal{A}_m{\bar V})_{\psi_r(a)}$ and the definition of non-singular point
% the sequence satisfies the conditions above.
Hence $\{
\mathcal{A}_m(t_{r,s}):A_m\Phi_r(V)_{\psi_r(a)} \to  A_m\Phi_s(V)_{\psi_s(a)}, r \geq s  \}$ is the $(\si,D)$-prolongation
sequence of a $(\si,D)$-subvariety $W$ of $\mathcal{A}_m\bar{V}_a$, where
$W(K)=\{x \in \mathcal{A}_m\bar{V}_a(K): \psi_r(x) \in \mathcal{A}_m\Phi_r(V)_{\psi_r(a)}(K), 
r \geq 0\}$ and $\mathcal{A}_m\Phi_r(V)_{\psi(a)}=\Phi_r(W)$ for all $r$. We define then
$\mathcal{A}_mV_a=W$.\\
%The map $\rho_m:\mathcal{A}_m\bar{V} \to \bar{V}$ restricts to a surjective map
%$\mathcal{A}_mV \to V$ that we shall denote also $\rho_m$. Let $a \in V(K)$ be
%as in the hypothesis of the second part of the theorem.
%We define $\mathcal{A}_mV_a$ to be the fiber of $\rho_m$ over $a$.
%Then $\mathcal{A}_mV_a(K)=\{b \in \mathcal{A}_m\bar{V}_a(K):
% \psi_r(b) \in \mathcal{A}_m\Phi_r(V)_{\psi_r(a)}(K), r \geq 0 \}.$ So for the
% "moreover" part it is enough to prove that, viewed as a sequence of
% algebraic subvarieties of $\Phi_r(\mathcal{A}_m\bar{V})_a$ for $r \geq 0$,
% $\{\mathcal{A}_m(t_{r,s}):\mathcal{A}_m\Phi_r(V)_{\psi_r(a)} \to
% \mathcal{A}_m\Phi_s(V)_{\psi_s(a)}, r \geq  s\}$
%is the $(\si,D)$-prolongation sequence of a $(\si,D)$-subvariety of
%$\mathcal{A}_m\bar{V}_a.$ But this is clear: By \ref{arc02} and our
% choice of $a$
%$\mathcal{A}_m(t_{r,s}):\mathcal{A}_m\Phi(V)_{\psi_r(a)} \to \mathcal{A}_m\Phi_s(V)_{\psi_s(a)} $ is surjective 
%and we know that they are dominant.\\
$\Box$

\begin{lem}\label{arc2}
Let $U,V$ be two $(\si,D)$-subvarieties of an algebraic variety $\bar{V}$.
Let $a \in U(K) \cap V(K)$ be a non-singular point of $U$. Then
$U=V$ if and only if $\mathcal{A}_m\Phi_l(U)_{\psi_l(a)}=\mathcal{A}_m\Phi_l(V)_{\psi_l(a)} $ for all $m,l$.
\end{lem}

{\it Proof}:\\

If $\mathcal{A}_mU_a(K)=\mathcal{A}_mV_a(K)$ for all $m$, then
$\Phi_r(\mathcal{A}_mU_a)(K)=\Phi_r(\mathcal{A}_mV_a)(K)$. Thus, by \ref{arc1},
$\mathcal{A}_m\Phi_r(U)_{\psi_r(a)}(K)=\mathcal{A}_m\Phi_r(V)_{\psi_r(a)}(K)$.
Hence, for all $r$ and $m$, we have
$\mathcal{A}_m\Phi_r(U)_{\psi_r(a)}=\mathcal{A}_m\Phi_r(V)_{\psi_r(a)}$.
Lemma \ref{arc03} implies that $U$ and $V$ have the same $(\si,D)$-prolongation sequence. Hence $U=V$.\\
$\Box$

\begin{defi}
Let $V$ be a variety and $a$ a non-singular point of $V$. We define the $(\si,D)$-tangent space
 $T_{\si,D}(V)_a$ of $V$ at $a$
as follows:

Let $P_r$ be a finite tuple of polynomials generating $I(\Phi_r(V)_{\psi_r(a)})$. Then
$T_{\si,D}(V)_a$ is defined by the equations $J_{P_r}(\psi_r(a)) \cdot (\psi_r(Y))= 0$.
In other words, the prolongation sequence of $T_{\si,D}(V)_a$ is 
$dt_{l,r}: T(\Phi_l(V))_{\psi_l(a)} \to T(\Phi_r(V))_{\psi_r(a)},l \geq r\}$,
where $T$ denotes the usual tangent bundle and $t_{l,r}$ the natural projection 
$\Phi_l(V)_{\psi_l(a)} \to \Phi_r(V)_{\psi_r(a)}$.
\end{defi}

\begin{rem}
Let $a$ be a non-singular point of the $(\si,D)$-variety $V$. Then $T_{\si,D}(V)_a$ is a subgroup
of $\gr_a^n(K)$, and by the same arguments as above, its prolongation sequence is
$(d(t_{l,r})_{\psi_l(a)}:T(\Phi_l(V))_{\psi_l(a)} \to T(\Phi_r(V))_{\psi_r(a)})_{l \geq r}$.
\end{rem}

%Note that $T_{\si,D}(V)_0$ is a nontrivial additive group: Indeed,
%$T_{\si,D}(V)_0=\{y:\psi_r(y)\in T(\Phi_r(V))_0 \forall r \in \na\}$.
% Let $T^{(0)}(V)_0=T(V)_0$
%and $T^{(r+1)}(V)_0=\{y:
%y \in T^{(r)}(V)_0 \wedge \psi_{r+1}(y) \in T(\phi_{r+1}(V))_0  \}$. 
%Then $T_{\si,D}(V)_0= \cap_{r \in \na}T^{(r)}(V)_0$. Since $T^{(r)}(V)_0$ is a
%decreasing sequence of varieties, by Noetherianity, there is
%$M \in \na$ such that $T_{\si,D}(V)_0=T^{(M)}(V)_0$.
\begin{lem}\label{arc3}
Let $V$ be a $(\si,D)$-variety in $\mathbb{A}^l$ and $a$ a non-singular point of $V$. Then $\mathcal{A}_1V_a$
is isomorphic to $T(V)_a$. Let  $\bar{V}$ be the Zariski closure of 
$V(\mathcal{U})$ in $\mathbb{A}^l$ and $m \in \na$; 
then the map given by lemma \ref{arc01} which identifies the fibers of
$\mathcal{A}_{m+1}\bar{V}_a \to \mathcal{A}_m\bar{V}_a$ with $T(\bar{V})_a$ restricts
to an isomorphism of the fibers of
$\mathcal{A}_{m+1}V_a \to \mathcal{A}_mV_a$ with $T(V)_a$.
\end{lem}

{\it Proof}:\\

We identify $\mathcal{A}_1\bar{V}$ with $T(\bar{V})$.
Let $b \in T(\bar{V})_a(K)$. By definition
$(a,b) \in \mathcal{A}_1V(\mathcal{U})$ if and only if $\psi_r(a,b) \in T(\Phi_r(V))(K)$
 for all $r$. We view $T(\Phi_r(V))_{\psi_r(a)}$ as an algebraic subvariety of
 $\Phi_r(T(\bar{V}))$ under the identification of $T(\Phi_r(\bar{V}))$ with
 $\Phi_r(T(\bar{V}))$; in particular we identify $\psi_r(a,b)$ with
 $(\psi_r(a),\psi_r(b))$. Hence $b\in \mathcal{A}_1V_a(K)$
 if and only if $b \in T(V)_b$ and the first part of the theorem is proved.

 Now we look at the map given in \ref{arc01}. In particular, if
 $c \in \mathcal{A}_mV_a(K)$ and $r \geq 0$, by \ref{arc1},
 $\psi_r(c) \in \mathcal{A}_m\Phi_r(V)_{\psi_r(a)}$ and the following diagram commutes

$$\xymatrix{
% U \ar@/_/[ddr]_y \ar@/^/[drr]^x
 %  \ar@{.>}[dr]|-{(x,y)}            \\
  & (\mathcal{A}_{m+1}{\bar V}_a)_c \ar[d] \ar[r]
                 & (\mathcal{A}_{m+1}\Phi_r({\bar V})_{\psi_r(a)})_{\psi_r(c)} \ar[d]       \\
  & T({\bar V})_{a} \ar[r]   & T(\Phi_r({\bar V}))_{\psi_r(a)}                }$$

where the horizontal arrows are $\psi_r$ and the vertical arrows are the maps given by \ref{arc01} applied to $\bar{V}$
 and $\Phi_r(\bar{V})$. So $(\mathcal{A}_{m+1}V_a)_c$ is
 identified with $T_{\si,D}(V)_a$.\\
 $\Box$\\

\begin{notadefi}
In analogy with the material of \cite{kolg}, section 0.3,
 since the $(\si,D)$-topology is Noetherian, given a difference-differential
 subfield $F$ of $K$ and $a \in K$ there is a numerical polynomial 
 $P_{a/F}(X) \in {\mathbb Q}[X]$ of degree at most 2, such that for sufficiently large $r \in \na$,
 $P_{a/F}(r)=\td(\psi_r(a)/F)$. We call the degree of $P_{a/F}$ the
 $(\si,D)$-type of $a$ over $F$, and the leading coefficient of $P_{a/F}$
 the $(\si,D)$-dimension of $a$ over $F$, it is denoted
 $dim_{\si,D}(a/F)$. For a $(\si,D)$-variety $V$ defined over $F$ we define
 $P_V=P_{a/F}$ where $a$ is a $(\si,D)$-generic of $V$ over $F$.
 We have that the $(\si,D)$-type of $a$ over $F$ is 2 if and only if $a$ contains an element which
is $(\si,D)$-transcendental over $F$.

Let $(\mathcal{U},\si,D)$ be a saturated model of {\it DCFA}, let $F=acl(F) \subset \mathcal{U}$.
Let $a \in \mathcal{U}$ and let $p=tp(a/F)$. We denote by $m(p)$ (or by $m(a/F)$) the $(\si,D)$-type of $a$ over $F$
and we write $dim_{\si,D}(p)$ for $dim_{\si,D}(a/F)$.
If $p'$ is a non-forking extension of $p$ then $m(p)=m(p')$ and $dim_{\si,D}(p) =dim_{\si,D}(p')$.
If $A$ is an arbitrary subset of $\mathcal U$ we write $m(a/A)$ instead of $m(a/acl(A))$.

If $V$ is a $(\si,D)$-variety over $K$, $m(V)$ denotes the $(\si,D)$-type of $V$. Then, if
$a$ is a $(\si,D)$-generic of $V$, $m(V)=m(qftp(a/F))$.

\end{notadefi}

 \begin{cor}\label{arc4}. 
Let $V$ be a $(\si,D)$-variety in $\mathbb{A}^l$, and $m \in \na$. Then
 for $a \in V(K)$ non-singular, the
$(\si,D)$-type of $V$ and $\mathcal{A}_mV_a$ are equal.
\end{cor}

{\it Proof}:\\

By \ref{arc1} $\Phi_r(\mathcal{A}_mV_a)=\mathcal{A}_m\Phi_r(V)_{\psi_r(a)}$. But if $b$ is a non-singular
point of a variety $U$, then we have $dim(\mathcal{A}_mU_b)=m dim(U)$.\\
$\Box$

\begin{rem}\label{arc41}
By \ref{arc3}, for $m=1$ and for $a \in V(K)$ non-singular, we have $P_V=P_{T(V)_a}$.
\end{rem}

\begin{lem}\label{arc5} Let $F=acl(F)$. Then
\begin{enumerate}
\item $m(a,b/F)=max\{m(a/F),m(b/F) \}$.
\item If $m(a/F)=m(b/F)$ then $dim_{\si,D}(a,b/F)=dim_{\si,D}(a/F)+dim_{\si,D}(b/Fa)$.
\item If $m(a/F)>m(b/F)$ then $dim_{\si,D}(a,b/F)=dim_{\si,D}(a/F)$.
\end{enumerate}
\end{lem}

{\it Proof}:\\

It suffices to compute the degree and the leading coefficient of the respective polynomials.\\
$\Box$

%We will not define weight of a type, all we need to know for our case is that a stationary type has finite weight, and
%regular types have weight 1. For more
%details we refer to \cite{wag}, Chapter 5, section 5.2.

\begin{defi}\label{arc8}
Let $p$ be a regular type. We say that $p$ is
$(\si,D)$-type minimal if for any type $q$, $p \not\perp q$ implies
$m(q) \geq m(p)$.
\end{defi}

\begin{defi}
A $(\si,D)$-variety $V$ is $(\si,D)$-type minimal if for every
proper $(\si,D)$-subvariety $U$, $m(V) < m(U)$.
\end{defi}

\begin{lem}\label{arc6}
Let $p$ be a  type and let $V$ be the $(\si,D)$-locus of $p$ over $K$
 (i.e. the Kolchin closure of the set of $a$ realizations of $p$) 
If $V$ is $(\si,D)$-type minimal then $p$ is regular and $(\si,D)$-type minimal.
\end{lem}

{\it Proof}:\\

Let $a$ be a realization of a forking extension of $p$ to some
$L=acl(L) \supset K$. Let $b$ realize a nonforking extension of $p$
to $L$. Let $U$ be the $(\si,D)$-locus of $(a,b)$ over $L$. Then the projection
on the second coordinate: $U \to V$ is dominant, thus $m(a,b/L) \geq m(V)$.
Now if $a \dfo_L b$, then the $(\si,D)$-locus of $b$ over
$acl(La)$  is a proper subvariety of $V$ and therefore  $m(b/La)<m(V)$; from $m(a/L)<m(V)$, we deduce
$m(a,b/L)<m(V)$ which is impossible.\\
$\Box$

\begin{lem}\label{arc7}
If $p$ is a type over $K$, there is a finite sequence of regular types $p_1, \cdots, p_k$ such that $m(p) \geq m(p_i)$ for all $i$ and $p$
is domination-equivalent to $p_1 \times \cdots \times p_k$.
\end{lem}

{\it Proof}:\\

By \ref{prt13} it suffices to show that given a regular type $q$, such that $p \not\perp q$, there
is a regular type $r$ such that $q \not\perp r$ and $m(r) \leq m(p)$.
Let $a$ be a realization of a nonforking extension of $p$ to some $L$ and let
$b$ be a realization of a nonforking extension of $q$ to $L$ such that
$a \dfo_L b$. Let $c=Cb(tp(a/L,b))$. Thus $c \not\in acl(L)$ and
$c \in acl(Lb)$. So $r=tp(c/L)$ is regular (because $c \in acl(Lb)$) and non-orthogonal to $q$.
On the other hand, there are $a_1, \cdots,a_l$ realizations of $p$ such that
$c \in dcl(La_1 \cdots a_l)$. Then, by \ref{arc5}, $m(r) \leq m(q)$.\\
$\Box$

%\begin{cor}\label{arc8}
%Let $p$ be a $(\si,D)$-type minimal regular type.
%Then for any type $q$, $p \not\perp q$ implies
%$m(p) \leq m(q)$.
%\end{cor}
%
%{\it Proof}:\\
%
%Suppose $p \not\perp q$ and $m(p) > m(q)$. Let $q_1,\cdots,q_k$ be the sequence given by \ref{arc7}
%applied to $q$. Then there is $i$ such that $p \not\perp q_i$ which contradicts the $(\si,D)$-type minimality of $p$ .\\
%$\Box$

\begin{lem}\label{arc9}
Let $G$ be a $(\si,D)$-vector group (that is, a $(\si,D)$-variety which is a subgroup of
$\gr_a^k$ for some $k$). Then $T_{\si,D}(G)_0$ is definably isomorphic to $G$.
Moreover, if $H$ is a $(\si,D)$-subgroup of $G$, then the
restriction of this isomorphism to $H$ is an isomorphism between $H$ and $T_{\si,D}(H)_0$.
\end{lem}

{\it Proof}:\\

Suppose that $G$ is a $(\si,D)$-subgroup of $\gr_a^k$. For each $r \in \na$, $\Phi_r(G)$ is a subgroup
of $\Phi_r(\gr_a^k)=\gr_a^{k(r+1)^2}$. Let $z_r:\Phi_r(\gr_a^k) \to T(\Phi_r(\gr_a^k))$ defined by
$x \mapsto(0,x)$; this map identifies $\Phi_r(\gr_a^k)$ and $T(\Phi_r(\gr_a^k))_0$.
Since $\Phi_r(G)$ is an algebraic subgroup of $\gr_a^{k(r+1)^2}$ , its defining ideal is generated by
linear  polynomials, and thus its tangent space at 0 is defined by the same polynomials.
This means that $z_r$ restricts to an isomorphism $\Phi_r(G) \to T(\Phi_r(G))_0$. Hence $(z_r:r \geq 0)$
identifies the prolongation sequence of $G$ and the prolongation sequence of $T_{\si,D}(G)_0$.
For the moreover part, it suffices to note that, by our construction above, the restriction
of $z_r$ to $\Phi_r(H)$ is an isomorphism between  $\Phi_r(H)$ and $T(\Phi_r(H))_0$. \\
$\Box$

We will see now to reduce some questions concerning groups definable in
a model of {\it DCFA} to questions on groups definable in
{\it DCF} or {\it ACFA}. These ideas are, actually, implicit in the axioms of {\it DCFA}.

Let $ G $ be a connected differential algebraic group defined over $E =acl(E) $.

For each $n \in {\mathbb N}  $ let $G^{(n)}=G \times \sigma(G) \times \cdots \times \sigma^n(G)  $,
and let $q_n$ be the group homomorphism from $G$ to $G^{(n)}$ defined by $q_n(g)=(g,\sigma(g), \cdots , \sigma^n(g))  $.

Let $g$ be a generic point of $G$ such that the tuples $ g,\sigma(g), \cdots , \sigma^n(g)$ are differentially  independent over
$E$; then $q_n(g)$ is a generic point of $G^{(n)}$; thus  $q_n(G)$ is dense in $ G^{(n)}$ (for the $D$-topology)
and $G^{(n)}$ is connected (in {\it DCF}).

Let $H$ be a definable subgroup of $G$. For each $n \in {\mathbb N}$ let $H^{(n)}$ be the differential Zariski closure of
$q_n(H)  $ in $G^{(n)}  $; then $H^{(n)}  $ is a differential algebraic  subgroup of $ G^{(n)} $.

Let $\tilde{H}^{(n)}= \{g \in G: q_n(g) \in  H^{(n)}    \}  $. These subgroups of $G$ form a decreasing sequence
of quantifier-free definable groups containing $H$. Let $ \tilde{H} =\bigcap_{n \in{\mathbb N}} \tilde{H}^{(n)}  $; 
since the
 ($\sigma,D $)-topology is Noetherian, there is $N \in {\mathbb N}  $ such that $ \tilde{H}  =\tilde{H}^{(N)} $.
Then $\tilde H$ is the
difference-differential Zariski closure of $H$. 

\begin{lem}\label{sidcl}
Let $G$ be a connected differential algebraic group and let $H$
be a definable subgroup of $G(\mathcal{U})$ defined over $E=acl(E)$, $\tilde H$
its difference-differential Zariski closure. Then $[\tilde H:H]<\infty$. 
\end{lem}

{\it Proof}:\\

Let $g, h \in G $. By definition, $g \dnfo_E h$ if and only if for every
$n \in \na$ $q_n(g)$ and $q_n(h)$ are independent over $E$ in the sense of
{\it DCF}. This implies easily that if $g \in H$, then $g$ is a generic
of $H$ if and only if for every $n \in \na$ $q_n(g)$ is a generic
of $H^{(n)}$ (in the sense of {\it DCF}). Thus a generic of $H$ will
be a generic of $\tilde{H}$ and, by \ref{GIII5}, $\s(H)=\s(\tilde{H})$
and $[\tilde{H}:H]< \infty$.\\
$\Box$

\begin{defi}\label{concomp}
Let $G$ be a quantifier-free definable group defined in a model of {\it DCFA}.
We say that $G$ is quantifier-free connected if it has no proper
quantifier-free definable subgroups of finite index. By Noetherianity, every quantifier-free definable group
$G$ has a smallest quantifier-free definable subgroup of finite index which we call the
quantifier-free connected component of $G$.
\end{defi}

\begin{rem}
In {\it DCFA} $H$ is quantifier-free-connected if and only if for all $n$ $q_n(H)$ is connected for the
$D$-topology.
\end{rem}

\begin{cor}\label{arccor1}\hspace{10cm}
\begin{enumerate}
\item Let $H$ be a quantifier-free definable subgroup of $\gr_a^n({\mathcal U})$. Then $H$
is a $Fix\si \cap {\mathcal C}$-vector space, so it is divisible and has therefore no subgroup of
finite index. This implies that every definable subgroup of $\gr_a^n(\mathcal{U})$ is
quantifier-free definable. 
\item Let $G$ be a definable subgroup of $\gr_a^n$, $H < G$. Then $G/H$ is definably isomorphic
to a subgroup of $\gr_a^l$ for some $l$.
\end{enumerate}
\end{cor}

{\it Proof}:\\

(1) Using the fact that every algebraic subgroup of a vector group is
defined by linear equations, it follows easily that every
differential subgroup of a vector group is defined by linear
differential equations. Hence, in the notation introduced above, each
$\tilde H^n$ is defined by linear differential equations, and this
implies that
$H$ is defined by linear $(\si,D)$-equations. Thus $H$ is stable by
multiplication by elements of $Fix \si\cap {\mathcal C}$, and is therefore
a $(Fix \si\cap {\mathcal C})$-vector space.

This proves the first assertion, and the others are clear, using the
fact that every definable group has finite index in its
$(\si,D)$-closure (by \ref{sidcl}).\\
(2) Let $L$ be an $l$-tuple of linear difference-differential equations such that $H=Ker(L)$.
Then $L$ defines a group homomorphism $G \to \gr_a^l$ with kernel $H$.
$L(G)$ is a definable subgroup of $\gr_a^l$.\\
$\Box$

\begin{cor}\label{arc10}
Let $G$ be a $(\si,D)$-subgroup of $\gr_a^k$. Suppose that for every proper definable subgroup $H$ of $G$,
$m(H)<m(G)$. Then $m(V)<m(G)$ for any proper $(\si,D)$-subvariety of $G$. In particular the generic type
of $G$ is regular.
\end{cor}

{\it Proof}:\\

Let $V$ be a  $(\si,D)$-type minimal  $(\si,D)$-subvariety of $G$ such that $m(V)=m(G)$. After possibly replacing
$V$ by a translate we may assume that $0 \in V$ is non-singular. By \ref{arc41}, $m(T(V)_0)=m(V)=m(G)$.
Since  $T(V)_0$ is a subgroup of $T(G)_0 \simeq G$,we obtain $T(V)_0=T(G)_0$.
By \ref{arc4}, $P_V=P_{T(V)_0}=P_{T(G)_0}=P_G$. Hence $V=G$. By \ref{arc6}, the generic type of $G$
is regular.\\
$\Box$

\begin{lem}\label{arc11}
Let $a,c$ be tuples of $\mathcal{U}$. Let $V$ be the $(\si,D)$-locus of $a$ over $K$. Assume that
$c=Cb(qftp(a/acl(Kc)))$. Then there is $m \in \na$ and a tuple $d$ in $\mathcal{A}_mV_a$ such that
$c \in K(a,d)_{\si,D}$.
\end{lem}

{\it Proof}:\\

Let $U$ be the $(\si,D)$-locus of $a$ over $acl(Kc)$.
 As {\it DCFA} eliminates imaginaries every definable set has a canonical parameter. 
Then $c$ is interdefinable with the canonical parameter of $U$ which,
by \ref{arc2}, is interdefinable over $K(a)_{\si,D}$ with the sequence of the canonical parameters of $\mathcal{A}_mU_a$ over
$K(a)_{\si,D}$. By quantifier-free stability $\mathcal{A}_mU_a$ is defined with parameters from 
$\mathcal{A}_mU_a \subset \mathcal{A}_mV_a$.\\
$\Box$

\begin{lem}\label{arc12}
Let $(K,\si,D)$ be a submodel of $(\mathcal{U},\si,D)$. Let $V$ be a $(\si,D)$-variety defined over $K$ and
let $a \in V(\mathcal{U})$ be a non-singular point. Let $b \in \mathcal{A}_mV_a$. Then there are $b_1,b_2, \cdots,b_m=b$,
such that
$b_i \in acl(Ka,b)$ and each $b_i$ is in some $K\cup \{a,b_{i-1}\}$-definable principal homogeneous space for
$T(V)_a$.
\end{lem}

{\it Proof}:\\

By \ref{arc01} and \ref{arc3} each fiber $\rho_{i+1,i}:\mathcal{A}_{i+1}V_a \to \mathcal{A}_iV_a$ is
a principal homogeneous space for $T(V)_a$. Then set $b_i=\rho_{m,i}(b)$.\\
$\Box$

\begin{lem}\label{arc13}
Let $(K,\si,D)$ be a submodel of $(\mathcal{U},\si,D)$. Let $p$ be a $(\si,D)$-type minimal
regular type over $K$ such that $m(p)=d$. If $p$ is not locally modular, then there are a vector group $G$ and a quantifier-free type
$q$  such that:
\begin{enumerate}
\item $m(q)=m(G)=d$.
\item $(x \in G) \in q$.
\item $p \not\perp q$.
\end{enumerate}
\end{lem}

{\it Proof}:\\

By \ref{prt011} and \ref{prt012} we may assume that $\s(p)=\omega^i$ where $i\in\{0,1,2\}$.
By \ref{prt10}, enlarging $K$ if necessary, there are tuples 
$a$ and $c$, with $a$ a tuple of realisations of $p$, $tp(c/K)$ $p$-internal, 
$c=Cb(a/acl(Kc))$, $tp(a/Kc)$ $p$-semi-regular and $c\notin cl_p(Ka)$.
Let $V$ be the locus of $a$ over $K$.
%Suppose that $c\not\in cl_p(Ka)$ and $\s(a/Kc)=\omega^ik$. Then $tp(a/Kc)$ is $p$-semi regular.

By \ref{arc11}, there is a $k$-tuple $d$ in $A_mV_a(\mathcal{U})$ such that $c\in
acl(K,a,d)$. For $i=1,\ldots,m$ let $d_i=\rho_{m,i}(d)$. Then for each
$i$, $d_i$ is in some $K(ad_{i-1})$-definable $T(V)_a^k$-principal
homogeneous space.

%By \ref{arc11} there is a tuple $d=(d_1,\cdots, d_k)$ of $\mathcal{A}_mV_a(\mathcal{U})$ such that  $c \in acl(K,a,d)$,
%and for each $i$, the realizations of $qftp(d_i/Kad_{i-1})$ form a coset of $T(V)_a$ by a subgroup $H$ defined
%over $K(a)$. This coset is defined over $K(a,d_{i-1})$.

Let $m=w_p(c/Ka)$. This means that for any $L=acl(L) \subset \mathcal{U}$
such that $L \dnfo_K c$, given a tuple $(g_1, \cdots, g_m)$
realizing $p^{(m)}$ we have that $g_i \dfo_L c$ for all $i$ if
and only if $g \subset cl_p(Lc)$.
As $c\in acl(K,a,d)$, $c \not\in cl_p(K,a)$, $tp(c/K)$ is $p$-internal,
there is $j \in \{1, \cdots, k\}$ such that $w_p(c/Kad_{j-1})=m$ and
$w_p(c/Kad_j) \leq m-1$.
 Let $L=acl(L) \subset \mathcal{U}$ contain $Kad_{j-1}$,
such that $L \dnfo_K C$, and $(g_1, \cdots, g_m)$ 
realizing  $p^{(m)}$ such that $g_i \dfo_L c$ for all $i$. 
Since  $w_p(c/Kad_{j-1})>w_p(c/Kad_j)$, either there  is $g_k$ such that
$g_k \dnfo_{Ld_j} c$, or $tp(g_k/Ld_j)$ forks over $L$.
In both cases,  $d_j$ and $g$ are dependent
over $L$. Hence
$tp(d_j/Kad_{j-1})\not\perp p$.

%$tp(c/La)$ is $p$-semi-simple we have that $g \dnfo_{La} d_{j-1}$
%and $g \dfo_{La} d_j$. 
%Hence $tp(d_j/Lad_{j-1}) \not\perp p$ and by
%our choice of $L$, $tp(d_j/Kad_{j-1}) \not\perp p$. 
Let $q=tp(d_j/Kad_{j-1})$.

Then we have $m(q)=m(H) \leq m(T(V)_a)=m(p)$, hence $m(p)=m(q)$.

%By \ref{prt14} there are tuples $a$ and $c$ of realizations of $p$ such that
%$c=Cb(a/Kc)$ and $r=tp(c/Ka)$ is $p$-semi-regular.
%Let $V$ be the locus of $a$ over $K$, by \ref{arc11} there is $m \in \na$ and a finite tuple $d=(d_1, \cdots, d_k)$ of
%$\mathcal{A}_mV_a$ such that $c \in dcl(K,a,d)$. As $c\not\in acl(K,a)$ there is a proper subtuple $e$ of $d$ and
%$i\in \{1, \cdots, k  \}$ such that $c \dnfo_{Ka} e$ and $c \dfo_{Kae} d_i$.
% Then $r$ is nonorthogonal to
%$tp(d_i/Kae)$.
%Let $b_1, \cdots,b_m=d_i$ be the elements given by \ref{arc12} (applied to $e_i$). Then there is $j$ such that
%$r \not\perp tp(b_j/Kaeb_{j-1})$. By \ref{arc12}, there is a type $q$ such that $r \not\perp q$ and
%$(x\in T(V)_a) \in r$. As $r$ is $p$-semi-regular $p \not\perp r$, thus
%$p \not\perp q$. Take $G=T(V)_a$. By \ref{arc41}, $m(G)=m$, then $m(q) \leq m$, and
%since $p$ is $(\si,D)$-type minimal, by \ref{arc8}, $m(q)=m$\\
$\Box$

\begin{lem}\label{arc14}
Let $p$ be a  regular $(\si,D)$-type minimal type. If there are a $(\si,D)$-vector group $G$
and a type $q$ that satisfy the conclusions of  \ref{arc13}, then there exists a $(\si,D)$-vector group
whose generic type is regular, $(\si,D)$-type minimal and non-orthogonal to $p$.
\end{lem}

{\it Proof}:\\

We order the triplets $ord(G)=\{m(G),dim_{\si,D}(G),\s(G)\}$ with the lexicographical
order. We proceed by induction on $ord(G)$.

{\bf Claim}:

We may assume that if $H$ is a proper quantifier-free connected,
 quantifier-free
definable subgroup of $G$,  then $m(H)<m(G)$.

{\bf Proof}:
Suppose that $m(H)=m(G)$. Let $\mu:G \to G/H$ be the quotient map. By \ref{arc5}, $ord(G)>ord(G/H)$.
If we replace $q$ by a nonforking extension of $q$ we may assume that $H$ is defined over the domain $A$ of $q$.
Let $a$ be a realization of $q$ with $tp(a/A) \not\perp p$. As
$q \not\perp p$, we have either $p \not\perp q_0=qftp(\mu(a)/A)$ or
$p \not\perp q'=qftp(a/A\mu(a))$.
If $p \not\perp q_0$  then $m(p) \leq m(q_0)$ by \ref{arc8},  and since $(x \in G/H) \in q_0$,
$m(q_0)\leq m(G/H)\leq m(G)=m(p)$. So $m(q_0)=m(p)$ and we apply induction hypothesis to $p,q_0$ and $G/H$.
If $p\not\perp q'$, let $b$ be a
realization of $qftp(a/A\mu(a))$ such that $b \dnfo_{A \mu(a)} a$.
Then  $a-b \in H$ and $p \not\perp q''=qftp(a-b/Ab)$ and the same
argument applies.\\

By \ref{arc10} and as $q$ is realized in $G$ and
$m(p)=m(q)=m(G)$, $q$ is a generic of $G$, and is regular and $(\si,D)$-type minimal.

%Write $q=qftp(a/A)$; let $q_0=qftp(\mu(a)/A)$ and $q'=qftp(a/A\mu(a))$.
%Let $b \in a+H $ such that $b \dnfo_{A\mu(a)} a$ and let $q''=qftp(a-b/Ab)$. Thus $q''$ is a translation of the
% nonforking extension of $q'$ to $Ab$. By transitivity we have that either $p \not\perp q_0$ or $p \not\perp q''$.
%If $p \not\perp q_0$  then $m(p) \leq m(q_0)$ by \ref{arc8},  and since $(x \in G/H) \in q_0$,
%$m(q_0)\leq m(G/H)\leq m(G)=m(p)$. So $m(q_0)=m(p)$ and we apply induction hypothesis to $p,q_0$ and $G/H$.The same
%argument applies for $p\not\perp q''$, and the claim is proved.\\
%
%now let $r$ be the generic type of $G$. Thus by the claim and \ref{arc10}, $r$ is regular.
%If we take nonforking extensions of $p$ and $q$ we may assume that $p \not\perp^a q$. Let $a$ and $b$ be two
%dependent realizations of $p$ and $q$ respectively. Since $r$ is the generic type of $G$ there are $c_1,c_2$
%such that $b \in dcl(c_1,c_2)$. Then $a \dfo c_1,c_2$. Suppose that $p \perp q$, this would imply that $a \dnfo c_1$,
%so $a \dfo_{c_1} c_2$, but $c_2$ realizes $r$, so $c_1 \dfo c_2$. Then the locus $V$ of $c_2$ over $c_1$ is a proper
%$(\si,D)$-subvariety of $G$, hence $m(V) < d$. Since $d=m(G)$, by \ref{arc8} we have $p \perp qftp(c_2/c_1)$ and this
%is absurd.\\
$\Box$

\begin{cor}\label{arc15}
Let $p$ be regular non locally modular  type. Then there is a $(\si,D)$-vector group $G$ whose generic type
is $(\si,D)$-type minimal and non-orthogonal to $p$.
\end{cor}

{\it Proof}:\\

%Let $q$ be a  regular type such that $p \not\perp q$, and such that it is of minimal $(\si,D)$-type.
%By \ref{orto}, $q$ is $(\si,D)$-type minimal.
By \ref{arc7} there is a regular type $q$ of minimal $(\si,D)$-type which is non-orthogonal to $p$.
By \ref{arc13}, $q$ satisfies the hypothesis of \ref{arc14}, then there is a $(\si,D)$-vector group $G$ whose generic type
$r$ is nonorthogonal to $q$; again by \ref{orto}, then there is such an $r$ which is non-orthogonal to $p$.\\
$\Box$

\begin{lem}\label{arc16}
Let $G$ be a $(\si,D)$-vector group and let $p$ be its generic type. If $p$ is regular there is a definable subgroup
of $\gr_a$ whose generic type is regular and non-orthogonal to $p$.
\end{lem}

{\it Proof}:\\

Suppose that $G <\gr_a^d$ for some $d \in \na$. One of the projections $\pi:G \to \gr_a$ must have an infinite image
in $\gr_a$.
Let $a$ realize $p$, then $\pi(a)$ realizes the generic type of $H=\pi(G)$; this type is $tp(\pi(a)/K))$
which is also regular. Hence $H$ satisfies the conclusion of the lemma.\\
$\Box$ 
  
\begin{teo}\label{arc17}
Let $p$ be a regular non locally modular type. Then there is a definable subgroup of the additive group
 whose generic type is regular and non-orthogonal to $p$.
\end{teo}

{\it Proof}:\\

By \ref{arc15} there is a $(\si,D)$-vector group $G$ whose generic type $q$ is regular and non-orthogonal to $p$,
 by \ref{arc16} there is a definable subgroup $H$ of the additive group  whose generic type $r$ is regular and
non-orthogonal to $q$. By transitivity $p \not\perp r$.\\
$\Box$

\begin{lem}\label{arc18}
Let $G$ be a definable subgroup of $\gr_a^n$. If $G$ has infinite dimension then $\s(G)\geq \omega$.
\end{lem}

{\it Proof}:\\

By \ref{arccor1}, $G$ is quantifier-free definable and is a 
$(Fix \si\cap \mathcal{C})$-vector space.
%If the dimension of $G$ as a $(Fix(\si)\cap \mathcal{C})$-vector space is $n$, then $G$
%is definably isomorphic to $(Fix(\si)\cap \mathcal{C})^n$, then it is finite-dimensional which is absurd. 
%Hence the dimension of $G$ as a $(Fix(\si)\cap \mathcal{C})$-vector space is infinite, and this implies that
%for  arbitrarely large $n$ $G$ contains a subgroup definable isomorphic to  $(Fix(\si)\cap \mathcal{C})^n$; as the 
%latter has $\s$-rank $n$, $\s(G)$ is infinite.\\
If $g_1,\ldots,g_n\in G$ are $(Fix \si\cap {\mathcal C})$-linearly independent, 
then the subgroup $H$ they generate is definable and has $\s$-rank $n$ (since it is definably isomorphic 
to $(Fix \si\cap {\mathcal C})^n$). Thus our hypothesis implies that $G$ contains elements of arbitrarily 
high finite $\s$-rank, 
and therefore that $\s(G)\geq \omega$. \\
%As $G$ is a definable subgroup of $\gr_a^n$ it is defined by a system of linear $(\si,D)$-equations, then for
%any $c \in Fix \si \cap {\mathcal C}$, and for any $g \in G$, $cg \in G$.
%Let $g \in  G$ be non-zero and
%let $H=\{cg: c \in Fix\si \cap {\mathcal C}\}$. 	
%Then $H$ is a definable subgroup of $G$ and $\s(H)=1$. 
%This shows that $G$ contains subgroups of arbitrarily high finite $\s$-rank. By Lascar's inequalities,
%$\s(G)$ is greater than or equal to its dimension as a $Fix\si \cap \mathcal{C}$-vector space.
%The iteration of this proceedure implies 
%Hence $\s(G) \geq \omega$.\\
$\Box$

\begin{teo}\label{arc19}
Let $p$ be a regular type of $\s$-rank 1. If $p$ is non locally modular then it is
non-orthogonal to $Fix \si \cap {\mathcal C}$.
\end{teo}

{\it Proof}:\\

By \ref{arc17} there is a definable subgroup $G$ of $\gr_a$ whose generic type $q$ is regular,
$(\si,D)$-type minimal and non-orthogonal to $p$.
$p \not \perp q$ implies that $\s(q)=\alpha+1$ for some $\alpha$.
Then, by 5.4.3 of \cite{wag}, $G$ contains a definable subgroup $N$ 
such that $\s(G/N)<\omega$, and by \ref{arc18} and  \ref{concomp}, $G$ must be 
finite-dimensional.
%Then, by \ref{arc18} and \ref{arc10}, $G$ is finite-dimensional.
 Thus, by \ref{J313}, $p \not\perp Fix \si \cap {\mathcal C}$.\\
%Let $q$ be the generic of a definable subgroup $G$ of $\gr_a$, and assume that $q$ is non-orthogonal to $p$, $\s(p)=1$.
 %Then $\s(q)=\alpha+1$ for some $\alpha$, then $G$ has a definable normal subgroup $N$
 %such that $\s(G/N)<\omega$, and the generic of $G/N$ is non-orthogonal to $p$. By ???, $G/N$ is a finite-dimensional
 %vector space over $Fix(\si)\cap {\cal C}$, and therefore $p$ is finite-dimensional. 
$\Box$

\cleardoublepage
\chapter{Definable Groups}
\label{chap:grupos}

This chapter is devoted to the study of definable groups in {\it DCFA}.
The fact that in difference-differential fields, having infinite $(\si,D)$-transcendence degree does not characterisize the
$(\si,D)$-generic type, represents a difficulty in the treatment of definable groups,
so we shall try different ways to describe certain kind of definable groups departing from properties of groups definable in
differential and difference fields.
In the first section %we give some results that will help us to reduce some questions to simpler definable groups.
%In section 2
 we follow the work of Kowalski and Pillay (\cite{kopi}) to show that a definable group 
is embedded in an algebraic group. %In section 2 we define the generalized logarithmic derivative, 
%and with its help we give some conditions satisfied by subgroups of the cartesian product of a $D$-group and 
%an additive group. 
Section 2 is devoted to the study $1$-basedness, stability and stable embeddability 
of commutative groups.

%\section{Preliminaries}

 %Then,  by the way independence in {\it DCFA} is
%defined,  $p_n(g) \dnfo_E p_n(h)$ for {\it DCF} if and only if $g \dnfo_E h$ for {\it DCFA};
%thus $ tp(g/E) $ is generic type of $H$ if and only if for every $n \in {\mathbb N} $ ,
% $tp(p_n(g)/E)   $ is generic type of $H^{(n)}  $ in {\it DCF}.

%In particular, if $h \in \tilde{H}    $ such that $g \dnfo_E h  $;
 %then, for every $n \in {\mathbb N}  $, $p_n(g) \dnfo_E p_n(h)  $ in {\it DCF}; then, for every $n \in {\mathbb N}  $,
%$p_n(g ) \cdot p_n(h) \dnfo_E p_n(h) $ in {\it DCF}. Thus $g \cdot h \dnfo_E h $.
%This implies that any generic type of $\tilde{H}$ is a generic type of $H$. \\

%Applying \ref{GIII5} we have:

%\begin{prop}\label{GLD3}
%$\s(\tilde{H})=\s(H)$ and $ [\tilde{H}:H] < \infty .$
%\end{prop}

%Similar reductions can be made using $\tau_m$ or $\Phi_m$.

\section{A Definable Group is Embedded in an Algebraic Group}

We introduce $*$-definable groups in stable theories.
Suppose that $T$ is a complete theory and $M$ a saturated model of $T$.
A $*$-tuple is a tuple $(a_i)_{i \in I}$, where $I$ is an index set of cardinality less than the cardinality of $M$,
and $a_i \in M^{eq}$ for all $i \in I$.
Let $A \subset M$. A $*$-definable set is a collection of $*$-tuples, indexed by the same set of parameters $I$,
which is the set of realizations of  a partial type $p(x_i)_{i \in I}$ over $A$.
A $*$-definable group is a group with $*$-definable domain and multiplication.

The following propositions are proved in \cite{kopi}.
\begin{prop}\label{GIII1}
Let $T$ be a stable theory; $M$ a saturated model of $T$.
Let $a,b,c,x,y,z$ be $*$-tuples of $M$ of length strictly less than the cardinal of $M$, such that:

\begin{enumerate}

\item $acl(M,a,b)=acl(M,a,c) = acl(M,b,c)       $

\item $acl(M,a,x)  =  acl(M,a,y)$ and $Cb(stp(x,y/M,a)) $ is interalgebraic with $a$ over $M $.

\item As in 2. with $b,z,y$ in place of $a,x,y $

\item As in 2. with $c,z,x$ in place of $a,x,y $

\item Other than $\{a,b,c  \},  \{a,x,y  \},  \{b,z,y  \},\{c,z,x  \} $, any 3-element subset of $\{a,b,c,x,y,z    \}  $ is independent over $\mathcal M$.
\end {enumerate}

Then there is a $*$-definable group $H$ defined over $M$ and $a',b',c'  \in H $ generic independent
 over $M$ such that $a$  is interalgebraic with $a'$ over $M$,   $b$  is interalgebraic with $b'$ over $M$
 and  $c$  is interalgebraic with $c'$ over $M$.

\end{prop}

\begin{prop}\label{GIII2}
Let $T$ be a simple theory; $M$ a saturated model of $T$. Let $G,H$ be type-definable groups,
 defined over $K \prec {M}$, and let $a,b,c \in G $  and $a',b',c' \in H  $ such that

\begin{enumerate}
\item $a,b$ are generic independent over $M$.
\item $a \cdot b = c  $ and $a' \cdot b' =c'  $.
\item $a$ is interalgebraic with $a'$ over $M$,  $b$ is interalgebraic with $b'$ over $M$ and
 $c$ is interalgebraic with $c'$ over $M$
\end{enumerate}
Then there is a type-definable over $M$ subgroup $G_1$ of bounded index in $G$,
and a type-definable over $M$ subgroup $H_1$ of $H$ and a  type-definable over $M$
isomorphism $f$ between $G_1/N_1  $ and $H_1 / N_2$ where $N_1  $ and $N_2$  are finite normal subgroups of $G_1 $ and $H_1$ respectively.

\end{prop}

\begin{rem}
If $T$ in \ref{GIII2} is supersimple and $G,H$ are definable,
then we can choose $G_1$ definable of finite index in $G$ and $f$ definable.
\end{rem}
The following result is proved in \cite{Hru}:
\begin{prop}\label{GIII3}
Let $G$ be a $*$-definable group in a stable structure.
Then there is a projective system of definable groups with inverse limit $G'$,
and a $*$-definable isomorphism between $G$ and $G'$.
\end{prop}

\begin{teo}\label{GIII6}
Let $({\mathcal U},\si,D)$ be a model of {\it DCFA}, $K \prec {\mathcal U} $ and $G$ a $K$-definable group.
Then there is an algebraic group $H$, a definable subgroup $G_1$  of $G$ of finite index,
and a definable isomorphism between $G_1/N_1$ and $H_1 / N_2$, where  $H_1$ is a definable subgroup of $H(\mathcal{U})$,
 $N_1$ is a finite normal subgroup of $G_1$, and  $N_2$ is a  finite normal subgroup of $H_1$.
\end{teo}

{\it Proof:}\\

Let $a,b,y   $ be generic independent elements of $G$ over $K$. Let $x = a \cdot y,  z = b^{-1} \cdot y,  c = a \cdot b  $, so $x = c \cdot z  $.
Let $\bar{a}=(D^i\si^j (a): i  \in {\mathbb N},   j \in {\mathbb Z})  $, and similarly for $\bar{b}, \bar{c}, \bar{x},\bar{y},\bar{z} $.
Then by \ref{DCFA421}, working in {\it ACF},
 $\bar{a},\bar{b},\bar{c},\bar{x},\bar{y},\bar{z} $
 satisfy the conditions of \ref{GIII1}.
Thus there is a $*$-definable group $H$ over $K$, and generic $K$-independent elements $a^*,b^*, c^* \in H    $ such that
 $\bar{a}  $ is interalgebraic with  $a^*  $ over $K$, $\bar{b}  $ is interalgebraic with  $b^*  $ over $K$,
 $\bar{c}  $ is interalgebraic with  $c^*  $ over $K$  and
 $ c^*=  a^* \cdot  b^*   $ (the interalgebraicity is in the sense of {\it ACF}).

Since  {\it ACF} is $\omega$-stable, by \ref{GIII3}, $H$ is the inverse limit of $H_i, i \in \omega $,
 where the $H_i$ are algebraic groups.

Let $\pi_i : H \longrightarrow H_i  $ be the $i$-th canonical epimorphism.
Let $a_i = \pi_i (a^*) $, $b_i = \pi_i (b^*) $ and $c_i = \pi_i (c^*) $. Then $a^*   $ is
interalgebraic with $(a_i)_{i \in \omega } $ over $K$ ,
 $b^*   $ is interalgebraic with $(b_i)_{i \in \omega}  $ over $K$ and
$c^*   $ is interalgebraic with $(c_i)_{i \in \omega}  $ over $K$, all interalgebraicities in the
sense of {\it ACF}.

Since for $i < j$, $a_i \in K(a_j)$  , $b_i \in K(b_j)$  and  $c_i \in K(c_j)$,
there is $i \in \omega$ such that $a$ is interalgebraic with $a_i$ over $K$,   $b$ is interalgebraic
 with $b_i$ over $K$ and  $c$ is interalgebraic with $c_i$ over $K$  in the sense of {\it DCFA}.
So we can apply  \ref{GIII2} to $a,b,c \in G$ and $a_i,b_i,c_i \in H_i$.\\
$\Box$

\section{Abelian Groups}

In this section, we study abelian groups defined over some 
subset $K=acl(K)$ of a model $({\mathcal U},\si,D)$ of {\it DCFA}. 
We investigate whether they 
are $1$-based, and whether they are stable stably embedded (i.e., stable 
with the structure induced by $\mathcal U$). By \ref{GIII6} and \ref{wgr}, we 
may reduce to the case when the group $H$ is a quantifier-free 
definable subgroup of some commutative algebraic group $G$, and $G$ has no proper (infinite) algebraic subgroup, i.e. $G$ is either $\gr_a$, 
$\gr_m$, or a simple Abelian variety $A$. 

%In this section we study abelian groups defined over difference-differential fields.\\ 
%Theorem \ref{wgr} allows us to reduce our study to subgroups of $\gr_a$, $\gr_m$ or 
%a simple abelian variety.
{\bf From now on we suppose all the groups are quantifier-free definable}.\\

{\bf The additive group}
\begin{prop}\label{adgr}
No infinite  definable subgroup of  $\gr_a^n(\mathcal{U})$ is $1$-based.
\end{prop}

{\it Proof}:\\

Let $H<\gr_a^n$ be a definable infinite group. By \ref{arccor1},
$H$ is quantifier-free definable and contains a definable subgroup $H_0$ which is definably isomorphic to 
$Fix \si \cap {\mathcal C}$. Hence $H$ is not $1$-based.\\
%It is clear that $H$ is defined by linear difference-differential equations, and
%therefore it is a $(Fix \si \cap {\mathcal C})$-vector space.\\
$\Box$\\

{\bf The multiplicative group}\\
The logarithmic derivative $lD:\gr_m \to \gr_a$, $x\mapsto Dx/x$ is a group epimorphism with
$Ker(lD)=\gr_m(\mathcal{C})$. 

Given a polynomial $P(T)= \sum_{i=0}^na_iT^i\in {\mathbb Z}[T]$, we denote by $P(\si)$ the homomorphism defined by
$x \mapsto \sum_{i=0}^na_i\si^i(x)$.

\begin{prop}
Let $H$ be a quantifier-free ${\mathcal L}_{\si,D}$-definable subgroup of $\gr_m$. 
If $lD(H) \neq 0$ then $H$ is not $1$-based. If $lD(H)=0$ then there is a polynomial
$P(T)$ such that $H=Ker(P(\si))$. Then we have that $H$ is $1$-based if and only if $P(T)$ 
is relatively prime to all cyclotomic 
 polynomials $T^m-1$ for all $m\in \na$ 
\end{prop}

{\it Proof}: \\

By \ref{adgr}, if $lD(H) \neq0$ then $H$ is not $1$-based.
If $lD(H)=0$, as $Ker(lD)=\gr_m(\mathcal{C})$, $H$ is ${\mathcal L}_{\si}$-definable in $\mathcal{C}$. 
Hence there is a polynomial 
$P(T)=\sum_{i=0}^na_iT^i \in \mathbb{Z}[T]$ such that $H$ is defined by $\Pi_{i=0}^n\si^i(X^{a_i})=1$.
In {\it ACFA}, $H$ is  1-based, stable, stably embedded if and only if $P(T)$ is relatively prime to all cyclotomic 
 polynomials $T^m-1$ 
for $m \geq 1$ (see \cite{HMM}). By \ref{st1} the same holds for {\it DCFA}.\\
$\Box$

{\bf Abelian varieties}\\

%This relation has lead to model theoretic proofs of Manin-Mumfurd and
%Mordell Lang conjectures both by E. Hrushovski (\cite{HMM} and \cite{HML}). Recently Pillay found
%more direct proofs to this conjectures, for the case of characteristic zero (\cite{PMM} and \cite{PML}).
First we mention some facts about Abelian varieties in difference and differential fields.
For a detailed exposition on Abelian varieties the reader may consult \cite{langabelian}.
                             
\begin{defi}\label{ab1}
An Abelian variety is a connected algebraic group $A$ which is complete, that is, for any variety $V$ the
projection $\pi:A \times V \to V$ is a closed map.
\end{defi}

As a consequence of the definition we have that an Abelian variety is commutative.

Let $B$ be an algebraic subgroup of an Abelian variety $A$. Then $A/B$ is an Abelian variety. If in addition $B$
is connected $B$ is an Abelian variety.
An Abelian variety is called simple if it has no infinite proper Abelian subvarieties.
Let $A$ and $B$ be two Abelian varieties. Let $f:A \to B$ be a homomorphism. We say that $f$ is an isogeny if
$f$ is surjective and $Ker(f)$ is finite. We say that $A$ and $B$ are isogenous if  there are isogenies
$f:A \to B$ and $g:B \to A$.

\begin{prop}\label{ab5}\textup{({\it ACF})}
There is no nontrivial algebraic homomorphism from a vector group into an Abelian variety.
\end{prop}

Now we mention some properties concerning 1-basedness of Abelian varieties in difference and differential fields.
Consider a saturated model $({\mathcal U},\si)$ of {\it ACFA}.
% Let $P(T)=\sum_{i=0}^na_iT^i \in {\mathbb Z}[T]$.
%Define $Ker(P(\si))(A)=\{a \in A:\sum_{i=0}^na_i\si^i(a)=0\}$.

In \cite{HMM}, Hrushovski gives a full description of definable subgroups 
of $A({\mathcal U})$ when $A$ is a simple Abelian variety defined over ${\mathcal U}$.
 When $A$ is defined over $Fix\si$, this description is particularly simple, at least up to commensurability. 
Let $R=End(A)$ (the ring of algebraic endomorphisms of $A$). If $P(T)=\sum_{i=0}^n e_iT^i\in R[T]$,
 define $Ker(P(\si))=\{a\in A({\mathcal U})\mid \sum_{i=0}^n e_i(\si^i(a))=0\}$.

\begin{prop} \textup{({\it ACFA}, \cite{HMM})}\label{ab2} 
Let $A$ be a simple Abelian variety defined over $\mathcal U$, and let $B$ be a definable subgroup of 
$A({\mathcal U})$ of finite $\s$-rank. 
\begin{enumerate}
\item If $A$ is not isomorphic to an Abelian variety defined over $(Fix\si)^{alg}$, 
then $B$ is $1$-based and stable stably embedded. 
\item Assume that $A$ is defined over $Fix\si$. 
Then there is $P(T)\in R[T]$ such that $B\cap Ker(P(\si))$ has finite index in $B$ and in $Ker(P(\si))$. 
Then $B$ is $1$-based if and only if the polynomial $P(T)$ is relatively prime to all cyclotomic polynomials 
$T^m-1$, $m\in\na$. If $B$ is $1$-based, then it is also stable stably embedded.
\end{enumerate} 
\end{prop}

%\begin{prop}\label{ab2}(ACFA, see \cite{HMM})\
%If $A$ is defined over $Fix \si$, then
%$Ker(P(\si))(A)$ is stable,stably embedded and 1-based if and only if $P(T)$ is relatively prime to all cyclotomic
%polynomials $T^n-1$, $n \in \na$.
%\end{prop}

%We make an observation. When we work in characteristic 0, in {\it ACFA} if $tp(a/E)$ is hereditary
%orthogonal to $Fix\si$ then it is stable and stably embedded. This comes from the fact that,
%given $a,F$ and $E' \supset E$ such that $a \dnfo_{E'} F $, $acl_{\si}(E'a)F$ has no finite $\si$-stable extension,
% and this
%implies that $tp_{ACFA}(a/E')$ is stationary.
%From this we have that if a type is 1-based in $\mathcal{C}$ then it is 1-based in $\mathcal{U}$
%(see Chapter 2, section 5).
%\begin{defi}\label{ab3}
%Let $G$ be an algebraic group. By an extension of $G$ by a vector group we mean %an algebraic group $H$
%together with an epimorphism $p: H \to G$ such that $Ker(p)$ is a vector group.
%\end{defi}

%\begin{prop}\label{ab4}
%Let $A$ be an Abelian variety. Then there is a universal extension of $A$ by a vector group. That is,
%there is an extension of $A$ by a vector group $({\hat A},p)$, such that if $(B,q)$ is an extension
%of $A$ by a vector group, then there is $j:{\hat A} \to B$ such that $q = p \circ j$.
%\end{prop}
%As a consequence of \ref{ab4} we have that there is an algebraic section $s:{\hat A} \to \tau({\hat A})$.

We work now in a saturated model $({\mathcal U},D)$ of {\it DCF}.

\begin{prop}\label{ab6}
Let $A$ be an Abelian variety. Then there is a ${\mathcal L}_D$-definable (canonical)homomorphism
$\mu: A \to \gr_a^n$, for $n=dim(A)$, such that $Ker(\mu)$ has finite Morley rank.
\end{prop}
the kernel of this canonical homomorphsim, $Ker(\mu)$, is known as the Manin kernel of $A$, we denote it by $A^{\sharp}$. 

\begin{prop}\textup{(Properties of the Manin Kernel, see \cite{manin} for the proofs)}\label{propmanin}\\
Let $A$ and $B$ be Abelian varieties. Then
\begin{enumerate}
\item $A^{\sharp}$ is the Kolchin closure of the torsion subgroup $Tor(A)$ of $A$.
\item $(A\times B)^{\sharp}=A^{\sharp}\times B^{\sharp}$, and if $B<A$ then $B\cap A^\#=B^\#$.
%hence if $B<A$, then $B\cap A^{\sharp}=B^{\sharp}$ 
%(since $A$ is isogenous to the direct product of $B$ with an Abelian variety). 
\item A differential isogeny  between
$A^{\sharp}$ and $B^{\sharp}$ is the restriction of an algebraic isogeny from $A$ to $B$.
\end{enumerate}
\end{prop}

\begin{defi}
We say that an Abelian variety descends to the constants if it is isomorphic to an Abelian variety defined over the constants.
\end{defi}

\begin{prop}\label{ab8}\textup{({\it DCF}, see \cite{manin})}
Let $A$ be a simple Abelian variety. If $A$ is defined over ${\mathcal C}$, then $A^{\sharp}=A({\mathcal C})$.
If $A$ does not descend to the constants, then $A^{\sharp}$ is strongly minimal and 1-based.
\end{prop}

We now return to {\it DCFA} and fix a saturated model $({\mathcal U,\si,D})$ of {\it DCFA} 
and a simple Abelian variety $A$ defined over $K=acl(K) \subset \mathcal{U}$.

Let $H$ be an ${\mathcal L}_{\si,D}$-definable connected subgroup of $A$ defined over 
the difference-differential field $K$.
 Since $H$ is $1$-based if and only if  $\tilde{H}$
is $1$-based,  we can suppose that $H$ is quantifier-free definable and quantifier-free connected.
%Then there is $k \in \na$ and an ${\mathcal L}_D$-definable subgroup $S$ of
%$A \times \cdots \times A^{{\si}^k}$ such that
%$H=\{x \in A:(x,\si(x), \cdots, \si^k(x)) \in S  \}$.

%If $G=\gr_a^n$, we know that there are no 1-based subgroups of an additive group.
%If $G=\gr_m^n$, then the $D$-definable subgroups of $G$ are defined by linear
%equations in $x,Dx, \cdots,D^mx$ for some $m$, then it is non 1-based because
%it is a $(Fix \si \cap \mathcal{C})$-vector space.\\
Let $\mu:A \to \gr_a^d$ as in \ref{ab6}. If $H \not\subset Ker\mu$ then by \ref{adgr}
$H$ is not $1$-based.

Assume that $H \subset A^{\sharp}$. 
%Then $S \subset A^{\sharp} \times (A^{\sharp})^{\si} \times \cdots \times (A^{\sharp})^{\si^k}$.
We first show a very useful lemma.

\begin{lem}\label{ab20}
Let $H$ be a quantifier-free definable subgroup of $A^{\sharp}$ 
which is quantifier-free connected. Then $H=H'\cap A^{\sharp}$ for some 
quantifier-free ${\cal L}_\si$-definable subgroup $H'$ of $A$.
%Let $H$ be a quantifier-free ${\mathcal L}_{\si,D}$-definable subgroup of 
%$A^{\sharp}$. Then there is $k$ and a differential 
%subgroup $S$ of $A^{\sharp}\times\cdots\times (A^{\sharp})^{\si^k}$ such that 
%$H=\{a\in A :(a,\si(a),\ldots,\si^k(a))\in S\}$. We 
%know that $S=U^{\sharp}$ for some algebraic subgroup $U$ of $A\times 
%\cdots\times A^{\si^k}$. By \ref{propmanin}, the results of \cite{HMM} give a full 
%description of the quantifier-free ${\mathcal L}_{\si,D}$-definable 
%subgroups of $A^{\sharp}$. 
\end{lem}

{\it Proof}:\\

%Let $H$ be a quantifier-free definable  subgroup of $A^{\sharp}$ which is connected for the $(\si,D)$-topology. 
%Then there is an integer
%$k$ and  a connected algebraic subgroup $S$ of $A\times A^\si\times 
%A^{\si^k}$ such that $H=\{a\in A^{\sharp} : (a,\si(a),\cdots,\si^k(a)\in S\}$.
Our hypotheses imply that there is an integer $k$ and a differential 
 subgroup $S$ of $A\times A^\si\times \cdots \times A^{\si^k}$ such that 
$H=\{a\in A: (a,\si(a),\cdots,\si^k(a))\in S\}$. By \ref{propmanin}.2, replacing $S$ 
by its Zariski closure $\bar S$ we get $H=\{a\in A^{\sharp}: 
(a,\si(a),\cdots,\si^k(a))\in \bar S\}$. 
Thus $H=H'\cap A^{\sharp}$, with $H'=\{a\in A : (a,\si(a),\cdots,\si^k(a)\in
\bar{S}\}$.\\
$\Box$

%{\it Proof}:\\
%Otherwise, $A^{\sharp}$ contains an infinite definable proper subgroup; as in \cite{HMM}, we get a $D$-definable 
%isogeny $f:A^{\sharp} \to (A^{\sharp})^{\si^k}$. By \ref{ab6.1} this gives us an algebraic isogeny $f':A \to A ^{\si^k}$.
%
%$\Box$

%\begin{rem}
%With the help of \ref{ab6.1} and \ref{J313}, we can prove that if $A$ is not isomorphic to an Abelian variety defined over 
%$Fix \si^k$ then every finite-dimensional subgroup of $A^{\sharp}$ is stable (the proof is the same as in \cite{HMM}).
%\end{rem}      
%\begin{prop}\label{ab7}
%Let $A$ be a simple Abelian variety. Then $A^{\sharp}$ has no infinite $D$-definable subgroups.
%\end{prop}

{\bf Case 1}: $A$ is isomorphic to a simple Abelian variety $A'$ defined over $\mathcal{C}$.\\
We can suppose that $A$ is defined over $\mathcal{C}$. Then, by \ref{ab8},
$A^{\sharp}=A(\mathcal{C})$.
% Hence $S \subset A(\mathcal{C}) \times A^{\si}(\mathcal{C}) \times \cdots \times A^{\si^k}(\mathcal{C})$. 
Then $H$ is $1$-based for {\it DCFA} if and only if it is $1$-based for $ACFA$, by \ref{st1},
and in that case it will also be stable stably embedded (by 
\ref{st3})

 If $H=A(\mathcal{C})$ then we know that $H$ is not 1-based in {\it ACFA}.
%and
% this implies that $H$ is 1-based for {\it DCFA}: From the way independence is defined we have that if $tp_{DCF}(a)$ then $tp(a)$ is 1-based, similarly, if
 %$tp(a)$ is not 1-based then one of the types $tp_{DCF}(D^ka/a \cdots D^{k-1}a)$
 %is not 1-based.

If $H$ is a proper subgroup of $A(\mathcal{C})$, 
\ref{ab2} gives a precise 
description of 
that case.
%, hence, by \ref{ab2}, $H$ is not 1-based if and
%only if some type in the semi-minimal analysis of $S$ is non orthogonal to $Fix \si \cap {\mathcal C}$ (see \cite{HMM}),
%and again the results of that paper give a precise description of this case.
\\\\
{\bf Case 2}: $A$ does not descend to $\mathcal{C}$.

Then, by \cite{manin}, $A^{\sharp}$ is strongly minimal and 1-based for {\it DCF}. By \ref{dcfdcfa} it is 1-based
for {\it DCFA}.

Let us first note an immediate consequence of \ref{ab20} : 

\begin{cor}\label{ab9}
If for all $k \in \na$, $A$ and $A^{\si^k}$ are not isogenous,
then $\s(A^{\sharp})=1$.
\end{cor}
%
%
%If $H \neq A^{\sharp}$ then $dim(H) < \infty$. Then if $A$ is not isomorphic to an Abelian veriety defined over
%$Fix \si^k$ then $H$ is stable and stably embedded.
%Then $H$ is not 1-based if and
%only if some type in the semi-minimal analysis is nonorthogonal to $Fix\si \cap \mathcal{C}$.
%If $H=A^{\sharp}$, as $A^{\sharp}$ is 1-based for {\it DCF} then it is 1-based for {\it DCFA}.

%\begin{teo}
%$\s(A^{\sharp})=1$.
%\end{teo}
%
%{\it Proof}:\\
%
%Suppose that there is an isogeny between $A$ and $A^{\si^k}$ for some $k$, then as in \cite{HMM},
%there is an $D$-definable isogeny
%$f:A^{\sharp} \to (A^{\sharp})^{\si^k}$. Without lose of generality we may
%suppose that $k=1$.
%
%So $f$ gives us a $D$-definable subgroup $H$ of
%$A^{\sharp} \times (A^{\sharp})^{\si}$. Then we have that there is $m \in \na$ and an algebraic
%subgroup $H'$ of $\tau_m(A) \times \tau_m(A^{\si})$ projecting onto $A$ and onto $A^{\si}$, such that
%$H=\{x \in A: (x,Dx,\cdots,D^m x) \in H'\}$. 

We will now investigate stability and stable 
embeddability of $H$. By $1$-basedness and quantifier-free 
$\omega$-stability, we know that if $X\subset A^{\sharp}$ is quantifier-free 
definable, then $X$ is a Boolean combination of cosets of 
quantifier-free definable subgroups of $A^{\sharp}$.  

Assume first that $H\neq A^{\sharp}$, and let $a$ be a generic of $H$ over $K$. 
Then $H$ is finite-dimensional, and therefore $\s(H)<\omega$.
As $H$ is $1$-based, there is an increasing sequence 
of subgroups $H_i$ of $H$ with  $\s(H_{i+1}/H_i)=1$.

By \ref{ab20}, we may assume that $H_i=U_i\cap A^{\sharp}$ for some quantifier-free 
${\mathcal L}_\si$-definable subgroups $U_i$ of $A$. Note that \ref{ab20} also 
implies that each quotient $U_{i+1}/U_i$ is $c$-minimal (i.e., all
quantifier-free definable ${\mathcal L}_\si$-definable subgroups are 
either 
finite or of finite index). Furthermore, by elimination of imaginaries 
in 
{\it ACFA}, $acl_{\si}(Ka)$ contains tuples $a_i$ coding the cosets $a+U_i$. 
Hence $tp(a/K)$ satisfies the conditions of \ref{st5} and we obtain that 
if $tp_{ACFA}(a/K)$ is stable stably embedded then so is $tp(a/K)$. 

For the other direction, observe that if $tp_{ACFA}(a/K)$ is not 
stable stably embedded, then for some $i$, the generic {\it ACFA}-type of 
$U_{i+1}/U_i$ is non-orthogonal to $\si(x)=x$, and there is a (${\mathcal 
L}_\si$)-definable morphism $\psi$  with finite kernel $U_{i+1}/U_i \to 
 B(Fix \si^k)$ 
for some $k$ and Abelian variety $B$ (see \cite{HMM}). 
But, returning to {\it DCFA}, no non-algebraic type realized in $Fix \si^k$ can 
be stable stably embedded, since for instance the formula 
$\varphi(x,y)=\exists z\ z^2=x+y \ \land \ \si(z)=z$ is not definable (\ref{st1},3). 
This proves the other implication.

\smallskip

%Assume first that $H\neq A^{\sharp}$. Then $H$ is finite-dimensional, and therefore has finite $\s$-rank. 
%By the discussion above,  we may assume that $\s(H)=1$. Let $a,m$ be as above, and consider 
%$tp_{ACFA}(a,Da,\cdots,D^ma/A)$. By \ref{st2}, $tp(a/K)$ will be stable and 
%stably embedded if and only if it is stationary. Assume that 
%$tp_{ACFA}(a,Da,\cdots,D^ma)$ is not stationary. By 4.11 of \cite{salinas}, this 
%implies that there is an algebraically closed  difference field $L$ 
%containing $K$ and linearly disjoint from $K(a,Da,\cdots,D^ma)_\si$ 
%over $K$, and an element $b\in Fix\si$ such that $b\in 
%acl_{\si}(L,a,\cdots,D^ma)$. Looking at the coeeficients of the minimal 
%polynomial of $b$ over $L(a,\cdots,D^ma)_{\si}$, we may assume that 
%$b\in L(a,\cdots,D^ma)_{\si}$. Hence, we may also assume that 
%$L=acl(L)$. Since $\s(a/L)=1$, this implies that $a\in 
%acl(L,b)=L(b)_D^{alg}$. Since all $D^ib\in Fix\si$, we get 
%$tp_{ACFA}(a/L)\not\perp \si(x)=x$. 
Thus we have shown: 

\smallskip 
If $H$ is finite dimensional, then $tp(a/K)$ is stable stably embedded 
if and only if $tp_{ACFA}(a/K)$ is stable stably embedded. 

%The results in \cite{HMM} then give a complete description of that case. 
Using \ref{ab20}, \ref{ab2} gives us a full description of that 
case.

In particular, we then have that if $H$ is not stable 
stably embedded, then $A$ is isomorphic to an Abelian 
variety defined over $Fix\si^k$ for some $k$. 

\smallskip 
Let us now assume that $H=A^{\sharp}$. Let $a$ be a generic of $H$ over $K$. Then 
$tp_{ACFA}(a,\cdots,D^ma/K)$ is the generic type of an algebraic 
variety $V$, and is therefore stationary (by 2.11 of \cite{salinas}). Thus, using the finite dimensional case, if 
$A$ is not isomorphic to an Abelian variety defined over 
$(Fix\si)^{alg}$, then $H$ is stable stably embedded. If $A$ is 
isomorphic to a variety $B$ defined over $Fix \si^k$, via an 
isomorphism $\psi$, then the subgroup $\psi^{-1}(Ker (\si^k -1))\cap A^{\sharp}$ is 
not stable stably embedded. 
%Indeed, we may assume that  $k$ is large 
%enough so that $B[2]\subset B(Fix(\si^k))$. Let $b\in B$, $\si^k(b)=b$, $b$ 
%differentially transcendental over $K$. Then 
%$tp(\psi(a)/K)$ has several distinct non-forking extensions to $acl(Kb)$: 
%those containing 
%the formula $\exists y \in B\, 2y=x+b \land \si^k(y)=y$, and those 
%containing its negation. 

We summarize the results obtained:

\begin{teo} \label{absum} 
Let $A$ be a simple Abelian variety, and let $H$ be a 
quantifier-free definable subgroup of $A({\mathcal U})$ defined over $K=acl(K)$. If 
$H\not\subset A^{\sharp}({\mathcal U})$, then $H$ is not 1-based. Assume 
now that $H\subset A^{\sharp}({\mathcal U})$, and let $a$ be a generic of 
$H$ over $K$. Then 
\begin{enumerate}

\item If $A$ is defined over the field  $\mathcal C$ of constants, then 
$H$ is $1$-based if and only if it is stable stably embedded, if and 
only if $tp_{ACFA}(a/K)$ is hereditarily orthogonal to 
$(\si(x)=x)$. The results in \cite{HMM} yield a complete description of 
the subgroups $H$ which are not $1$-based. 

\item If $A$ does not descend to the field  $\mathcal C$ of constants, 
then $H$ is $1$-based. Moreover 

\begin{enumerate}
\item If $A$ is not isomorphic to an Abelian variety defined over 
$Fix \si^k$ for some $k$, then $H$ is stable stably embedded. 

\item Assume that $A$ is defined over $Fix(\si)$. Then $H$ is stable 
stably embedded if and only $tp_{ACFA}(a/K)$ is stable stably 
embedded. Again, the results in  \cite{HMM} give a full description of 
this case. 
\end{enumerate}
\end{enumerate}
\end{teo}

\clearpage

\def\rasp{\leavevmode\raise.45ex\hbox{$\rhook$}} \def\cprime{$'$}

\end{document}